\documentclass{amsart}
\usepackage{amssymb}
\usepackage{verbatim,tikz}
\usepackage{mathrsfs, amsthm}
\usepackage{dsfont, color}
\usepackage[all]{xy}

\usepackage[colorlinks, citecolor=blue, linkcolor=red]{hyperref}
\numberwithin{equation}{section}

\newcommand\eT{{}^{e}T}
\newcommand\eW[1]{{}^{e}\wedge^{#1}M}
\newcommand\eS{{}^{e}\Sym^{2}}

\newcommand\mrurl[1]{MR. \url{http://www.ams.org/mathscinet-getitem?mr=#1}}
\newcommand\aurl[1]{\url{file:///data/ms/archive/#1}}
\newcommand\durl[1]{doi. \url{http:///dx.doi.org/#1}}
\newcommand\purl[1]{\url{file:///data/ms/print/#1}}
\newcommand\surl[1]{\url{file:///data/ms/screen/#1}}
\newcommand\zurl[1]{Zbl. \url{http://www.emis.de/zmath-item?#1}}
\newcommand\Aurl[1]{Arxiv. \url{http://arxiv.org/abs/#1}}

\newcommand\boxb[1]{\square_b}

\newcommand\ff{\operatorname{ff}}

\numberwithin{equation}{section}

\newcommand\paperbody%
        {}

\newtheorem{theorem}{Theorem}
\newtheorem*{theorem*}{Theorem}
\newtheorem{lemma}{Lemma}
\newtheorem{proposition}[lemma]{Proposition}
\newtheorem{corollary}[lemma]{Corollary}
\newtheorem{definition}[lemma]{Definition}

\theoremstyle{remark}
\newtheorem{remark}[lemma]{Remark}

\numberwithin{lemma}{section}

\newcommand\cV{\mathcal{V}_{\operatorname{c}}}

\newcommand\cFTs{{}^{\Phi}\overline{T}\kern-1pt{}^*}

\newcommand\Nul{\operatorname{Nul}}

\newcommand\Tr{\operatorname{Tr}}

\newcommand\Vol{\operatorname{Vol}}

\newcommand\SO{\operatorname{SO}}

\newcommand\cC{\mathcal C}

\newcommand\cK{\mathcal{K}}

\renewcommand\cV{\mathcal{V}}

\newcommand\BB{\mathbb B}
\newcommand\CC{\mathbb C}

\newcommand\HH{\mathbb H}

\newcommand\NN{\mathbb N}

\newcommand\RR{\mathbb R}
\renewcommand\SS{\mathbb S}
\newcommand\ZZ{\mathbb Z}

\newcommand\bbB{\mathbb B}

\newcommand\bbH{\mathbb H}

\newcommand\bbR{\mathbb R}
\newcommand\bbS{\mathbb S}

\newcommand\CI{{\mathcal{C}}^{\infty}}

\newcommand\Diff[1]{\operatorname{Diff}^{#1}}

\newcommand\cFNs{{}^{\Phi}\overline N\kern-1pt{}^*}

\newcommand\tr{\operatorname{tr}}

\newcommand\Hom{\operatorname{Hom}}

\newcommand\Id{\operatorname{Id}}
\newcommand\Lap{\varDelta}

\newcommand\Sym{\operatorname{Sym}}

\newcommand\dCI{\dot{\mathcal{C}}^{\infty}}

\newcommand\Ric{\operatorname{Ric}}

\newcommand\cl{\operatorname{cl}}

\renewcommand\Re{\operatorname{Re}}
\renewcommand\Im{\operatorname{Im}}

\begin{document}

\begin{abstract}
We study the eleven dimensional supergravity equations which describe a low energy approximation to string theories and are related to M-theory under the AdS/CFT correspondence. These equations take the form of a non-linear differential system, on
$\bbB^7\times\bbS^4$ with the characteristic degeneracy at the boundary of an edge system, associated to the fibration with fiber $\bbS^4.$ We compute the indicial roots of the linearized system from the Hodge decomposition of the
4-sphere following the work of Kantor, then using the edge calculus and scattering theory we prove that the moduli space of solutions, near the Freund--Rubin states, is parametrized by three pairs of data on the bounding 6-sphere.
\end{abstract}

\title{The eleven dimensional supergravity equations on edge manifolds}
\author{Xuwen Zhu}
\address{Department of Mathematics, Stanford University}
\email{xuwenzhu@stanford.edu}
\date{\today}

\maketitle
\section{Introduction}

Supergravity is a theory of local supersymmetry, which arises in the representations of super Lie algebras~\cite{van1981supergravity}. A supergravity system is a low energy approximation to string theories~\cite{Nastase:2011aa}, and can be viewed as a generalization of Einstein's equation. Nahm~\cite{nahm1978supersymmetries} showed that the dimension of the system is at most eleven in order for the system to be physical, and in this dimension if the system exists it is unique. The existence of such systems was later shown~\cite{MR0441129} by constructing a specific solution. More special solutions were constructed by physicists later~\cite{castellani1984bosonic,  van1985complete, bergshoeff1987supermembranes}. Witten~\cite{Witten:1996md} showed that under the AdS/CFT correspondence M-theory is related to 11-dimensional supergravity, and there are more recent results~\cite{blau2002penrose}. Under dimensional reduction the fields break into many subfields and there are many such lower dimensional systems~\cite{nahm1978supersymmetries}. The full eleven dimensional case,  with only two fields, is in many ways the simplest to consider.

We are specifically interested in the bosonic sectors in the supergravity theory, which is a system of equations on the 11-dimensional product manifold $M=\BB^7 \times \bbS^4$, the product of a 7-dimensional ball and a 4-dimensional sphere. The fields are a metric, $g$, and a 4-form, $F$. Derived as the variational equations from a Lagrangian, the supergravity equations are
\begin{equation}\label{eq:supergravity}
\begin{array}{l}
R_{\alpha\beta}=\frac{1}{12}(F_{\alpha \gamma_1\gamma_2 \gamma_3}F_{\beta}^{\gamma_1 \gamma_2 \gamma_3}-\frac{1}{12}F_{\gamma_1 \gamma_2 \gamma_3 \gamma_4}F^{\gamma_1 \gamma_2 \gamma_3 \gamma_4}g_{\alpha \beta})\\
d*F=-\frac{1}{2}F \wedge F\\
dF=0
\end{array}
\end{equation}
where $R_{\alpha\beta}$ is the Ricci tensor, and we are using the Einstein summation notation here. For simplicity later, we also use the following notation 
$$
F\circ F:=\frac{1}{12}(F_{\alpha \gamma_1\gamma_2 \gamma_3}F_{\beta}^{\gamma_1 \gamma_2 \gamma_3}-\frac{1}{12}F_{\gamma_1 \gamma_2 \gamma_3 \gamma_4}F^{\gamma_1 \gamma_2 \gamma_3 \gamma_4}g_{\alpha \beta}).
$$

The nonlinear supergravity operator has an edge structure in the sense of Mazzeo \cite{MR1133743}, which is a natural generalization in this context of the product of a conformally compact manifold and a compact manifold. We consider solutions as sections of the edge bundles, which are rescalings of the usual form bundles. The Fredholm property of certain elliptic edge operators is related to the invertibility of the corresponding normal operator $N(L)$, which is the lift of the operator to the front face of the double stretched space $X_e^2$ which appears in the resolution of these operators.  The invertibility of the normal operator is in turn related to its action on appropriate polyhomogeneous functions at the left boundary of $X_e^2$, the form of the expansion of the solution is determined by the indicial operator $I_{s}(L)$. The inverse of the indicial operator $I_{s}(L)^{-1}$ exists and is meromorphic on the complement of a discrete set $\{s\in \operatorname{spec}_bL\}$, the set of indicial roots of $L$ which are the exponents in the expansion. In this way the indicial operator as a model on the boundary determines the form of the leading order expansion of the solution. 

One solution for this system is given by the metric which is the product of the round sphere with a Poincar\'{e}--Einstein metric on $\bbB^7$ with a volume form on the 4-sphere as the 4-form, in particular the Freund--Rubin solution~\cite{freund1980dynamics} is contained in this class. Recall that a Poincar\'{e}--Einstein manifold is one that satisfies the vacuum Einstein equation and has a conformal degeneracy at the boundary. In the paper by Graham and Lee~\cite{MR1112625}, solutions are constructed which are $C^{n-1,\gamma}$ close to the hyperbolic metric on the ball $\BB^n$ near the boundary. They showed that every such perturbation is determined by the conformal data on the boundary sphere. We will follow a similar idea here for the equation~\eqref{eq:supergravity}, replacing the Ricci curvature operator by the nonlinear supergravity operator, considering its linearization around one of the product solutions, and using a perturbation argument to show that all the solutions nearby are determined by three pairs of data on the boundary. 

Kantor studied this problem in his thesis~\cite{kantor2009eleven}, where he computed the indicial roots of the system and produced one family of solutions by varying along a specific direction of the 4-form. In this paper we use Hodge decomposition to get the same set of indicial roots and show that all the solutions nearby are prescribed by boundary data for the linearized operator, more specifically, the indicial kernels corresponding to three pairs of special indicial roots.

\subsection{Equations derived from the Lagrangian}\label{Lagrangianintro}
%Write the Lagrangian and say the three equations are derived. Then show that the product solution with constants satisfies the equations

The 11-dimension supergravity theory contains the following information on an 11-dimensional manifold $M$: a metric $g \in \Sym^2(M)$ and a 4-form $F \in \bigwedge^4(M)$.  The Lagrangian $L$ is
\begin{equation}
L(g,A)=\int_M R\Vol_g-\frac{1}{2}\bigg(\int_M F \wedge *F+\int_M \frac{1}{3}A\wedge F \wedge F\bigg).
\end{equation}
Here $R$ is the scalar curvature of the metric $g$, $A$ is a 3-form such that $F$ is the field strength $F=dA$. The first term is the classical Einstein--Hilbert action, and the second and the third one are respectively of Yang--Mills and Maxwell type. Note here we are only interested in the equations derived from the variation of Lagrangian and the variation only depends on $F=dA$, therefore we only need $F$ to be globally defined but with only locally defined $A_{i}$ on coordinate patches $U_{i}$. Then
\begin{equation}
\delta_{A_i}\bigg(\int_{U_i} A_i\wedge F \wedge F\bigg)=3\int_{U_i}\delta A_i \wedge F\wedge F
\end{equation}
shows that the variation of this term is $F\wedge F$.

The supergravity equations~\eqref{eq:supergravity} are derived from the Lagrangian above (see section~\ref{gauge}). Since the Ricci operator is not elliptic, we follow~\cite{MR1112625} and add a gauge breaking term 
\begin{equation}
\Phi(g,t)=\delta_g^*g\Lap_{gt}Id
\end{equation}
to the first equation. Then we apply $d*$ to the 2nd equation, and combine this with the third equation to obtain the system:
$$
Q: S^2(T^*M) \oplus  {\textstyle\bigwedge^4(M)} \rightarrow S^2(T^*M) \oplus {\textstyle\bigwedge^4(M)}
$$
\begin{equation}\label{gaugedequ}
\left(\begin{array}{c}
g\\
F
\end{array}\right)\mapsto
\left( \begin{array}{c}
\Ric(g)-\Phi(g,t)-F\circ F\\
d* (d*F+\frac{1}{2}F\wedge F)
\end{array} \right)
\end{equation}
which is the nonlinear system we will be studying.

\subsection{Edge metrics and edge Sobolev spaces }
%Edge vector filed V_e, edge operator definition, Q is a nonlinear edge operator
Edge differential and pseudodifferential operators were formally introduced by Mazzeo \cite{MR1133743}. The general setting is a compact manifold with boundary, M, where the boundary has in addition a fibration
$$
\pi:\partial M\rightarrow B.
$$
In the setting considered here, $M=\bbB^7\times \bbS^4$ is the product of a seven dimensional closed ball (identified as the hyperbolic space) and a four-dimensional sphere. The relevant fibration has the four-sphere as fibre:
$$
\xymatrix{
\bbS^{4} \ar[r] & \partial M=\bbS^6\times \bbS^4\ar[d]^{\pi} \\
 & \bbS^6.
}
$$

The space of edge vector fields $\mathcal{V}_e(M)$ is the Lie algebra consisting of those smooth vector fields on $M$ which are tangent to the boundary and such that the induced vector field on the boundary is tangent to the fibre of $\pi$. Another Lie algebra of vector fields we will be using is $\cV_b(M)$ which is the space of all smooth vector fields tangent to the boundary~\cite{melrose1993atiyah}. As a consequence, 
\begin{equation}\label{LieBracket}
\cV_e\subset \cV_b, \quad [\cV_e,\cV_b] \subset \cV_b.
\end{equation}

In terms of local coordinates, let $(x,y_1,y_2,...y_6)$ be the coordinates on the upper half space model for hyperbolic space $\HH^7$, and $(z_{1}, \dots, z_{4})$ be the coordinates on the sphere $\bbS^4$. Then locally $\cV_b$ is spanned by $\{x\partial_x, \partial_{y_{i}}, \partial_{z_{j}}\}$,  while $\cV_e$ is spanned by $\{x\partial_x, x\partial_{y_{i}}, \partial_{z_{j}}\}.$ The edge forms are the wedge product of the dual form to the edge vector fields $\cV_e$, with a basis:
$
\{\frac{dx}{x},\frac{dy^{i}}{x}, dz^{j}\},
$
and the edge 2-tensor bundle is formed by the tensor product of the basis. We will work on the following edge vector bundle:
\begin{definition}
Let $K$ be the edge bundles, the sections of which are symmetric 2-tensors and 4-forms:
\begin{equation}\label{K}
K:=\eS(M) \oplus \eW{4}.
\end{equation}
\end{definition}

Edge differential operators form the linear span of products of edge vector fields over smooth functions. We denote the set of $m$-th order edge operator as $\Diff{m}_e(M)$ and we will see that the supergravity operator $Q$ is a nonlinear edge differential operator, so a nonlinear combination of elements of $\Diff{*}_e(M)$.

Let $L^{2}_{e}(M)$ be the $L^{2}$ space with respect to the edge volume form which locally is given by $ x^{-7}{dx dy_{1}\dots dy_{6}} dz_{1}\dots, dz_{4}$. 
The edge-Sobolev spaces are given by
\begin{equation}
H_e^s(M)=\{u \in L^2_{e}(M)| V_1 \dots V_{k} u\in L^2_{e}(M),0\leq k \leq s, V_{i} \in \cV_{e}(M)\}.
\end{equation}
%$$
%u\in H_b^k(M,I) \Leftrightarrow V_b^k u\in L^2(M,I)
%$$ 
However, for purpose of regularity we are also interested in hybrid spaces with additional tangential regularity, as the existence of solutions with infinitely smooth b-regularity gives polyhomogeneous expansions. Therefore we define
the hybrid Sobolev space with both boundary and edge regularity:
\begin{equation}
 H_{e,b}^{s,k}(M) =\{ u\in H_e^s(M)|V_{1}\dots V_{i} u \in H_e^s(M), 0\leq i \leq k, V_{j} \in \cV_{b}(M) \}.
%\Leftrightarrow v_1...v_{s'}V_b^k u \in L_0^(M,I)
\end{equation}
By the commutation relation~\eqref{LieBracket}, $H_{e,b}^{s,k}(M)$ is well defined, that is, independent of the order in which edge and b-vector fields are applied, and the proof is given in proposition~\ref{A1}. For the vector bundle $K$ over $M$, the Sobolev spaces $H_{e,b}^{s,k}(M;K)$ can be similarly defined by choosing an orthonormal basis and are independent of choices.

These Sobolev spaces are defined so that edge operators map between them, i.e., for any $m$-th order edge operator $P \in \Diff{m}_e(M)$,  
\begin{equation}
P: H_{e,b}^{s,k}(M) \rightarrow H_{e,b}^{s-m,k}(M), m\leq s.
\end{equation}
which is proved in proposition~\ref{A2}.

\subsection{Poincar\'{e}--Einstein metrics on $\BB^7$}
The product of an arbitrary Poincar\'{e}--Einstein metric with the spherical metric provides a large family of solutions to this system.  For any Poincar\'{e}--Einstein metric $h$ with curvature $-6c^2$ with $c>0$, the following metric and 4-form gives a solution to equations~\eqref{eq:supergravity}:
\begin{equation}\label{PErelation}
u=\Big(h \times \frac{9}{c^2}g_{\bbS^4}, c\Vol_{\bbS^4}\Big).
\end{equation}

According to~\cite{MR1112625}, the Poincar\'{e}--Einstein metrics near the hyperbolic metric can be obtained by perturbation of the conformal boundary data. More specifically, there is the following result:

\begin{proposition}[\cite{MR0441129}] Let $M=\BB^{n+1}$ be the unit ball and $\hat h$ the standard metric on $\bbS^n$. For any smooth Riemannian metric $\hat g$ on $\bbS^n$ which is sufficiently close to $\hat h$ in $C^{2,\alpha}$ norm if $n>4$ or $C^{3,\alpha}$ norm if $n=3$, for some $0<\alpha <1$, there exists a smooth metric $g$ on the interior of $M$, with a $C^0$ conformal compactification satisfying 
$$
\operatorname{Ric}(g)=-ng, \quad g \text{ has conformal infinity } [\hat g].
$$
\end{proposition}

We are mainly interested in the solutions that are perturbations of such product solutions, in particular, we will focus on the solutions with $c=6$ in~\eqref{PErelation} and $h$ being the hyperbolic metric on the ball, i.e. on $X=\HH^7 \times \bbS^4$: 
\begin{equation}\label{u0}
u_{0}=(t,W)=\Big(g_{\bbH^{7}} \times \frac{1}{4}g_{\bbS^{4}}, 6\Vol_{\bbS^4}\Big),
\end{equation}
which is also known as the Freund--Rubin solution.

\subsection{Main theorem}
Besides the boundary conformal data that prescribes the Poincar\'{e}--Einstein metric, we will show that there are additionally three pairs of data on $\bbS^{6}$ that together parametrize the solution to~\eqref{eq:supergravity}.

First we define three bundles on $\bbS^6$ that correspond to the incoming and outgoing boundary data for the linearized supergravity operator.
\begin{definition}\label{bundle}
Let $V_1^{\pm}$ to be the space of 3-forms with $*_\bbS$ eigenvalue $\pm i$:
$$
V_1^{\pm}:=\{v_1^{\pm}\in C^{\infty}(\bbS^6;\bigwedge\nolimits^{\!3}T^*\bbS^6): *_{\bbS^6} v_1=\pm i v_1\}.
$$
Let $V_2^\pm$ and $V_3^\pm$ be the smooth functions on the 6-sphere tensored with eigenforms on 4-sphere:
$$
V_{2}^{+}=V_2^{-}:=\{v_2=f_{2}\otimes \xi_{16}: f_2\in C^{\infty}(\bbS^6), \xi_{16} \in E_{16}^{cl}(\bbS^4)\},
$$
$$
V_{3}^{+}=V_{3}^{-}:=\{v_3=f_{3}\otimes \xi_{40}: f_3 \in C^{\infty}(\bbS^6), \xi_{40} \in E_{40}^{cl}(\bbS^4)\},
$$
where $E_{\lambda}^{cl}(\bbS^{4})$ is the space of closed 1-forms with eigenvalue $\lambda$ on $\bbS^{4}$.
We also set
$$
V=\oplus_{i=1}^{3}V_{i}^{+}.
$$
\end{definition}

\begin{remark}
Note the dimension of the closed 1-forms with the first and second eigenvalues are determined by the degree 2 and 3 spherical harmonics in 4 variables, which, respectively, are 5 and 14 dimensional.
\end{remark}

We also require three numbers that define the leading term in the expansion of the solution, which come from indicial roots computation listed in Appendix~\ref{indicial}: 
\begin{equation}
\theta_1^{\pm} = 3\pm 6i,\ \theta_2^{\pm}=3 \pm i \sqrt{21116145}/1655,\ \theta_3^{\pm}=3 \pm i3 \sqrt{582842}/20098.
\end{equation}
From the indicial calculation, we also see that the real parts of all other indicial roots are at least distance 1 away from 3. And from here we fix a number 
$$\delta\in (0,1).$$

Then we define three scattering operators to relate the incoming and outgoing data:

\begin{definition}
Let $S_{i}, i=1,2,3$, be the scattering operators defined as
$$
S_{i}: V_{i}^{+}\rightarrow V_{i}^{-}
$$
for equation~\eqref{V1}, and equation~\eqref{V23} with eigenvalues 16 and 40, such that the linearized operator $dQ$ acting on the leading part of~\eqref{expansion} is contained in $O(x^{3+\delta})$.
\end{definition}

We also introduce the following notation of splitting of forms and tensors. For the product manifold $M=\HH^{7}\times \SS^{4}$, use $\pi_{\HH^{7}}$ and $\pi_{\SS^{4}}$ for the projection of $M$ onto the two components, then we can
define 
$$\eW{i,j}=\pi_{\HH^{7}}^{*}({}^{0}\wedge^{i}\HH^{7}) \wedge \pi_{\SS^{4}}^{*}(\wedge^{j}\SS^{4})$$
and we have the following identification
\begin{equation}
\eW{\ell}=\oplus_{i+j=\ell}\eW{i,j}
\end{equation}
That is, for $H\in \eW{\ell}$, we have
$$
H=\sum_{i+j=\ell} H_{(i,j)}
$$
where 
$H_{(i,j)}$ has the form $\sum_{k}f_{k}\alpha_{k} \wedge \beta_{k}, \ \alpha_{k} \in {}^{0}\wedge^{i}\bbH^{7}, \ \beta_{k} \in \wedge^{j}\bbS^{4}$. More specifically, locally near the boundary $H_{(i,j)}$ is given by the form
$$
x^{-i} dx \wedge dy_{I-1} \wedge dz_{J}+ x^{-i} dy_{I}\wedge dz_{J'}
$$
where $dy_{I-1}$ is an $(i-1)$-form in the span of $\{dy^{I_{1}}\wedge\dots dy^{I_{i-1}}\}$, similarly for $dy_{I}$ (an $i$-form) and $dz^{J}, dz^{J'}$ ($j$-forms in $z$ variables). There is also a similar notation for the 2-tensor $k$, where we take the linear map as an identification of edge $\ell$-tensors on $M$ with
$$
\sum_{i+j=\ell} \pi_{\HH}^{*}({}^{0}\wedge^{i}\HH^{7}) \wedge\pi_{\SS}^{*}(\wedge^{j}\SS^{4}),
$$
then we have the following decomposition for $k \in \eS(M)$
$$
k=\sum_{i+j=2}k_{(i,j)}.
$$

Now we introduce the boundary data to parametrize the solution:
\begin{definition}\label{d:expansion}[Leading expansion for the linear operator]
When we say the leading expansion is given by the outgoing data $(v_{1}^{+}, v_{2}^{+}, v_{3}^{+}) \in V$, this means that $(k,H)$ has an expansion with the following leading terms:
\begin{equation}\label{expansion}
\begin{aligned}
H_{(4,0)} &= \frac{dx}{x}\wedge (v_1^+  x^{\theta_1^+}+S_1( v_1^+)  x^{\theta_1^-}) + O(x^{3+\delta})\\
\operatorname{Tr}_{\HH^7}k&=7\delta_{\bbS^{4}}(v_2^+ x^{\theta_2^+}+S_2( v_2^+)  x^{\theta_2^-}+v_3^+ x^{\theta_3^+}+S_3( v_3^+) x^{\theta_3^-})+ O(x^{3+\delta})\\
\operatorname{Tr}_{\bbS^4}k&=4\delta_{\bbS^{4}}(v_2^+  x^{\theta_2^+}+S_2( v_2^+)  x^{\theta_2^-}+v_3^+ x^{\theta_3^+}+S_3( v_3^+)  x^{\theta_3^-})+ O(x^{3+\delta})\\
k_{(1,1)}&=d_{\bbH^{7}}\delta_{\bbS^{4}}( v_2^+ x^{\theta_2^+}+S_2( v_2^+)  x^{\theta_2^-}+v_3^+x^{\theta_3^+}+S_3( v_3^+) x^{\theta_3^-})+O(x^{3+\delta})\\
H_{(1,3)} &= d_{\bbH^{7}}*_{\bbS^{4}}(v_2^+ x^{\theta_2^+}+S_2( v_2^+)x^{\theta_2^-}+v_3^+ x^{\theta_3^+}+S_3( v_3^+)x^{\theta_3^-})+O(x^{3+\delta}) \\
H_{(0,4)} &=d_{\bbS^{4}} *_{\bbS^{4}} (v_2^+x^{\theta_2^+}+S_2( v_2^+)  x^{\theta_2^-}+v_3^+ x^{\theta_3^+}+S_3( v_3^+) x^{\theta_3^-})+O(x^{3+\delta})
\end{aligned}
\end{equation}
and the other components in $(k,H)$ are all in $O(x^{3+\delta})$.
\end{definition}

Now the main result is the characterization of the solution:
\begin{theorem*}
For  $k\gg 0$, $\delta\in (0,1)$, there exists $\rho>0$ and $\epsilon>0$, such that, for a Poincar\'{e}--Einstein metric $h$ that is sufficiently close to $g_{\bbH^{7}}$ with small difference in conformal boundary data $\|\hat h-g_{\SS^{6}}\|_{H^{k}(\SS^{6})}<\epsilon$, and any boundary value perturbation $v\in V$ with $\|v\|_{H^k(\bbS^{6};V)}<\rho$, there is a unique solution $u=(g, H) \in D_{v,h} \subset x^{-\delta}H_{e,b}^{s,k}(M;K)$ satisfying the supergravity equations~\eqref{eq:supergravity}, with the leading expansion of $(g-h\times \frac{1}{4}g_{\SS^{4}},H-6 \Vol_{\bbS^{4}})$ given by~\eqref{expansion}.
\end{theorem*}
\begin{remark}
For the base metric, we consider all $h\times \frac{1}{4}g_{\SS^{4}}$ where $h$ is a nearby Poincar\'{e}--Einstein metric close to $g_{\HH^{7}}$. And the difference of $h$ and $g_{\HH^{7}}$ is measured by their difference at the conformal infinity $\|\hat h-g_{\SS^{6}}\|_{H^{k}(\SS^{6})}$ for sufficiently large $k$, which by~\cite{MR1112625} implies that they are close as two Poincar\'{e}--Einstein metrics on $\BB^{7}$.
\end{remark}

Our approach is based on the implicit function theorem. From the boundary data $v$ we construct a perturbation term using the Poisson operator $P$, then consider a translation of the gauged operator: $Q^{h,v}(\cdot)=Q(\cdot+Pv)$ where the dependence on $h$ is in both the construction of the gauged term in the operator $Q$ and the Poisson operator $P$. A right inverse of the linearization, denoted $(dQ)^{-1}$, is constructed, and we show that $Q^{h,v}\circ (dQ)^{-1}$ is an isomorphism on the Sobolev space $x^{\delta}H_{e,b}^{s,k}(M;K)$ which is the range space of $dQ$. From here we deduce there is a unique solution to $Q$ for each boundary parameter set $v$.

To get the isomorphism result on $Q^{h,v}\circ (dQ)^{-1}$, we note that the model operator  on the boundary is $SO(5)$-invariant, and therefore utilize the Hodge decomposition of functions and forms on $\mathbb{S}^4$ to decompose the equations into blocks that allows us to compute the indicial roots for each block. 
Indicial roots are those $s$ that the indicial operator has a nontrivial kernel, and as mentioned above these roots are related to the leading order of the solution expansions near the boundary.

Once the indicial roots are computed, we construct the right inverse $(dQ)^{-1}$. The operator exhibits different properties for large and small eigenvalues. For large ones, the projected operator is already invertible by constructing a parametrix in the small edge calculus. For small eigenvalues, two resolvents $R_\pm=\lim_{\epsilon\downarrow 0} (dQ \pm i\epsilon)^{-1}$ are constructed and we combine them to get a real-valued right inverse.  We show that those elements corresponding to indicial roots with real part equal to 3 are the boundary perturbations needed in the theorem.

The paper is organized as follows. In section~\ref{gauge} we show the derivation of the equations from the Lagrangian and discuss the gauge breaking condition. In section~\ref{linearsection} we compute the linearization of the operator and its indicial roots. In section~\ref{Fredholm} we analyze the linearized operator, prove it is Fredholm and construct the boundary data. In section~\ref{fit}  we construct the solutions for the nonlinear equations using the implicit function theorem. 

\textbf{Acknowledgement:} I would like to thank Richard Melrose for many helpful discussions, ideas, and suggestions. I would also like to thank Robin Graham, Colin Guillarmou, and Rafe Mazzeo for many valuable comments on this project.

\paperbody
 
\section{Gauged operator construction}\label{gauge}

\subsection{Equations derived from Lagrangian}
As mentioned in Section~\ref{Lagrangianintro}, the supergravity system arises as the variational equations for the Lagrangian
$$
L(g,A)=\int_M R\Vol_g-\frac{1}{2}\bigg(\int_M F \wedge *F+\int_M \frac{1}{3}A\wedge F \wedge F\bigg)
$$
where $R$ denotes the scalar curvature of metric $g$, which is different from the Ricci curvature $R_{\alpha\beta}$. 
Now we compute its variation along two directions, namely, the metric and the form direction.
The first term is the Einstein-Hilbert action, for which the variation in $g$ is 
\begin{equation}
\delta_g \bigg(\int R\Vol_g \bigg)=\int\bigg(R_{\alpha \beta}-\frac{R}{2}g_{\alpha \beta}\bigg)\delta g^{\alpha \beta}\Vol_g.
\end{equation} 

The variation of the second term $F\wedge *F$ in the metric direction is  
\begin{multline}
\delta_g\bigg(\frac{1}{2}\int F\wedge *F\bigg)\\
=\frac{2}{4!}\int F_{\eta_1...\eta_4}F_{\xi_1...\xi_4}g^{\eta_2 \xi_2}g^{\eta_3 \xi_3} g^{\eta_4\xi_4}\delta g^{\eta_1\xi_1}\Vol_g
-\frac{1}{4}\int F\wedge *F g_{\alpha\beta}\delta g^{\alpha \beta}\Vol_{g}.
\end{multline}
 
Combining these we get the first equation on the metric
\begin{equation}\label{sg.1}
R_{\alpha\beta}-\frac{1}{2}Rg_{\alpha \beta}=\frac{1}{12}F_{\alpha \eta_1\eta_2\eta_3}F_{\beta}^{\eta_1\eta_2\eta_3}-\frac{1}{4}\langle F,F\rangle g_{\alpha \beta}.
\end{equation} 
Here $\langle\bullet,\bullet\rangle$ is the inner product on forms:
\begin{equation}
\langle F,F\rangle=\frac{1}{4!}F_{\eta_1...\eta_4}F^{\eta_1...\eta_4}.
\end{equation}
Taking the trace of the equation~\eqref{sg.1}, we get 
\begin{equation}
R=\frac{1}{6}\langle F,F\rangle.
\end{equation}
Finally, substituting $R$ into~\eqref{sg.1}, we get 
\begin{equation}
R_{\alpha\beta}=\frac{1}{12}(F_{\alpha \gamma_1\gamma_2 \gamma_3}F_{\beta}^{\gamma_1 \gamma_2 \gamma_3}-\frac{1}{12}F_{\gamma_1 \gamma_2 \gamma_3 \gamma_4}F^{\gamma_1 \gamma_2 \gamma_3 \gamma_4}g_{\alpha \beta}),
\end{equation}
which gives the first equation in~\eqref{eq:supergravity}.

The variation with respect to the 3-form $A$ is 
\begin{equation}
\delta_A L=\int \delta F\wedge * F-\frac{1}{6} \delta A \wedge F \wedge F-\frac{1}{3} A \wedge \delta F\wedge F=-\int \delta A\wedge (d*F+\frac{1}{2}F\wedge F),
\end{equation}
which gives the second supergravity equation:
\begin{equation}
d*F +\frac{1}{2}F\wedge F=0.
\end{equation}

Since $F$ is locally exact, we have the third equation
\begin{equation}
dF=0.
\end{equation}

\subsection{Poincar\'e--Einstein metric} Product solutions to the supergravity equations are obtained as follows: let $X$ be a 7-dimensional Einstein manifold with negative scalar curvature $\alpha<0$ and $K$ be a 4-dimensional Einstein manifold with positive scalar curvature $\beta>0$. Consider $X\times K$ with the product metric; then we have 
\begin{equation}
R_{\alpha\beta}=\left(\begin{array}{cc}
6\alpha g_{AB}^X & 0\\
0& 3\beta g_{ab}^K 
\end{array}
\right)
\end{equation}

Let $F=c\Vol_K$. A straightforward computation shows
\begin{equation}
(F\circ F)_{\alpha \beta}=\frac{c^2}{12}\left(\begin{array}{cc}
-2 g_{AB}^X & 0\\
0& 4 g_{ab}^K 
\end{array}
\right).
\end{equation}
Therefore any triple $(c,\alpha, \beta)$ satisfying 
\begin{equation}
-c^2/6=6\alpha, c^2/3=3\beta
\end{equation}
corresponds to a solution to the supergravity equation. 

\subsection{Edge bundles}
Such product metrics fit into the setting of edge bundles. As introduced in~\cite{MR1133743}, edge tangent bundles $\eT M$ are defined by declaring 
$\cV_{e}$ to be its smooth sections and its dual bundle, $\eT^{*}M$, is the edge cotangent bundle. We denote the edge form bundle,
\begin{equation}
{}^{e}\wedge^{m}(T^{*}M)=:\eW{m},
\end{equation}
of which the local sections can be written as the $\cC^{\infty}(M)$ combinations of 
$$
\frac{dx}{x}\wedge \frac{dy^{I_{1}}}{x} \dots \wedge \frac{dy^{I_{k}}}{x} \wedge dz^{j_{1}}\dots \wedge dz^{j_{l}}, 1+k+l=m
$$
and 
$$
\frac{dy^{I_{1}}}{x} \dots \wedge \frac{dy^{I_{k}}}{x} \wedge dz^{j_{1}}\dots \wedge dz^{j_{l}}, k+l=m.
$$
Similarly the edge symmetric 2-tensor bundle $\eS(M):=\Sym^{2}(\eT^{*} M)$ is spanned by 2-tensors with local forms of 
$$
\left(\frac{dx}{x} \ \frac{dy^{I}}{x} \ dz^{j}\right)
\left(\begin{array}{ccc}
k_{00} & k_{0J} &k_{0j} \\ 
k_{I0} & k_{IJ} & k_{Ij} \\
k_{i0} & k_{iJ} & k_{ij}
\end{array}\right)
\left(
\begin{array}{c}
\frac{dx}{x}\\
\frac{dy^{J}}{x}\\
dz^{j}
\end{array}
\right)
$$
with smooth coefficients $k_{**}$.

It is easy to check that the supergravity operator $S$ is an edge operator
\begin{equation}\label{S}
\begin{array}{c}
S: {}^{e}\Sym^{2}(M)\oplus \eW{4} \rightarrow {}^{e}\Sym^{2}(M)\oplus \eW{8} \oplus  \eW{5} \\
\left(\begin{array}{c}
g\\
F 
\end{array}\right) 
\mapsto
\left(
\begin{array}{c}
\Ric g - F \circ F\\
d*F +\frac{1}{2} F\wedge F \\
dF
\end{array}
\right).
\end{array}
\end{equation}

\subsection{A square system} 
To get a square system, we apply $d*$ to the second equation. Because of the closed condition $dF=0$, $d*d*F$ is the same as $\Lap F$. This leads to the following square system (here ${}^{e}\wedge^{4}_{cl}M$ denotes the bundle of closed edge 4-forms on $M$):
\begin{equation}\label{SS}
\begin{array}{c}
\tilde S:  \eS\oplus {}^{e}\wedge^{4}_{cl}M \rightarrow \eS\oplus {}^{e}\wedge^{4}_{cl}M\\
\left(\begin{array}{c}
g\\
F 
\end{array}\right) 
\mapsto
\left(\begin{array}{c}
\Ric g-F\circ F\\
\Lap F+\frac{1}{2}d*(F\wedge F)
\end{array} \right)
\end{array}
\end{equation}

\begin{proposition}\label{p:S}
The kernel of the square supergravity operator $\tilde S$~\eqref{SS} is the same as the original supergravity operator $S$~\eqref{S}:
\begin{equation}
\Nul(S) \cap x^{\delta}H_{e,b}^{2,k}(M;K) = \Nul(\tilde S)\cap x^{\delta}H_{e,b}^{2,k}(M;K)
\end{equation}
\end{proposition}
\begin{proof}
We only need to show that, in $x^{\delta}H^{2,k}_{e,b}(M;K), $ the null space of $\tilde S$ does not have extra elements that are not in $\Nul(S)$.
We show this by proving that if 
\begin{equation}\label{harmonic}
\begin{aligned}
d*\left(d*F+\frac{1}{2}(F\wedge F)\right)=0\\
dF=0
\end{aligned}
\end{equation}
then 
\begin{equation}\label{vanish}
\omega:=d*F+\frac{1}{2}(F\wedge F)=0.
\end{equation}
Note that~\eqref{harmonic} implies that $\omega$ is a harmonic form on $\bbH^{7}\times \bbS^{4}$. 

Consider the Hodge decomposition of forms on $\bbS^{4}$ by taking $\alpha_{i}$ and $\beta_{i}$ to be the basis of coclosed and closed forms for the $i$-th eigenvalue 
%(assume the multiplicity of each eigenvalue is one):
\begin{equation}\label{omega}
\begin{aligned}
d_{\SS^{4}}\alpha_{i}=\lambda_{i}\beta_{i}, \ \delta_{\SS^{4}} \alpha_{i}=0\\
 \delta_{\SS^{4}} \beta_{i}=\lambda_{i} \alpha_{i}, \ d_{\SS^{4}}\beta_{i}=0
\end{aligned}
\end{equation}
This implies $\Lap_{\SS^{4}} \alpha_{i}=\lambda_{i}^{2}\alpha_{i}, \Lap_{\SS^{4}}\beta_{i}=\lambda_{i}^{2}\beta_{i}$. Write $\omega=\sum_{i=1}^{\infty}u_{i}\alpha_{i} + v_{i} \beta_{i}$, where $u_{i},v_{i}$ are forms on $\bbH^{7}$. The decay condition on $F$ implies that $u_{i}, v_{i}$ are $L^{2}$ forms on $\bbH^{7}$. Using the eigenspace decomposition to rewrite the closed and coclosed condition for $\omega$ and combining with~\eqref{omega}, we get
\begin{align}
d_{\HH^{7}}u_{i}=0, \ \lambda_{i} u_{i} + d_{\HH^{7}} v_{i}=0\\
\delta_{\HH^{7}}v_{i}=0, \ \delta_{\HH^{7}}u_{i}+\lambda_{i} v_{i}=0.
\end{align}
Then we get 
$$
\Lap_{\HH^{7}}u_{i}=\lambda_{i}^{2} u_{i}, \ \Lap_{\HH^{7}}v_{i}=\lambda_{i}^{2} v_{i}.
$$
Since there are no $L^{2}$ eigenforms on $\bbH^{7}$~\cite{MR961517}, we get $u_{i}=0, v_{i}=0$, which proves~\eqref{vanish}. 
\end{proof}

\subsection{Gauge condition} 
%The solution to the gauge-modified equations also satisfies the original equations.
Following~\cite{MR1112625} in the setting of a Poincar\'{e}--Einstein metric, we add a gauge operator to the curvature term where $t=g_{\HH^{7}}\times \frac{1}{4}g_{\SS^{4}}$ is the background metric:
\begin{equation}\label{e:gaugeterm}
\Phi(g,t)=\delta_g^*gt^{-1}\delta_gG_gt.
\end{equation}
Here
$$
[G_gt]_{ij}=t_{ij}-\frac{1}{2}t_k^kg_{ij}, \quad [\delta_gt]_i=-t_{ij}^j,
$$
$\delta_g^*$ is the formal adjoint of $\delta_g$, which can be written as
$$
[\delta^*_g w]_{ij}=\frac{1}{2}(w_{i,j}+w_{j,i}),
$$
and $gt^{-1}$ is the endomorphism of $T^{*}M$ given by
$$
[gt^{-1}w]_i=g_{ij}(t^{-1})^{jk}w_k.
$$
Note that another way to write the gauge term in~\eqref{e:gaugeterm} is given in~\cite[(4.1)]{kantor2009eleven} as
$$
\Phi(g,t)=\delta_{g}^{*}g \Delta_{gt}Id.
$$

By adding the gauge term to the first equation of $\tilde S$ we get an operator $Q$, which is a map from the space of symmetric 2-tensors and closed 4-forms to itself:
$$
Q:  \eS\oplus {}^{e}\wedge^{4}_{cl}M \rightarrow \eS\oplus{}^{e}\wedge^{4}_{cl}M$$
\begin{equation} \label{operatorQ}
\left(\begin{array}{c}
g\\
F
\end{array}\right)\mapsto
\left( \begin{array}{c}
\Ric(g)-\Phi(g,t)-F\circ F\\
\Lap F+\frac{1}{2}d*(F\wedge F)
\end{array} \right)
\end{equation}
which is the main object of study below.

As discussed in Lemma 2.2 in~\cite{MR1112625}, $\Ric(g)+ng-\Phi(g,t)=0$ holds if and only if $id: (M,g)\rightarrow (M,t)$ is harmonic and $\Ric(g)+ng=0$ when $(t,g)$ satisfies certain regularity restrictions. We will show that the gauged equations here yield the solution to the supergravity equations in a similar manner.

We first prove a gauge elimination lemma for the linearized operator $dQ$ which is computed in Proposition~\ref{Linearization}. As can be seen from (\ref{operatorQ}), only the first part (the map on 2-tensors) involves the gauge term, therefore we restrict the discussion to the first part of $dQ$. We use $dQ_g(k,H)$  to denote the linearization of the tensor part of $Q$ given by
$$
\Ric(g)-\Phi(g,t)-F\circ F
$$ 
along the metric direction at the point $(t,W)$, which acts on $(k,H)\in \Gamma(K)$. Also $d\tilde S_{g}(k,H)$ is defined similarly to $dQ_{g}(k,H)$. And we define $d\Phi(k)_{t}$ as the linearization of $\Phi(g,t)$ along the direction of first variable at $t$ while the second variable is fixed at $t$. 

First we give the following gauge-breaking lemma for the linearized operator, which is adapted from Theorem 4.1 and Theorem 4.2 in~\cite{kantor2009eleven}. For a 1-form $v$, we define $v^{\sharp}$ to be the dual vector field of $v$ with respect to $t$, and define $L_{v^{\sharp}}g$ to be Lie derivative of $g$ along $v^{\sharp}$.

\begin{proposition}\label{lineargauge}
For fixed $\delta\in (-1,0)$, there exists $\epsilon>0$, such that for any $\|g-t\|_{x^{\delta}H^{2}_{e}(M;\eS(M))}<\epsilon$ and $k\in x^{\delta}H^{2}_{e}(M;\eS(M))$ satisfying $dQ_{g}(k,H)=0$, there exists a 1-form $v$ and $\tilde k=k+L_{v^\sharp}g$ such that $d\tilde S_{g}(\tilde k, H)=0$.
\end{proposition}

To prove the proposition, we first determine the equation to solve for such a 1-form $v$, which appeared in Theorem 4.1 in~\cite{kantor2009eleven}.
\begin{lemma}\label{klemma} 
Given $k\in \eS(M)$, if a 1-form $v$ satisfies
\begin{equation}\label{gaugeop}
((\Lap^{rough}-\Ric)v)_\lambda= \frac{1}{2}(2\nabla^\alpha k_{\alpha\lambda}-\nabla_\lambda Tr_g(k))
\end{equation}\label{kcondition}
then the 2-tensor $\tilde k=k+L_{v^\sharp}g$ satisfies the gauge condition
$$
d\Phi_{t}(\tilde k)=0.
$$
\end{lemma}

\begin{proof}
Following the proof of \cite[Theorem 4.2]{kantor2009eleven}, let $\Psi(v,g)$ be the map
$$
\Psi(v,g)^k=(\phi^*_{v^\sharp}g)^{\alpha\beta}(\Gamma_{\alpha\beta}^k(\phi^*_{v^\sharp}g)-\Gamma_{\alpha\beta}^k(t)),
$$
where $\phi^{*}_{v^{\sharp}}$ is the diffeomorphism $\exp^{t} (v^{\sharp},1)$
(following
the exponential flow for the metric $t$ to time 1 in the direction of $v^{\sharp}$).
Let $D_{1}\Psi(0,t)$ and $D_{2}\Psi(0,t)$ be the linearization of $\Psi$ along the first and second variable at $(v,g)=(0,t)$, then they satisfy
$$
\delta_t^*t D_1 \Psi(0,t)(v)=d\Phi_{t}(L_{v^\sharp}g), \ \delta_t^*t D_2\Psi(0,t)(k)=d\Phi_{t}(k),
$$
which are proved in \cite[Theorem 4.2]{kantor2009eleven}. 
Therefore in order to get $d\Phi_{t}(\tilde k)=0$, we only need
$$
-D_1 \Psi(0,t)(v)=D_2\Psi(0,t)(k).
$$
The left hand side can be reduced to 
$$
-g^{\alpha\beta}\nabla_\alpha \nabla_\beta v^{k}-R^k_{\mu}v^{\mu}=(\Lap^{rough}-\Ric)v^k
$$
and right hand side is
$$
\frac{1}{2}g^{\alpha\beta}g^{k\lambda}(\nabla_\alpha k_{\beta\lambda}+\nabla_\beta k_{\alpha\lambda}-\nabla_\lambda k_{\alpha\beta}). 
$$
Lowering the index on both side, we get
$$
((\Lap^{rough}-\Ric)v)_\lambda= \frac{1}{2}(2\nabla^\alpha k_{\alpha\lambda}-\nabla_\lambda Tr_g(k)).
$$
\end{proof}

Next we discuss the solvability of the operator defined in the left hand side of~\eqref{gaugeop}.
\begin{lemma}
If $|\delta|<1$, then at the point $t=g_{\HH^{7}}\times \frac{1}{4}g_{\SS^{4}}$ the operator
$$
\Lap^{rough}-\Ric:x^{\delta}H^2_{e}({}^eT^*M) \rightarrow x^\delta L^2_{e}({}^eT^*M)
$$
is an isomorphism.
\end{lemma}

\begin{proof}
Using the splitting
\begin{equation}\label{e:split}
{}^eT^*M \cong \pi_\bbH^*{}^e T^*\bbH^7 \oplus \pi_\bbS^* {}^e T^* \bbS^4
\end{equation}
and the product structure of the metric, we write the operator as
$$
\Lap^{rough}-\Ric=\Lap_{\HH^7}^{rough} +\Lap_{\SS^4}^{rough}-\operatorname{diag}(-6,12).
$$
From~\cite[Propositon 4.1]{kantor2009eleven}, the operator $L$ is diagonal with respect to the splitting~\eqref{e:split}, so we only need to consider the following two operators: $L_{H}=L|_{\pi_\bbH^*{}^e T^*\bbH^7}$ and $L_{S}=L|_{\pi_\bbS^* {}^e T^* \bbS^4}$.
Decomposing them with respect to eigenfunctions on the 4-sphere, we get the following two families of operators (see~\cite[Corollary 4.1]{kantor2009eleven} for derivation):
\begin{equation}
L_{H}^{\lambda}=\Lap_\bbH+\lambda_{H}-24: C^{\infty}(\bbH^7) \rightarrow C^{\infty}(\bbH^7),
\end{equation}
\begin{equation}
L_{S}^{\lambda}=\Lap^{rough}_\bbH+\lambda_{S}+6:{\textstyle{}^e\wedge^*\bbH^7} \rightarrow {\textstyle{}^e\wedge^*\bbH^7}
\end{equation}
Consider the smallest eigenvalue in each case: $\lambda_{H}=16, \lambda_{S}=0$. The indicial radius, defined as the smallest number $R$ such that the operator has an indicial root $3+R$ (see~\cite[Definition 5]{kantor2009eleven}), for $L_{H}$ is 1 and for $L_{S}$ is 4 (again see~\cite[Corollary 4.1]{kantor2009eleven} for the computation). By the same argument as in proposition~\ref{small}, $L_{H}: x^{\delta}H^{2}_{e}(\HH^{7})\rightarrow x^{\delta}L^{2}_{e}(\HH^{7})$ is an isomorphism.
The same argument holds for $L_{S}$, which is also Fredholm and an isomorphism on $x^{\delta}H^2_{e}({}^e\wedge^*\bbH^7) \rightarrow x^{\delta}L^2_{e}({}^e\wedge^*\bbH^7 )$ for $|\delta|<4$.

Combining the statements for $L_{H}$ and $L_{S}$, we conclude that $\Lap^{rough}-\Ric$ is an isomorphism between $x^{\delta}H^2_{e}({}^eT^*M) \rightarrow x^\delta L^2_{e}({}^eT^*M)$ for $|\delta|<1$.

\end{proof}

The isomorphism holds true for metrics nearby, by a simple perturbation argument.
\begin{corollary}\label{gaugeiso}
There exists $\epsilon>0$, such that for any metric $g$ with\\$\|g-t\|{}_{x^{\delta}H^{2}_{e}(M;\eS(M))}<\epsilon,$ for $|\delta|<1$, $\Lap^{rough}-\Ric$ is an isomorphism as a map 
$$
\Lap^{rough}-\Ric:x^{\delta}H^2_{e}({}^eT^*M) \rightarrow x^\delta L^2_{e}({}^eT^*M).
$$
\end{corollary}
\begin{proof}
Let $A_{g}=\Lap^{rough}-\Ric$, by writing the coefficients out, we have that for any $u \in x^{\delta}H^2({}^eT^*M) $,
$$\|(A_{g}-A_{t})u\|_{x^\delta L^2_{e}({}^eT^*M)} \leq C\|g-t\|_{x^{\delta}H^{2}_{e}(M;\eS(M))}\|u\|_{x^{\delta}H^{2}_{e}({}^eT^*M)},$$
which shows that $A_{g}$ is also an isomorphism for $g$ sufficiently close to $t$.
\end{proof}

With the lemmas above, we can prove the proposition. 

\begin{proof}[Proof of Proposition \ref{lineargauge}] From Corollary~\ref{gaugeiso} we know $\Lap^{rough}-\Ric:x^{\delta}H^2_{e}({}^eT^*M) \rightarrow x^\delta L^2_{e}({}^eT^*M)$ is an isomorphism for $g$ close to $t$, therefore there exists a one-form $v$ satisfying~\ref{kcondition}. Then from Lemma~\ref{klemma}, $\tilde k=k+L_{v^\sharp}g$ satisfies $d\Phi_{g}(\tilde k)=0$. Putting it back to the linearized equation, we get $d\tilde S_g(\tilde k,H)=0$.
\end{proof}

Next we prove the nonlinear version of gauge elimination by using .

\begin{proposition}\label{gaugeeq}
If a metric and a closed 4-form $(g,V)$ satisfies the gauged equations $Q(g,V)=0$, then there is an diffeomorphism $g\mapsto \tilde g$ such that $\Phi(\tilde g,t)=0$ and $(\tilde g,V)$ is a solution to equation~\eqref{eq:supergravity} i.e. $S(\tilde g,V)=0$.
\end{proposition}

\begin{proof}
Consider the vector field in the affine Sobolev space $\{t\}+\ x^{\delta}H^{2}_{e}(M;\eS(M))$ with $\delta \in (-1,0)$, defined by
\begin{equation}
X_{g}=\tilde k_{g}=k_{g}+L_{v^{\sharp}}g,
\end{equation}
It is easy to check that for any $g\in  x^{\delta}H^{2}_{e}(M;\eS(M))$,  $\tilde k_{g} \in x^{\delta}H^{1}_{e}(M;\eS(M))$.
Now consider the integral curve $g(s)$ starting from $g(0)=g$. From the integral curve theory on an infinite dimensional manifold (see for example Proposition 1 in Chapter 4 of~\cite{MR772023}), since the map $g\mapsto \tilde k_{g}$ is Lipchitz, the integral curve exists and we take $\tilde g=g(1)$. From the construction of $\tilde k_{g}$ we know that $\tilde g$ satisfies
$
\tilde S(\tilde g, V)=0.
$
And hence it also satisfies $S(\tilde g, V)=0$ by Proposition~\ref{p:S}.
\end{proof}

\section{The linearized and the indicial operator}\label{linearsection}
We now consider the linearization of the gauged supergravity operator near the base solution $(t,W)=(g_{\HH^7}\times \frac{1}{4}g_{\bbS^4}, 6\Vol_{\bbS^4})$.  
The first step in proving the invertibility of $Q$ as an edge operator is to compute the indicial roots and indicial kernels of this linearized operator $dQ$, which is done with respect to the Hodge decomposition on the 4-sphere. There is a pair of indicial roots associated to each eigenvalue and as the eigenvalues becomes larger, the pairs move apart. The distribution of indicial roots are illustrated in Figure~\ref{figure:indicial}. Note that the operator $dQ$ is formally self-adjoint, and we will consider its spectral theory.

With respect to the volume form on $\bbH^7\times \bbS^4$, there is an inclusion of weighted functions and forms. Since the edge volume form is locally given by $$x^{-7}{dx dy_{1}\dots dy_{6}} dz_{1}\dots, dz_{4},$$ 
we can see that for $\Re(s)>3$, 
$$
x^{s}\mathcal{C}^{\infty}(M) \subset L^2_e(M),
$$
while for any $\Re(s)\leq 3$, $x^{s} \notin L^2_e(M)$. Therefore $\Re(s)=3$ line is the $L^2$ cutoff line. There are only three pairs of exceptional indicial roots corresponding to the lowest three eigenvalues that lie on the $L^{2}$ line, in the sense that they have real parts equal to 3 and the remainders are pure imaginary and symmetric around this line.

\subsection{Linearization of the operator $Q$}
The nonlinear supergravity operator contains two parts: the gauged curvature operator $\operatorname{Ric}(\cdot)-\Phi(t, \cdot)$ with a nonlinear part $F\circ F$, and the second order differential operator $d*d*F$ with a nonlinear part $d*(F\wedge F)$. Note that since the Hodge operator $*$ depends on the metric,  the linearized operator couples the metric with the 4-form in both equations. The computation below is similar to Proposition 5.1 and 5.2 in section 5 of~\cite{kantor2009eleven}.

\begin{proposition}\label{Linearization}
The operator $Q: K\rightarrow K$ has linearization at $(t, W)$:
$$
dQ_{t,W} : \Gamma(K) \rightarrow \Gamma(K)
$$
\begin{equation}\label{eq:DQ}
\left(\begin{array}{c}
k\\
H\
\end{array}\right)
\mapsto
 \left(\begin{array}{c}
\Lap_{t}^{rough} k+\operatorname{L}\\
 d*(d*H+6\Vol_{\bbS^{4}} \wedge H +6d*_{\bbH^{7}} k_{(1,1)}+3d(\tr_{\bbH^{7}}(k)-\tr_{\bbS^{4}}(k))\wedge \Vol_{\bbH^{7}})
\end{array}\right)
\end{equation}
where the lower order term $\operatorname{L}$ with respect to the product $\bbS^{4}\times \HH^{7}$ is given by:
\begin{multline}
\operatorname{L}=\\
\left(
\begin{array}{cc}
-k_{IJ}-6\tr_{\bbS^{4}}(k) t_{IJ}+\tr_{\bbH^{7}} (k)t_{IJ}+2*_{\bbS^{4}} H_{(0,4)}t_{IJ} & 
6k_{(1,1)}-3*_{\bbS^{4}} H_{(1,3)}\\
6k_{(1,1)}-3*_{\bbS^{4}} H_{(1,3)}&
4k_{ij}+8\tr_{\bbS^{4}} (k)t_{ij}-*_{\bbS^{4}} H_{(0,4)}t_{ij}
\end{array}\right).
\end{multline}
\end{proposition}
Here the notation $\cdot_{(i,j)}$ is the decomposition of forms and tensors with respect to the product structure, as explained before Definition~\ref{d:expansion}.
The computation involves a curvature part and a form part. First of all, the variation of the gauged curvature operator is given by the following lemma.

\begin{lemma}
For $k \in \eS(M)$, the linearization of the gauged Ricci operator at the base metric $t$ is 
\begin{equation}
d(\operatorname{Ric}(\cdot)-\Phi(t, \cdot))_{t}(k)=\frac{1}{2} \Lap_{t}^{rough} k + R(k),
\end{equation}
where
$$
R(k)=\left(\begin{array}{cc} -7k_{IJ} +\operatorname{Tr}_{\HH^7}(k) t_{IJ} & 0 \\
0 & 16 k_{ij}-  \operatorname{Tr}_{\bbS^4}(k) t_{ij} \end{array} \right).
$$
\end{lemma}

\begin{proof}
Following the result in~\cite{MR1112625}, the linearization of the gauged operator at the base metric $t$ is
\begin{equation}
d(\operatorname{Ric}-\Phi(t,\cdot))_{t}(k)=\frac{1}{2}\Lap_t^{rough}k+k^{\alpha\beta}R_{\beta\gamma\delta\alpha}+\frac{1}{2}(R^{\beta}_\gamma k_{\beta\delta}+R^\beta_\delta k_{\beta\gamma}).
\end{equation}
Specifically, if the metric $t$ has constant sectional curvature which is the case here, the curvature term is diagonal and can be written as
\begin{equation}\label{diagonalcurvature}
R_{\alpha \beta\delta \gamma}=-(t_{\alpha\delta}t_{\gamma\beta}-t_{\alpha\beta}t_{\gamma\delta}),
\end{equation}
so the linearization of this total operator is as above.
\end{proof}

\begin{lemma}
The linearization of the term $F\circ F$ acting on a 2-tensor $k \in \eS(M)$ and a 4-form $H \in \eW{4}$ are respectively: 
\begin{multline}\label{3.6}
d(F\circ F)_{t,W}(k)\\
=\left(
\begin{array}{lc}
 \frac{1}{36}Tr_{S}(k)t_{IJ}-\frac{1}{36}k  &\  \frac{1}{144}\langle W,W\rangle k_{Ij}\\
 \noalign{\vskip.2in}
\frac{1}{144}\langle W,W \rangle k_{Ij} & \ 
\begin{array}{l}
\frac{1}{12}\big(-3W_{a i_1i_2i_3}W_{b j_1}^{j_2j_3}t^{i_1l_1}k_{l_1l_2}t^{l_2j_1}\\
+\frac{1}{3}F_{i_1i_2i_3i_4}F_{j_1}^{i_1i_2i_3}t^{j_1l_1}k_{l_1l_2}t^{l_2i_1}t_{ab}-\frac{1}{12}\langle W,W \rangle k_{ab}\big) 
\end{array}
\end{array}
\right),
\end{multline}
and 
\begin{equation}\label{3.7}
d(F\circ F)_{t,W}(H)=\left(
\begin{array}{cc}
2*_{\SS^4} H_{(0,4)}t_{AB} & 3(*_{\SS^4}H_{(1,3)})_{Ab} \\
3(*_{\SS^4}H_{(1,3)})_{aB} &
\begin{array}{c}
\frac{1}{6}H_a^{i_1i_2i_3}W_{bi_1i_2i_3}\\
-\frac{1}{72}W_{i_1i_2i_3i_4}H^{i_1i_2i_3i_4}t_{ab}
\end{array}
\end{array}
\right).
\end{equation}
\end{lemma}
\begin{proof}
The proof is by direct computation. Note that
for metric variation $k \in \eS(T^*M)$
\[
D_{t,W}(F\circ F)(k)=\left(
\begin{array}{cc}
I & II\\
II & III 
\end{array}
\right),
\]
\begin{equation}
I_{AB}=\frac{1}{12}\Big(\frac{4}{12}F_{i_1i_2i_3i_4}F_{j_1}^{i_1i_2i_3}t^{j_1l_1}k_{l_1l_2}t^{l_2i_1}t_{AB}-\frac{1}{12}W_{1_1i_2i_3i_4}W^{i_1i_2i_3i_4}k_{AB}\Big),
\end{equation}
\begin{equation}
II_{Ab}=\frac{1}{144}W_{i_1i_2i_3i_4}W^{i_1i_2i_3i_4}k_{Ab},
\end{equation}
\begin{multline}
III_{ab}=\frac{1}{12}\Big(-3W_{a i_1i_2i_3}W_{b j_1}^{j_2j_3}t^{i_1l_1}k_{l_1l_2}t^{l_2j_1}+\frac{4}{12}F_{i_1i_2i_3i_4}F_{j_1}^{i_1i_2i_3}t^{j_1l_1}k_{l_1l_2}t^{l_2i_1}t_{ab}\\
-\frac{1}{12}W_{i_1i_2i_3i_4}W^{i_1i_2i_3i_4}k_{ab}\Big),
\end{multline}
which using inner product $W_{i_1i_2i_3i_4}W^{i_1i_2i_3i_4}=\langle W,W \rangle $ will give the expressions above.

And in the 4-form direction we have
\[
D_{t,W}(F\circ F)(H)=
\left(
\begin{array}{cc}
I & II\\
II & III 
\end{array}
\right),
\]
where 
\begin{equation}
I_{AB}=-\frac{1}{72}W_{i_1i_2i_3i_4}H^{i_1i_2i_3i_4}t_{AB}=2*_\bbS H_{(0,4)}t_{AB};
\end{equation}
\begin{equation}
II_{Ab}=\frac{1}{12}H_{Ai_1i_2i_3}W_b^{i_1i_2i_3}=3(*_\bbS H_{(1,3)})_{Ab},
\end{equation}
\begin{equation}
III_{ab}=\frac{1}{6}H_a^{i_1i_2i_3}W_{bi_1i_2i_3}-\frac{1}{72}W_{i_1i_2i_3i_4}H^{i_1i_2i_3i_4}t_{ab}.
\end{equation}
Combing those two we get the expression~\eqref{3.6} and~\eqref{3.7}.
\end{proof}

Next we compute the linearization for the second part:
\begin{lemma}
The linearization of the equation 
$$
d*F+\frac{1}{2}F\wedge F=0,
$$
in the form and tensor directions respectively are:
\begin{equation}
d(d*F+\frac{1}{2}F\wedge F)_{t,W}(H)=d*H+H\wedge F,
\end{equation}
\begin{equation}
d(d*F+\frac{1}{2}F\wedge F)_{t,W}(k)=6d*_{\bbH^{7}}k_{(1,1)}+3d(\tr_{\bbH^{7}}(k)-\tr_{\bbS^{4}}(k))\Vol_{\bbH^{7}}.
\end{equation}
\end{lemma}

\begin{proof}
The linearization along the form direction is straight-forward, as the terms are linear and quadratic in $F$. Along the metric direction, the linearization comes from the Hodge star:
\begin{equation}
\begin{array}{ll}
D(*F)_{\beta_1\beta_2..\beta_7}(k)&\\
=D(\frac{1}{4!}V^{\alpha_1..\alpha_4}_{\beta_1..\beta_7}W_{\alpha_1..\alpha_4})(k)=\frac{1}{4!}(\delta V)^{\alpha_1..\alpha_4}_{\beta_1..\beta_7}W_{\alpha_1..\alpha_4}+\frac{4}{4!}V^{\alpha_2..\alpha_4}_{\gamma_1\beta_1..\beta_7}(\delta g)^{\gamma_1\alpha_1}W_{\alpha_1..\alpha_4}\\
=\frac{1}{2}\frac{1}{4!}t^{\alpha \beta}k_{\alpha\beta}V^{\alpha_1..\alpha_4}_{\beta_1..\beta_7}W_{\alpha_1..\alpha_4}+\frac{1}{6}V^{\alpha_2..\alpha_4}_{\gamma_1\beta_1..\beta_7}t^{\gamma_1 \xi}k_{\xi \psi}t^{\psi \alpha_1}W_{\alpha_1..\alpha_4}\\
%=\frac{1}{2}Tr(k)(*W)_{\beta_1..\beta_7}+
=6d*_{\bbH^{7}} k_{(1,1)}+3d(\tr_{\bbH^{7}} (k)-\tr_{\bbS^{4}} (k))\Vol_{\bbH^{7}}
\end{array}
\end{equation}
which gives the expressions above.
\end{proof}

\begin{proof}[Proof of Proposition~\ref{Linearization}]
Combining the components above, the linearized equations are

\begin{equation}
\begin{array}{l}
\frac{1}{2}\Lap_{t}^{rough}k+k^{\alpha\beta}R_{\beta\gamma\delta\alpha}+\frac{1}{2}(R^{\beta}_\gamma k_{\beta\delta}+R^\beta_\delta k_{\beta\gamma})+L=0\\
d*(6d*_{\bbH^{7}} k_{(1,1)}+3d(\tr_{\bbH^{7}} (k)-\tr_{\bbS^{4}} (k))\Vol_{\bbH^{7}}+d*H+H\wedge \Vol_{\bbH^{7}})=0
\end{array}
\end{equation}
which  after rearrangement gives equation~\eqref{eq:DQ}.
\end{proof}

Here we only computed the linearization at the base product metric $t=g_{\bbH^7}\times \frac{1}{4}g_{\bbS^4}$. However, if the hyperbolic metric is changed to other Poincar\'{e}--Einstein metrics satisfying the relation~\eqref{PErelation}, the only difference is the equation~\eqref{diagonalcurvature} which holds up to terms vanishing at the boundary, in particular, \eqref{diagonalcurvature} is modified to be
$$
R_{\alpha \beta\delta \gamma}=-(t_{\alpha\delta}t_{\gamma\beta}-t_{\alpha\beta}t_{\gamma\delta})+x^{-3}\mathcal{E}^{1}(x^{2}t) + x^{-2}\mathcal{E}^{2}(x^{2}t),
$$
where $\mathcal{E}^{i}(x^{2}t)$ are tensors 
whose components in any coordinate system smooth up to $\partial M$ are polynomials, with coefficients in $C^{\infty}(\overline M)$, in the components of $x^{2}t, (x^{2}t)^{-1}$ and their partial derivatives, such that in each term the total number of derivatives that appear is at most $i$ (see~\cite{MR1112625}). 
This change does not affect computations later. The indicial operator of $dQ$ is the same, as the model operator acts on the tangent space at each boundary point (or consider the space obtained by rescaling a neighborhood), which is the standard hyperbolic upper space.

\subsection{Indicial roots computation}
Having obtained the linearized operator $dQ$, we next compute its indicial roots on the boundary of $\HH^7$, which together with the indicial kernels will give the parametrization of the kernel of this linear operator. Using the Hodge decomposition on the 4-sphere, the operator $dQ$ acts on sections of $\wedge^{*} \HH^7$ tensored with the finite dimensional eigenspaces of $\wedge^{*}\bbS^4$.

\begin{lemma}\label{Kdecomp}
Sections of the bundle $K$ decompose with respect to the Hodge decomposition of $\bbS^{4}$.
\end{lemma}
\begin{proof}
As mentioned before definition~\ref{d:expansion}, we identify the symmetric edge 2-tensor bundle with $$(\Sym^2({}^e T^*\HH^7))\oplus( {}^eT^*\HH^7 \otimes T^*\bbS^4 )\oplus \Sym^2(T^*\bbS^4),$$ and decompose the 4-form bundle according to its degree on $\HH^7$ and $\bbS^4$, i.e. 
$${}^e \wedge^4 T^*M = \oplus_{i+j=4} \wedge^i T^*\HH^7 \wedge \wedge^{j} T^*\bbS^4.$$ 
For each element of the form $u\wedge v$ with $$u \in \Gamma({}^e\wedge^*T^*\HH^7), v\in \Gamma(\wedge^*T^*\bbS^4)$$ the projection operator $\pi_\lambda$ maps it to $u\wedge \pi_{\lambda}v$, which by linearity extends to the whole bundle $K$. 
\end{proof}

We denote the projection by $\pi_{\lambda}$ on the sections of the bundle $K$ as the linear extension of the eigenvalue projection on $\bbS^{4}$. Note here we have a collection of eigenvalues on both functions and forms, specifically we have:
\begin{itemize}
\item on functions: $4k(k+3), k\geq 0$;
\item on closed 1-forms: $4k(k+3),k\geq 1$;
\item on co-closed 1-forms: $4(k+1)(k+2), k\geq 0$.
\end{itemize}

It follows from Lemma~\ref{Kdecomp} that the operator decomposes to a sum of infinitely many operators, each acting on a subbundle.
\begin{lemma}
The operator $dQ$ preserves the eigenspaces of $\bbS^4$, and decomposes as
$$
dQ=\sum_{\lambda\geq 0} dQ^{\lambda}:=\sum_{\lambda} \pi_\lambda \circ dQ \circ \pi_{\lambda}
$$
\end{lemma}
\begin{proof}
We only need to show that that Hodge laplacian $\Lap$ commutes with the linearized operator $dQ$. Since the linearized operator is composed from $\Lap^{hodge}$, $\Lap^{rough}$ (which are related by Bochner formula), Hodge $*$ operator, differential, and scalar operator, all of which commute with $\Lap$, $dQ$ therefore commutes with the eigenvalue projections.
\end{proof}

Here we write out the equations $dQ^{\lambda}u=0$ explicitly. For any element $(k,H)\in \Gamma(K)$ where $k$ is a symmetric two tensor and $H$ is a 4-form, the action of $dQ^{\lambda}$ is listed below with respect to the decomposition of $k$ and $H$. In particular, $\hat k_{ij}, \hat k_{IJ}$ are the trace free parts of $k$, and $Tr_{\SS^{4}}(k), Tr_{\HH^{7}}(k)$ are the traces of $k$; $k_{(1,1)}^{cl}$ and $k_{(1,1)}^{cc}$ correspond to cross terms $k_{Ij}$, and $H_{(i,j)}^{cl}$ and $H_{(i,j)}^{cc}$ are the $(i,j)$ component of $H$ (see discussion before Definition~\ref{d:expansion}). And $(*)^{cl}$ and $(*)^{cc}$ denote the $\SS^{4}$ projection to closed or coclosed forms, which also determines the choice of eigenvalue $\lambda$. 
\begin{equation}\label{e:318}
(\lambda+\Lap_{\HH^7}^{rough}-2)\hat k_{IJ}=0.
\end{equation}
\begin{equation}
(\lambda+\Lap_{\HH^7} +8)\hat k_{ij}=0.
\end{equation}

\begin{equation}\label{e:320}
\left\{
\begin{aligned}
6d_{\SS^{4}}*_{\SS^{4}}\Delta_{\HH^{7}}k_{(1,1)}^{cl}+3\lambda(Tr_{\HH^{7}}(k)-Tr_{\SS^{4}}(k))+*_{\SS^{4}}\lambda H_{(0,4)}^{cl}+d_{\SS^{4}}\Delta_{\HH^{7}}H_{(1,3)}^{cc}=0\\
\lambda k_{(1,1)}^{cl}+\Delta_{\HH^7}k_{(1,1)}^{cl}+12k_{(1,1)}^{cl}-6*_{\SS^4}H_{(1,3)}^{cc}=0\\
\lambda H_{(1,3)}^{cc}+d_{\HH^{7}}\delta_{\SS^{4}}H_{(0,4)}^{cl}=0\\
\lambda Tr_{\SS^{4}}(k)+\Lap_{\HH^7} Tr_{\SS^{4}}(k)+72Tr_{\SS^{4}}(k)-32*_{\SS^4}H_{0,4}^{cl} =0\\
\lambda Tr_{\HH^{7}}(k) +\Lap_{\HH^7}Tr_{\HH^{7}}(k) +12Tr_{\HH^{7}}(k) +28*_{\SS^4}H_{0,4}^{cl}-12Tr_{\SS^{4}}(k)=0
\end{aligned}
\right.
\end{equation}

\begin{equation}
\left\{
\begin{aligned}
-\lambda \delta_{\SS^{4}}H_{(3,1)}^{cl}-\Delta_{\HH^{7}}\delta_{\SS^{4}}H_{(3,1)}^{cl}+6*_{\HH^{7}}d_{\HH^{7}}\delta_{\SS^{4}}H_{(3,1)}^{cl}=0,\\
\lambda H_{(4,0)}^{cc}+d_{\HH^{7}}\delta_{\SS^{4}}H_{(3,1)}^{cl}=0.
\end{aligned}	
\right.
\end{equation}

\begin{equation}
\left\{
\begin{aligned}
6\Delta_{\HH^7} k_{(1,1)}^{cc}+\Delta_{\HH^7}*_{\SS^{4}}H_{(1,3)}^{cl}=0,\\
%\lambda *_{\HH^{7}} H_{(1,3)}^{cl}+\delta_{\HH^{7}}d_{\SS^{4}}H_{(2,2)}^{cc}+6\lambda k_{(1,1)}^{cc}=0\\
d_{\HH^7} H_{(1,3)}^{cl}+d_{\SS^4}H_{(2,2)}^{cc}=0,\\
\lambda k_{(1,1)}^{cc}+\Delta_{\HH^7}
k_{(1,1)}^{cc}+12k_{(1,1)}^{cc}-*_{\SS^4}H_{(1,3)}^{cl}=0.
\end{aligned}
\right.
\end{equation}

\begin{equation}\label{e:323}
\left\{
\begin{aligned}
\Delta_{\HH^{7}}\delta_{\SS^{4}}H_{(2,2)}^{cl}+\lambda\delta_{\SS^{4}}H_{(2,2)}^{cl}=0,\\
\lambda H_{(3,1)}^{cl}+d_{\HH^{7}}\delta_{\SS^{4}}H_{(2,2)}^{cl}=0.
\end{aligned}
\right.
\end{equation}
Notice that the leading order part of $dQ^{\lambda}$ is always the Hodge Laplacian $\Delta_{\HH^{7}}$ (or rough Laplacian only in the equation for two tensors on $\HH^{7}$). Moreover, in each of the system~\eqref{e:320}--\eqref{e:323}, some of the equations are essentially algebraic, which reduces the system to a smaller square one with diagonal term being $\Delta_{\HH^{7}}  \operatorname{Id}$.

Next we compute the indicial roots and kernels for the linearized operator  as an edge differential operator. Recall that $\partial M$ is the total space of fibration over $Y=\partial \BB^7$. 
\begin{definition}[Indicial operator]
Let $L:\Gamma(E_1)\rightarrow \Gamma(E_2)$ be an edge operator between two vector bundles over $M$. For any boundary point $p\in Y$, and $s\in \CC$, the indicial operator of $L$ at point $p$ is defined as
$$
I_p[L](s):\Gamma(E_1|_{\pi^{-1}(p)}) \rightarrow \Gamma(E_2|_{\pi^{-1}(p)})
$$
$$
(I_p[L](s))v=x^{-s}L(x^s \tilde v)|_{\pi^{-1}(p)}
$$
where $\tilde v$ is an extension of $v$ to a neighborhood of $\pi^{-1}(p)$. The indicial roots of $L$ at point $p$ are those $s\in \CC$ such that $I_p[L](s)$ has a nontrivial kernel, and the corresponding kernels are called indicial kernels.
\end{definition}
\begin{remark}
In the conformally compact case, the indicial operator is a bundle map from $E_1|_p$ to $E_2|_p$ (which is simpler than a partial differential operator as in the general edge case). In our case the indicial operator $I_{p}(L)(s)$ is an elliptic operator on a compact manifold $\SS^{4}$, hence the kernel of $I_{p}(L)(s)$ is independent of the choice of domain.   Moreover, since we have an $SO(7)$ symmetry for the operator, the indicial roots are invariant on $\bbS^6$. 
\end{remark}

\begin{proposition}\label{indicialroots}
The indicial roots of operator $dQ$ are symmetric around $\Re z=3$, with three special pairs of roots
$$
\theta_1^{\pm} = 3\pm 6i,\ \theta_2^{\pm}=3 \pm i \sqrt{21116145}/1655,\ \theta_3=3 \pm i3 \sqrt{582842}/20098.
$$
and all other roots lying in $\{\|\Re z-3\|\geq 1\}$.
\end{proposition}
\begin{proof}
With the harmonic decomposition on sphere $\bbS^4$, the linearized operator $dQ$ is block-diagonalized and we compute the indicial roots for the linear system $dQ$ in Appendix~\ref{indicial}. We summarize the results below and Figure~\ref{figure:indicial} is an illustration of the indicial roots distribution.
The indicial roots fall into the following three categories:

\begin{figure}
    \centering
    \includegraphics[width = 0.7\textwidth]{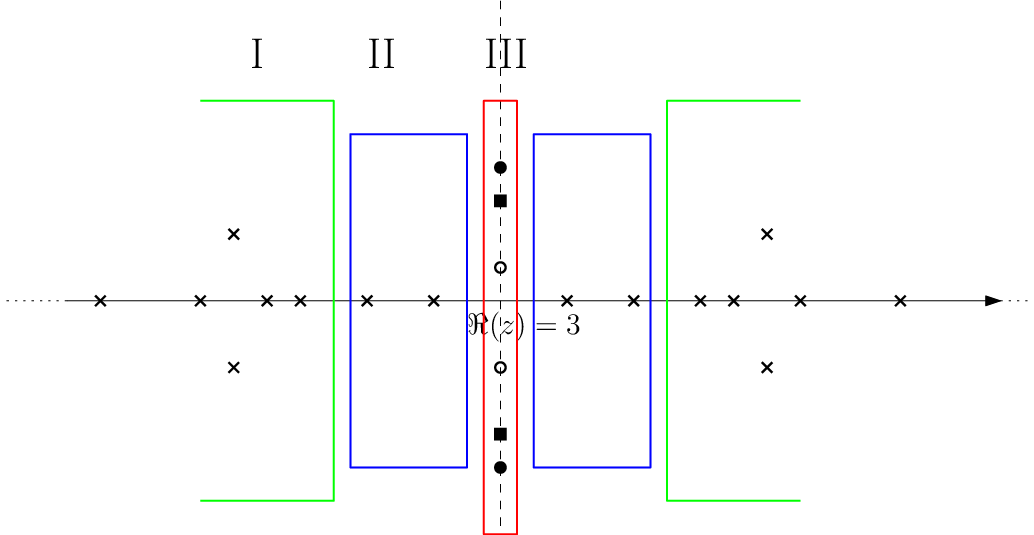}
    \caption{Indicial roots of the linearized supergravity operator on $\mathbb{C}$}\label{figure:indicial}
  \end{figure}

\begin{enumerate}
\item The roots corresponding to harmonic forms:
\begin{enumerate}
\item The equation for trace-free 2-tensors on $\HH^7$ arising from the first component
of~\eqref{eq:DQ} is
$$
(\Lap_{\SS^{4}}+\Lap_{\HH^{7}}-2)\hat k_{IJ}=0,
$$
and the corresponding indicial equation is
$$
(-s^2+6s)\hat k_{IJ}=0.
$$
We have indicial roots
$$
S_1^{+}=0, S_1^-=6.
$$
This corresponds to the perturbation of the hyperbolic metric to a Poincar\'{e}--Einstein metric.

\item  The equation for trace-free 2-tensors on $\SS^4$ is
$$
\Lap_{\SS^{4}}^{rough} \hat k_{ij}+\Lap_{\HH^7} \hat k_{ij}+8\hat k_{ij}=0
$$
where indicial equation is
$$
(-s^2+6s+8)\hat k_{ij}=0,
$$
and the indicial roots are
$$
S_2^{\pm}=3\pm \sqrt{17}.
$$
\item Equations for $H_{(4,0)}$:
\begin{equation}\label{V1}
\begin{array}{l}
d_{\HH^7}*H_{(4,0)}+W\wedge H_{(4,0)}=0\\
d_{\HH^7} H_{(4,0)}=0
\end{array}
\end{equation}
where the indicial equation is
$$
-(s-3)(*_{\SS^{6}} N)\wedge dx/x-6 dx/x\wedge N=0,
$$
with indicial roots
$$
\theta_1^{\pm}=3\pm 6i.
$$
This corresponds to a perturbation of the 4-form on hyperbolic space.
\end{enumerate}

\item The roots corresponding to functions / closed 1-forms / coclosed 3-forms / closed 4-forms
 
\begin{enumerate}
\item The equations for $7\sigma=Tr_{\HH^7}(k), 4\tau=Tr_{\SS^4}(k), k_{(1,1)}, H_{(1,3)}, H_{(0,4)}$ are
\begin{equation}\label{V23}
\begin{array}{l}
6d_{\HH^7} *_{\HH^7} k_{(1,1)}^{cl}+d_{\SS^4} (3Tr_{H}(k)-3Tr_{\SS^4}(k))\wedge {}^7V+d_{\SS^4}*H_{(0,4)}^{cl}+d_{\HH^7}*H_{(1,3)}^{cc}=0 \\
d_{\HH^7} H_{(0,4)}^{cl}+d_{\SS^4}H_{(1,3)}^{cc}=0\\
d_{\HH^7} H_{(1,3)}^{cc}=0\\
\Delta_{\SS^4} k_{(1,1)}^{cl}+\Delta_{\HH^7}
k_{(1,1)}^{cl}+12k_{(1,1)}^{cl}-6*_{\SS^4}H_{(1,3)}^{cc}=0\\
\Lap_{\SS^4} \tau+\Lap_{\HH^7} \tau+72\tau-8*_{\SS^4}H_{0,4}^{cl} =0\\
\Lap_{\SS^4} \sigma +\Lap_{\HH^7} \sigma +12\sigma +4*_{\SS^4}H_{0,4}^{cl}-48 \tau=0
\end{array}
\end{equation}
The indicial equations are
\begin{multline}
\lambda^4 - 4 S^2 \lambda^3 + 24 S*\lambda^3 - 90 \lambda^3 + 6 S^4 \lambda^2 - 72 S^3 \lambda^2 \\
+  342 S^2 \lambda^2 - 756 S*\lambda^2 + 1152 \lambda^2 - 4 S^6 \lambda + 72 S^5 \lambda - 
  414 S^4 \lambda \\
  + 648 S^3 \lambda + 1152 S^2 \lambda - 3024 S*\lambda + 10368 \lambda \\+ S^8 - 
  24 S^7 + 162 S^6 + 108 S^5 - 6192 S^4 + 31536 S^3 - 33696 S^2 - 
  155520 S = 0
\end{multline}
When $\lambda=16$
%, it gives
%$$
%157696 - 299136 S + 55904 S^2 + 23472 S^3 - 11280 S^4 + 1260 S^5 + 
 % 98 S^6 - 24 S^7 + S^8 =0
%$$
there is a pair of roots with real part 3
% that is approximately $3\pm i2.78..$
\begin{equation}
s=\theta_2^{\pm}=3 \pm i \sqrt{21116145}/1655
\end{equation}
and when 
$
\lambda=40
$ 
%the equation is
%$$
%-942080 + 49920 S + 303584 S^2 - 57744 S^3 - 13152 S^4 + 2988 S^5 + 
 % 2 S^6 - 24 S^7 + S^8 = 0
%$$
there is a pair of roots with real part 3
% and with an approximation $3\pm i0.11$
$$
\theta_3^{\pm}=3 \pm i3 \sqrt{582842}/20098
$$
And here the five variables are related by 
$$
H_{(0,4)}^{cl}=d_{\SS^4} *_{\SS^4} d_{\SS^4} \xi, H_{(1,3)}^{cc}=-d_{\HH^7} *_{\SS^4} d_{\SS^4} \xi, k_{(1,1)}^{cl}=-d_{\SS^4} \delta_{\HH^7} \xi, 4\sigma= 7\tau=\xi
$$
where 
$$
\xi \in \delta_{\SS^4} \wedge^{cl}_{16} \bbS^4 
$$
similarly we have another indicial kernel corresponding to $\theta^{\pm}_3$ with $
\xi \in \delta_{\SS^4} \wedge^{\cl}_{40} \bbS^4 
$.

\item The equations for $H_{(3,1)}, H_{(4,0)}$ are
$$
\begin{array}{l}
d_{\SS^4} * H_{(3,1)}^{cl}+d_{\HH^7}*H_{(4,0)}^{cc}+6{}^4V \wedge H_{(4,0)}^{cc}=0\\
d_{\HH^7}H_{(3,1)}^{cl}+d_{\SS^4}H_{(4,0)}^{cc}=0
\end{array}
$$
where the indicial equations are
$$
(s-3)^2\pm 6i(s-3)-16=0
$$
with indicial roots
$$
S_3^{\pm}=3\pm \sqrt{7}\pm 3i.
$$
\end{enumerate}
\item The roots corresponding to coclosed 1-forms / closed 2-forms / coclosed 2-forms / closed 3-forms
\begin{enumerate} 
\item The equations for $k_{(1,1)}, H_{(1,3)}, H_{(2,2)}$ are
$$
\begin{array}{l}
6d_{\HH^7} *_{\HH^7}k_{(1,1)}^{cc}+d_{\HH^7}*H_{(1,3)}^{cl}=0\\
d_{\SS^4} *H_{(1,3)}^{cl}+d_{\HH^7}*H_{(2,2)}^{cc}+6d_{\SS^4} *_{\HH^7} k_{(1,1)}^{cc}=0\\
d_{\HH^7} H_{(1,3)}^{cl}+d_{\SS^4}H_{(2,2)}^{cc}=0\\
\frac{1}{2}\Delta_{\SS^4} k_{(1,1)}^{cc}+\frac{1}{2}\Delta_{\HH^7}
k_{(1,1)}^{cc}+6k_{(1,1)}^{cc}-\frac{1}{2}*_{\SS^4}H_{(1,3)}^{cl}=0
\end{array}
$$
The indicial equation is
$$
\lambda^2-(36+(s-1)(s-5)+s^2-6s-1)\lambda -(s-1)(s-5)(-s^2+6s+1)=0.
$$
With the smallest eigenvalue for coclosed 1-forms being $\lambda=24$, the indicial roots are
$$
S_4^{\pm}=3\pm \sqrt{\pm 3\sqrt{97}+31}.
$$

\item 
The equations for $H_{(2,2)}, H_{(3,1)}$ are
$$
\begin{array}{l}
d_{\SS^4}*H_{(2,2)}^{cl}+d_{\HH^7}*H_{(3,1)}^{cc}=0\\
d_{\HH^7}H_{(2,2)}^{cl}+d_{\SS^4} H_{(3,1)}^{cc}=0
\end{array}
$$
The indicial equations are 
$$
(\Lap_{\SS^4}^{Hodge} -(2-s)(4-s))H_{(3,1)}=0,
$$
and for $\lambda=24$ we have
$$
S_5^{\pm}=3\pm \sqrt{17}.
$$
\end{enumerate}
\end{enumerate} 

\end{proof}

\section{Fredholm property of the linearized operator}\label{Fredholm}
Once we identify these indicial roots, we proceed using different strategies according to whether the indicial roots land on the $L^2$ line or not. We show that, for all sufficiently large indicial roots, the linearized operator after projection is invertible on suitable edge Sobolev spaces. This is done by using small edge calculus and $SO(5)$ invariance with respect to the boundary. We discuss the remaining finitely many indicial roots individually. For the three exceptional pairs we use the scattering theory to construct two generalized inverses, which encode the boundary data that parametrizes the kernel of the linearized operator.  

We then describe the kernel of this linearized operator in terms of the two generalized inverses, and a scattering matrix construction that gives the Poisson operator. Near any Poincar\'{e}--Einstein metric product that is close to the base metric $t$, a perturbation argument shows that the space given by the difference of the two generalized inverses is transversal to the range space of the linearized operator and therefore this space gives the kernel of the linearized operator, which later will provide the parametrization of the kernel for the nonlinear operator.

First of all, we define the domain for the linearized operator:
\begin{definition}\label{Def:domain}
Fix $\delta  \in (0,1)$, take any small $\epsilon>0$, and define the domain as 
$$
D_k(\delta)=\{u\in x^{-\delta}H^{2,k}_{e,b}(M;K):(dQ\pm i\epsilon) u \in x^{\delta}H^{0,k}_{e,b}(M;K)\}. 
$$
\end{definition}
\begin{remark}
The domain is well defined independent of any sufficiently small $\epsilon$. This follows from Proposition~\ref{LAP}.
\end{remark}

Using the projection operator $\pi_{\lambda}$ defined above, the domain can be decomposed: 
$$
D_k(\delta)=\oplus_{\lambda \in \Lambda} D_k(\lambda, \delta),
$$
where $\Lambda=\{4k(k+3), k\geq 0\} \cup \{4(k+1)(k+2), k\geq 0\}$ is the set of eigenvalues on the 4-sphere. We denote the operator acting on each subbundle as 
$$
dQ^{\lambda}:=\pi_{\lambda} \circ dQ \circ \pi_{\lambda}, \ dQ^{\lambda>M}=\sum_{\lambda>M} dQ^{\lambda}.
$$

This eigenvalue decomposition extends to the hybrid Sobolev spaces in Definition~\ref{Def:domain}. Consider the bundle $K$ over $M=\bbB^7\times\bbS^4$ which carries a unitary linear action of $\SO(5)$
covering the action on $\bbS^4$. There is an induced action of SO(5) on the space of smooth sections vanishing to all orders at the boundary denoted as $\dCI(\bbB^{7}\times\bbS^4;K),$ which extends to all the weighted hybrid Sobolev spaces $x^sH^{k,l}_{\text{e,b}}(\bbB^7\times\bbS^4;K)$ since the group acts by isometries. The linearized operator 
$dQ\in\Diff 2_{\text{e}}(\bbB^7\times\bbS^4;K)$ is an elliptic edge operator for
the product edge structure and we have shown that $dQ$ commutes with the induced
action of $\SO(5)$ on $\dCI(\bbB^7\times\bbS^4;K).$

The Sobolev spaces of sections of $K$ decompose according to the
irreducible representations of $\SO(5),$ all finite dimensional and forming
a discrete set. In particular these may be labelled by the eigenvalues, $\lambda,$ of
the Casimir operator for $\SO(5)$ with a finite dimensional span when $\lambda$
is bounded above. The $\SO(7,1)$ action on $\bbH^7$ commutes with the
$\SO(5)$ action on $K$ and acts transitively on $\bbH^7$, so the
multiplicity of the $\SO(5)$ representation does not vary over $\bbH^7.$
The individual representations of $\SO(5)$ in the decomposition of $K$
therefore form bundles over $\bbH^7.$ Therefore we have
\begin{lemma}
The group $\SO(5)$ acts on $x^{\delta}H_{e,b}^{s,k}(M;K)$ transitively, and those Sobolev spaces decomposes to Sobolev spaces of sections of subbundles on $\HH^7$.
\end{lemma}
\begin{proof}	
This follows from the transitivity of the $SO(7,1)$ and $SO(5)$ actions discussed above.
\end{proof}

We will separately discuss three parts. 
\subsection{Large eigenvalues}
One part is the infinite dimensional subspace corresponding to large eigenvalues 
$$
\oplus_{\lambda>\lambda_{0}} D_k(\lambda, \delta),
$$
on which the operators $dQ^{\lambda>\lambda_{0}}\pm i\epsilon$ are isomorphisms, and their inverses approach limits $R^{\lambda>\lambda_{0}}_\pm$ uniformly as $\epsilon$ goes to zero. This is shown by using ellipticity and a parametrix construction. 

\begin{proposition}\label{largeInd}
There is $\lambda_{0}>0$, such that for $\lambda>\lambda_{0}$ and any small $\epsilon>0$ the two operators 
$$dQ^{\lambda>\lambda_{0}} \pm i\epsilon: \oplus_{\lambda>\lambda_{0}} D_k(\lambda,\delta) \rightarrow \oplus_{\lambda>\lambda_{0}} \pi_{\lambda}x^{\delta}H^{0,k}_{e,b}(M;K)$$ 
are both isomorphism and their inverses have limits as $\epsilon \downarrow 0$.
\end{proposition}

To prove this proposition, we will bundle all the large eigenvalues together.
\begin{definition}
For $\lambda \in [0,\infty)$, let $\pi_{\geq \lambda}: K\rightarrow K$ be defined as the projection off the span of the eigenspaces of the Casimir operator for $\SO(5)$ with eigenvalues smaller than $\lambda$, i.e. $\pi_{\geq \lambda}:=\operatorname{Id}-\sum_{\lambda'<\lambda} \pi_{\lambda'}$.
\end{definition}

\begin{proposition}\label{large}
For any weight $s\in \mathbb{R}$ and any orders $p,k$, the bounded operator defined as
$$
dQ: x^sH^{p+2,k}_{e,b} (M;K) \rightarrow x^sH_{e,b}^{p,k} (M;K)
$$
is such that $\pi_{\geq \lambda_{0}} dQ$ is an isomorphism onto the range of $\pi_{\geq\lambda_{0}}$ for some $\lambda_{0}  \in [0,\infty)$ (depending on $s$ but not on $p$ and $k$). Moreover, the range of $\operatorname{Id}-\pi_{\geq \lambda_{0}}$ on $C^{\infty}(M;K)$ is the space $C^{\infty}(M;\oplus_{\lambda'<\lambda_{0}}\pi_{\lambda'}K)$ of sections of a smooth vector bundle over $M$ and $dQ$ restricts to it as an elliptic element of $\operatorname{Diff}_0^2(M; \oplus_{\lambda'<\lambda_{0}}\pi_{\lambda'}K)$.
\end{proposition}

The second part of the proposition is from the definitions of $0$ and edge operators and the fact that ellipticity in edge symbols implies ellipticity in zero symbols once fiber directions are removed. To prove the first part of the proposition, we first construct an SO(5)-invariant parametrix in the small edge calculus by finding a appropriate kernel on the edge stretched product space $M_e^2$ which is defined from $M^2$ by blowing up the fiber diagonal over the boundary of $M$~\cite{MR1133743}.

\begin{definition}
The edge stretched product $M_e^2$ for an edge manifold $M$ is defined as the blow up $[M^2;S]$ where $S$ consists of all fibres of the product fibration $\pi^2:(\partial M)^2\rightarrow Y^2$ which intersect the diagonal of $(\partial M)^2$. Let $\beta: M_{e}^{2}\rightarrow M^{2}$ be the blow down map, then we denote $\overline{\beta_{-1}(S)}$, the closure of the preimage of $S$ under the blow down map, as the front face. 
\end{definition}
Notice that from the definition of fiber diagonal, the blow up actually preserves the product structure of $\bbH^7 \times\bbS^4$, i.e. the fiber diagonal contained in $M^2$ is just the product $\Lap \times \bbS^4 \times \bbS^4$, and the manifold after the blow up is actually the product of two 4-spheres and the 0-double space $(\HH^{7})^{2}_{0}=[(\HH^7)^2,\partial \Delta_{0}]$ as defined in~\cite{MR916753}. 
\begin{lemma}
For $M=\HH^7\times \bbS^4$, the edge stretched product is actually a product: $M_e^2=[(\HH^7)^2,\partial \Delta_{0}]\times (\bbS^4)^2$.
\end{lemma}
As a result, the front face also has a product structure $\overline{\beta^{-1}(\partial \Lap_{0})} \times (\SS^{4})^{2}$.
The elliptic element $dQ \in \operatorname{Diff}_{e}^{2}(M;K)$ lifts to be transversely elliptic to the fiber diagonal down to the front face. Therefore we have a parametrix construction in the small edge calculus denoted as 
$$
\Psi_e^{*}(M;K)=\cup_{m\in \ZZ} \Psi_e^{m}(M;K).
$$ 
Here we recall that $\Psi_e^{m}(M;K)$ is defined to be the set of $m$-th order pseudodifferential operators whose kernel is a classical conormal distribution that vanishes to infinite order at both side boundary faces of $M_{e}^{2}$ and smooth across the front face,  see~\cite[Definition 3.3]{MR1133743}. And the remainder will be a smoothing operator contained in $\Psi^{-\infty}_e(M;K) := \cap \Psi_e^{m}(M;K)$.
\begin{lemma}
The SO(5)-invariant elliptic operator $dQ\in \operatorname{Diff}_e^2(M;K)$ has an $\SO(5)$-invariant parametrix $\tilde E$ in $\Psi_e^{-2}(M;K),$ such that
$$
\Id-dQ \circ \tilde E, \ \Id- \tilde E\circ dQ \in \Psi^{-\infty}_e(M;K)
$$
are also $\SO(5)$-invariant.
\end{lemma}
\begin{proof}
Any elliptic edge differential operator has a parametrix in the small edge
calculus, following Theorem 3.8 in~\cite{MR1133743}. The construction gives the kernel of
$E$ as a classical conormal distribution with respect to the `lifted
diagonal' of the stretched edge produce $M^2_e.$ Because of the product structure of $M_{e}^{2}$, in fact the action of $\SO(5)$ on the kernel $E,$ through the
product action on $M^2,$ lifts smoothly to $M^2_e$ and preserves the lifted
diagonal (which is the closure of the diagonal in the interior). So we may
average under the product action and define
$$
\tilde E=\int_{g \in SO(5)} g\cdot E. 
$$
Since $dQ$ is $SO(5)$ invariant, $\tilde E$ is also a parametrix,
$$
dQ \circ \tilde E =\Id+\tilde R,
$$
and the averaged remainder $\tilde R$ is also SO(5) invariant.
\end{proof}

As a consequence, now $\tilde E$ and $\tilde R$ both commute with the spherical eigenvalue projection $\pi_{\geq \lambda}$. Because of the special product structure, the edge small calculus can be characterized by 0-small calculus $\Psi^{*}_0(\HH^7)$ defined in~\cite{MR916753}, which is again filtered by order $\Psi^{*}_0(\HH^7)=\cup \Psi^{m}_0(\HH^7)$. We recall the definition here that $\Psi^{m}_0(\HH^7)$ is the class of pseudodifferential operators of order $m$ whose distribution is classical conormal on $(\HH^{7})^{2}_{0}$ vanishing to infinite order at both left and right faces. Similarly the smoothing operators are given by $\Psi^{-\infty}_0(\HH^7)$. Now the remainder $\tilde R$ can be characterized as:
\begin{lemma}
For any $\lambda$, the Schwartz kernel of $\pi_{\lambda}\tilde R$  is in $C^{\infty}((\bbS^4)^2, \Psi^{-\infty}_0(\HH^7)\otimes \operatorname{Hom}(\pi_{\lambda}K)) \subset C^{\infty}(M_e^2; K)$. In consequence it is a smooth map from $(\bbS^4)^2$ to bounded operators on $x^sH_0^p(\bbH^{7};\pi_{\lambda}K)$ for any $s,p$; and for any bounded range of $s$, the operator norm acting on $x^sH_0^p(\bbH^{7};\pi_{\lambda}K)$ is uniformly bounded by the Schwartz kernel norm of $\pi_{\lambda}\tilde R$ in $C^{k'}\left((\bbH^{7})_{0}^{2}; \operatorname{Hom}(\pi_{\lambda}K)\right)$ for some $k\in \NN$.
\end{lemma}
\begin{proof}
As an element in $\Psi _e^{-\infty}(M;K)$, the Schwartz kernel of $\tilde R$ is smooth on the double edge space $M_e^2$, with values in the bundle $\Hom(K) \otimes \cK$ where $\cK$ is the kernel density bundle. From the properties of the small calculus, the Schwartz kernel of $\tilde R$ vanishes to infinite order at the left and right boundary faces of $M_{e}^{2}$. Because $M_e^2$ has the product structure $(\bbH^{7})_0^2\times (\bbS^4)^2 $, the Schwartz kernel of $\pi_{\lambda}\tilde R$ is in $C^\infty\left((\bbS^4)^2, C^\infty((\bbH^{7})_0^2, \operatorname{Hom}(\pi_{\lambda}K)\otimes \cK)\right)$ where $C^{\infty}((\bbH^{7})_0^2,\operatorname{Hom(\pi_{\lambda}K) \otimes \cK } )$ gives the kernel of a $\Psi^{-\infty}_0(\bbH^{7};\pi_{\lambda}K)$ operator acting
on $\cK$.

From~\cite[Corollary 3.24]{MR1133743},  if $A$ is any element in the edge small calculus $\Psi_{e}^{-\infty}(M)$,  then $A:x^{\delta}H_{e}^{s}(M)\rightarrow x^{\delta}H_{e}^{s'}(M)$ is bounded for any $\delta, s,s'\in \RR$. As a special case of edge calculus, the same proof can be used to show that $\pi_{\lambda}\tilde R$ which is an element in $\Psi_0^{-\infty}(\bbH^{7};\pi_{\lambda}K)$ acts on $x^sH_0^p(\bbH^{7};\pi_{\lambda}K)$ as a bounded operator for any $s,p\in \RR$. Now consider the map $\phi$ from $\Psi_0^{-\infty}(\bbH^{7};\pi_{\lambda}K)$ to bounded operators on $x^sH_0^p(\bbH^{7};\pi_{\lambda}K)$. The space of $\Psi_0^{-\infty}(\bbH^{7};\pi_{\lambda}K)$ operators corresponds to Schwartz kernels smooth on the double space $(\HH^{7})^{2}_{0}$ and vanishing to infinite order at the left and right boundaries. This space is a Fr\'echet space with the usual $C^{\infty}-$topology on the double space $(\HH^{7})^{2}_{0}$~\cite[Remark 2.2.2(b)]{Lauter}
with semi-norms given by $C^{k}\left((\bbH^{7})_{0}^{2}; \operatorname{Hom}(\pi_{\lambda}K)\right)$ bounds of Schwartz kernels of such operators.
Since $\phi$ is a continuous map from a Fr\'echet space to a normed space, the norm is bounded by some norm on $\Psi_0^{-\infty}(\bbH^{7};\pi_{\lambda}K)$, i.e. the operator norm of $\pi_{\lambda}\tilde R$ on $x^s H_0^p(\bbH^{7};\pi_{\lambda}K)$ is bounded by $C(s)\| \pi_{\lambda}\tilde R\|_{C^{k}\left((\bbH^{7})_{0}^{2}; \operatorname{Hom}(\pi_{\lambda}K)\right)}$ where $C(s)$ is a constant only depending on $s$. Therefore for any bounded interval $s \in [-S,S]$, the bound of $\| \pi_{\lambda}\tilde R\|_{B(x^sH_0^p(\bbH^{7};K))}$ is uniform. 
\end{proof}

We can use the following interpolation result to show that $\pi_{\lambda}\tilde R$ rapidly decays as $\lambda$ tends to infinity.

\begin{lemma}\label{l:l2}
$
x^sH^{p,k}_{e,b}(M;K)\subset L^2(\bbS^4;x^sH^{p+k}_{0}(\bbH^{7};K))\cap H^{p+k}(\bbS^4;x^sL^2_{0}(\bbH^{7};K)).
$
%$
%x^sH^{p,k}_{e,b}(M;K)=L^2(\bbS^4;x^sH^{p+k}_{0}(\bbH^{7};K))\cap H^{p+k}(\bbS^4;x^sL^2(\bbH^{7};K)).
%$
\end{lemma}
\begin{proof}
We only prove the case $s=0$ since the weight on the boundary defining function $x$ transfer to the 0-Sobolev spaces on $\bbH^{7}$ directly. The space $L^2(\bbS^4;H^{p+k}_{0}(\bbH^{7};K))$ gives $p+k$ order of edge regularity on the $\bbH^{7}$ direction while the latter space $H^{p+k}(\bbS^4;L^2(\bbH^{7};K))$ gives $p+k$ order of regularity in the $\bbS^{4}$ direction. Together they give $p+k$ edge regularity on $M$. Since we have the inclusion $\cV_{b}\supset \cV_{e}$, the inclusion in the statement follows from $H^{p,k}_{e,b}(M;K) \subset H^{p+k}_{e}(M;K)$.
%On the other hand, if a function $f$ is in $L^2(\bbS^{4};H^p_0(\bbH^{7};K))\cap H^{p+k}(\bbS^4;x^sL^2(\bbH^{7};K))$, that is $V_{1}\dots V_{p}(f) \in L^2_{e}(M;K)$ for any $V_{i} \in \cV_{e}$. Then applying an elliptic k-th order differential b-operator to $f$ we obtain an element in $H^p(\bbS^4;L^2(\bbH^{7};K))$, therefore by elliptic regularity $f \in H_{e,b}^{p,k}(M;K)$. 
\end{proof}

\begin{lemma}
As $\lambda$ tends to infinity, the bounded operators $\pi_{\geq \lambda} \tilde R$ decay in operator norm on any Sobolev space $x^sH^{p,k}_{e,b}(M;K)$, i.e.
$$
\lim_{\lambda \rightarrow \infty} \|\pi_{\geq \lambda}\tilde R\|_{x^sH^{p,k}_{e,b}(\bbH^{7}\times \bbS^{4};K)\rightarrow x^sH^{p,k}_{e,b}(\bbH^{7}\times \bbS^{4};K)}=0.
$$
\end{lemma}
\begin{proof}
%We first need to show the rapid convergence of the Schwartz kernel as a smooth map into the bounded operators on the zero space, i.e. 
Using Plancherel it follows that the Schwartz kernel of $\pi_{\geq \lambda}\tilde R$ rapidly converges to 0 in $H^{p+k}((\bbS^4)^2, B(x^sL_0^2(\bbH^{7};K)))$ and $L^2((\bbS^4)^2, B(x^sH_0^{p+k}(\bbH^{7};K)))$.
Then we obtain $\|\pi_{\geq \lambda}\tilde R\|\rightarrow 0$ as bounded operators on $x^sH^{p,k}_{e,b}(M;K)$ by Lemma~\ref{l:l2}.
\end{proof}

As a consequence, for any fixed $s,k,l$, there is a $\lambda_0$ such that $\|\pi_{\geq \lambda_0}\tilde R\|_{x^sH^{k,l}_{e,b}(M;K)} \leq \frac{1}{2}$, and this $\lambda_0$ only depends on some $C^k$ norm of the Schwartz kernel on the double space. In the case that $\pi_{\geq \lambda_0}\tilde R$ is small, we get that $\pi_{\geq \lambda_0}dQ \pi_{\geq \lambda_0} \tilde E$ is a perturbation of the identity, which is therefore an isomorphism, that is,
\begin{lemma}
For any $s,k,l$, there is a $\lambda_0$ depending only on s, such that
$$
\pi_{\geq \lambda_0} dQ \pi_{\geq \lambda_0} \tilde E =Id_{\pi_{\geq \lambda_0} K} +\pi_{\geq \lambda_0} \tilde R
$$
where the right hand side is an isomorphism from $x^sH^{k,l}_{e,b}(M;K)$ to itself.
\end{lemma}
\begin{proof}
The norm of the operator on the right hand side acting on $x^sH^{k,l}_{e,b}(M;K)$ is bounded away from 0.
\end{proof}

The same argument applies to $\tilde E\circ \pi_{\geq \lambda_0} dQ \pi_{\geq \lambda_0} $. Then from the above lemma, we get that $\pi_{\geq \lambda_0} dQ$ is an isomorphism mapping from $\pi_{\geq \lambda_0}x^sH^{k+2,l}_{e,b}(M;K)$ to $\pi_{\geq \lambda_0}x^sH^{k,l}_{e,b}(M;K)$, proving the first part of Proposition~\ref{large}. And this implies Proposition~\ref{largeInd}:
\begin{proof}[Proof of Proposition \ref{largeInd}]
From Proposition~\ref{large}, there is a constant $C$ such that 
$$
C^{-1}\|u\|_{x^{\delta}H_{e,b}^{2,k}(M;K)}\leq \|dQ^{\lambda>\lambda_{0}}(u)\|_{x^{\delta}H_{e,b}^{0,k}(M;K)}\leq C\|u\|_{x^{\delta}H_{e,b}^{2,k}(M;K)}.
$$
With any sufficiently small $\epsilon$, from the triangle inequality we have, for another constant $\tilde C$ 
$$
\tilde C^{-1}\|u\|_{x^{\delta}H_{e,b}^{2,k}(M;K)}\leq \|(dQ^{\lambda>\lambda_{0}}\pm i\epsilon)u\|_{x^{\delta}H_{e,b}^{0,k}(M;K)}\leq \tilde C\|u\|_{x^{\delta}H_{e,b}^{2,k}(M;K)}.
$$
Therefore $(dQ^{\lambda>\lambda_{0}}\pm i\epsilon)^{-1}$ converge in the operator norm of $x^{\delta}H_{e,b}^{0,k}(M;K) \rightarrow x^{\delta}H_{e,b}^{2,k}(M;K)$ as $\epsilon \rightarrow 0^{+}$.

\end{proof}

\subsection{Individual eigenvalues with $\lambda \neq 0, 16, 40$}
Now we consider those eigenvalues smaller than $\lambda_0$. Consider the projected operator
$$
dQ^{\lambda}:D_{k}(\lambda,\delta) \rightarrow x^{\delta}H^{0,k}_{e,b}(M;\pi_{\lambda}K),
$$
which is viewed as a 0-problem on $\HH^7$ (each tensored with fixed eigenforms on $\bbS^4$).

For the purpose of simplicity, we denote the set of special indicial roots as
$$
\Lambda=\{0,16,40\}.
$$
From Proposition~\ref{indicialroots}, except for $\lambda\in \Lambda$,  the indicial roots of $dQ^{\lambda}$ are contained in the range $(-\infty, 3-\bar\delta] \cup [3+\bar\delta, \infty)$ with $\bar \delta=1$.
Moreover the pairs of indicial roots separate further as $\lambda$ becomes larger. With this information, we show that,
\begin{proposition}\label{small}
For $\lambda\notin\Lambda$, $dQ^\lambda: \pi_\lambda x^\delta H_{e,b}^{s,k}(M;K) \rightarrow \pi_{\lambda} x^\delta H_{e,b}^{s-2,k}(M;K)$ is Fredholm for any $|\delta|<\bar\delta$. Moreover, when $\delta>0$, this map is injective; when $\delta<0$, it is surjective.
\end{proposition}

The idea of the proof essentially follows the proof of the proposition below:
\begin{proposition}[Theorem 6.1 from~\cite{MR1133743}]
Suppose $L\in \operatorname{Diff}^{m}_{e}(M)$ is elliptic and satisfies 
\begin{enumerate}
\item constant indicial roots over the boundary;
\item unique continuation property of $N(L)$;
\end{enumerate}
and the weight $\delta$ satisfies $|\delta|<\bar \delta$ where $\bar\delta$ is the indicial radius, then $L: x^{\delta}H^{l+m}_{e}(M)\rightarrow x^{\delta}H^{l}_{e}(M)$
is an isomorphism.
\end{proposition}
\begin{remark}
For property (1), it can be seen from the previous computation of indicial roots or just by commuting with the group action of $SO(7)$. For a fixed $\lambda$, the operator $dQ^{\lambda}$ acts on sections of vector bundles on $\HH^{7}$ as a 0-operator in the sense of~\cite{MR916753}. Since the operator commutes with the group $SO(7)$, we can decompose the operator $dQ^{\lambda}$ further using spherical harmonics which will give a system of ODEs each acting on one of the spherical modes. The unique continuation property of ordinary differential operators can then be applied to show that, if there is a solution that vanishes to infinite order at the boundary, then the solution must vanish everywhere. Combine all the spherical modes, we see that  $dQ^{\lambda}$ satisfies (2).  Another way to show this property is to adapt the proof from~\cite{MazUni} on unique continuation on scalar Laplacian on the asymptotic hyperbolic manifolds. Since the leading order part of each $dQ^{\lambda}$ is a scalar Laplacian and the rest is at most 1st order differential, one can prove the unique continuation by Carleman estimates where all the lower order terms are estimated just by  $|\nabla u|+|u|$.
\end{remark}

We first introduce two operators related to $dQ^{\lambda}$: the normal operator and reduced normal operator.

\begin{definition}[Normal operator]
For $L \in \operatorname{Diff}_e^*(M)$ the normal operator $N(L)$ is defined to be the restriction to the front face $B_{11}$ of the lift of $L$ to $M_e^2$. In terms of the local coordinate, if 
$$
L=\sum_{j+|\alpha|+|\beta|\leq m} a_{j,\alpha,\beta}(x,y,z)(x\partial_x)^j(x\partial_{y})^\alpha \partial_z^{\beta}
$$
then
$$
N(L)=\sum_{j+|\alpha|+|\beta|\leq m} a_{j,\alpha, \beta}(0, \tilde y, z)(s\partial_s)^j(s\partial_u)^\alpha \partial_z^\beta,
$$
where $s,u,\tilde x, \tilde y, z,\tilde z$ is the lifted coordinate system on $M_e^2$ covering $B_{11}$ such that
$$
s=\frac{x}{\tilde x}, \ u=\frac{y-\tilde y}{\tilde x}.
$$  
By the identification $\HH^{7}\simeq (0,\infty)_{s}\times \RR^{6}_{u}$, the normal operator $N(L)$  is naturally an operator on $\HH^{7}\times \SS^{4}$. 
\end{definition}
Since $dQ^{\lambda}$ does not have a $\partial_{z}$ component, $N(dQ^{\lambda})$ acts on the half space $\RR^{+} \times \RR^{6} \ni (s,u)$ which is the tangent space of $\HH^{7}$ for each fixed boundary point $\tilde y \in Y=\SS^{6}$. So it is the operator with coefficients frozen at the boundary.
This may be further reduced to be the reduced normal operator which is a family of differential b-operators. 
\begin{definition}[Reduced normal operator]
The reduced normal operator $N_0(L)$ is defined by conjugating $N(L)$ by the Fourier transform in the $\RR^6$ direction then rescaling. Specifically, if we denote $\eta$ the dual variable to $u$, and set $t=s|\eta|, \hat \eta=\frac{\eta}{|\eta|}$, then
$$
N_0(L)=\sum_{j+|\alpha|+|\beta|\leq m}a_{j,\alpha,\beta} (0,\tilde y, z)(t\partial_t)^j (it\hat \eta)^\alpha \partial_z^\beta, \quad t \in \RR_+, \quad  \hat \eta \in S_{\tilde y}^* Y.
$$
\end{definition}
Since the reduced normal operator of $dQ^\lambda$ is independent of the fiber variables $z$, it is for fixed $(\tilde y, \hat \eta)$ an ordinary differential operator on $\RR^{+}_{t}$ and has the following mapping property:

\begin{lemma}
For $\lambda\notin\Lambda$, given any fixed $(\tilde y,\hat \eta) \in \SS^{6}\times S_{\tilde y}^* \SS^{6}$ and $|\delta|<\bar \delta$, the reduced normal operator 
$$
N_0(dQ^\lambda): t^{\delta}H^{2}(\RR_{t}^{+})\rightarrow t^\delta L^2(\RR_t^+)
$$
is an isomorphism.
\end{lemma}
\begin{proof}
The proof is contained in Lemma 5.5--5.12 in~\cite{MR1133743}.
For each $\lambda \notin \Lambda$, $N_0(dQ^\lambda)$ is a regular singular second order ordinary differential operator, and has a pair of indicial roots $3\pm \delta_{\lambda}$ with $|\delta_{\lambda}|\geq \bar \delta$.  Near $t=0$ the operator is an ordinary b-operator controlled by the pair of indicial roots; near $t=+\infty$, the operator is of the form $t^{2}E+O(t)$ where the leading order term of $E$ is given by $\partial_{t}^{2}-\hat \eta^{2}$ hence elliptic in the sense of b-operators by considering the transformation $s=\frac{1}{t}$. This operator is of the ``Bessel type''~\cite[Definition 5.3]{MR1133743}. Following the proof there, parametrices can be constructed near the two ends, $t=0$ and $t=\infty$. Near $t=0$, only the b-structure matters and the parametrix $H_{0}^{0}$ is constructed in Theorem 4.4 of~\cite{MR1133743} and it is bounded between $t^{\delta}H_{b}^{\ell}(\bbR^{+})$ and $t^{\delta}H_{b}^{\ell+2}(\bbR^{+})$ for $|\delta|<\bar \delta$. Near $t=\infty$, the ``Bessel type'' structure of $N_{0}=a_{2,2}t^{2}\partial_{t}^{2}+a_{2,0}t^{2}+\sum_{j+m\leq 1} a_{j,m} t^{j}\partial_{t}^{m} $ is used. It is shown that the partial principle symbol is given by $\tilde \sigma(N_{0})=a_{2,2}t^{2}\tau^{2}+a_{2,0}t^{2}$ (where $\tau$ is dual to $\partial_{t}$) such that $\langle\tilde \sigma(N_{0})u,u\rangle \geq Ct^{2}(1+\tau^{2})\|u\|^{2}$ for sufficiently large $t$, hence a parametrix near $t=\infty$ can be constructed by integration $H_{0}^{\infty}(u):=\int e^{it\tau} \tilde \sigma(N_{0})^{-1} \hat u(\tau)d\tau$.
 The two parametrices are patched together using a cutoff functon $\phi(t)$ to give $H_{0}=\phi H_{0}^{0}+(1-\phi)H_{0}^{\infty}$ that gives the Fredholm property:
\begin{equation}
H_{0}\circ N_{0}=I-P_{01}, \ N_{0}\circ H_{0}=I-P_{02}
\end{equation}
such that $P_{01}$ and $P_{02}$ are compact. For an ordinary differential operator on $\RR^{+}$ acting on weighted space with weights between the two indicial roots, the operator is bijective.
\end{proof}

\begin{lemma}\label{l:NH}
For $\lambda \notin \Lambda$ and $ |\delta|<\bar \delta $, the normal operator $N(dQ^\lambda)$ is Fredholm on $x^{\delta}H_e^0(M;\pi_\lambda K)$. 
\end{lemma}
\begin{proof}%need work
This proof is contained in~\cite[Theorem 5.16]{MR1133743}. The reduced normal operator is obtained by Fourier transform and normalization of the operator $N(dQ^{\lambda})$, so we can reverse this process to do first rescaling then an inverse Fourier transform to get the parametrix for $N(dQ^{\lambda})$ from the parametrix $H_{0}$ above. Specifically, 
$$
\hat H(s,\tilde s,\eta )=H_{0}(s|\eta|, \tilde s |\eta|, \eta)|\eta|
$$
is a bounded operator from $s^{\delta}\hat H_{e}^{0}(M;K)$ to $s^{\delta}\hat H_{e}^{2}(M;K)$ where in the definition of $\hat H_{e}^{*}$ the differentiation $s\partial_{u}$ is replaced by multiplication of $s\eta$. And this gives $\hat P_{01}$ and $\hat P_{02}$ with the correct bounds. Then by doing an inverse Fourier transform
$$
N(H)(s,\tilde s, u, \tilde u)=\int e^{i(u-\tilde u)\eta}\hat H(s, \tilde s, \eta) d\eta
$$
we obtain the normal operator for the generalized inverse
$$
N(H): s^{\delta}H_{e}^{0}(M;\pi_{\lambda}K) \rightarrow s^{\delta}H_{e}^{2}(M;\pi_{\lambda}K)
$$
with corresponding compact errors $N(P_{0i}), \ i=1,2$.
\end{proof}
We may then use representation theory to show this operator is injective on any space contained in $L^{2}$:
\begin{lemma}\label{l:inj}
The kernel of the normal operator $N(dQ^\lambda)$ on $x^{\delta}H_{e,b}^{2,k}(M;\pi_\lambda K)$ is trivial for $\delta>0$.
\end{lemma}
\begin{proof}
From the expression of $dQ^{\lambda}$ in~\eqref{e:318}--\eqref{e:323}, we know that each operator $dQ^{\lambda}$ acting on some bundles on $\HH^{7}$ is given by leading order $\Delta_{\HH^{7}}\operatorname{Id}$ on a direct sum of forms, or $\Delta^{rough}_{\HH^{7}}$ on two-tensors. And off diagonal terms are lower order. So to study the normal operator, we only need to consider the injectivity and surjectivity of the operators $\Delta_{\HH^{7}}+L$ and $\Delta_{\HH^{7}}^{rough}+L$ with $L\in \RR$ out of the spectrum (which corresponds to the assumption that $\lambda\notin\Lambda$ hence the indicial roots are away from the $L^{2}$ cutoff line).

If the kernel of $N(dQ^{\lambda})$ inside $L^{2}$ is nontrivial, take such an element $u$, and based on the discussion above, $u$ lies in the eigenspace of $\Delta_{\HH^{7}}$ or $\Delta_{\HH^{7}}^{rough}$. Take the subspace formed by action of $SO(7,1)$ on $u$, which is a finite-dimensional $L^{2}$ eigenspace. However,  
there are no finite dimensional $L^2$ invariant subspaces of forms on $\HH^7$~\cite{MR961517};  and there are no $L^2$ eigentensors, 
from Delay's result~\cite{delay2002essential}. It follows that the kernel of $N(dQ^{\lambda})$ must be trivial.
 \end{proof}

As a result we have:
\begin{lemma}
For any $\bar \delta>\delta>0$, the normal operator $N(dQ^\lambda)$ is injective on $x^{\delta}H_{e,b}^{2, k}(M;\pi_\lambda K)$ and surjective on $x^{-\delta} H_{e,b}^{2,k}(M;\pi_\lambda K)$.
\end{lemma}
\begin{proof}
From the discussion above, the kernel of map $N(dQ^{\lambda})$ on $x^{\delta}H_{e,b}^{2, k}(M;\pi_\lambda K)$ is contained in the $L^2$ space of forms and tensors on $\HH^7$ such that for each fixed $\lambda$ it is contained in some eigenspace of $\Delta_{\HH^{7}}$ or $\Delta_{\HH^{7}}^{rough}$, which from the lemma above does not have any nontrivial elements. Therefore the map $N(dQ^\lambda)$ is injective. By considering $N(dQ^{\lambda})^{*}$ and duality with respect to the $L^{2}$ space, we can see that the operator $N(dQ^{\lambda})$ is surjective on the bigger space $x^{-\delta} H_{e,b}^{2,k}(M;\pi_\lambda K)$.
\end{proof}

We now return to the original operator $dQ^\lambda$ and show that it is Fredholm.

\begin{proof}[Proof of Proposition~\ref{small}]
From $N(H)$ constructed in Lemma~\ref{l:NH} which is defined on the front face of $M_{e}^{2}$, we extend it to $M_{e}^{2}$ and solve off the errors on the left boundary faces using the indicial operator. In this way we get a left parametrix $H_{L}$. Similarly one get a right parametrix $H_{R}$. This process is described in~\cite[Theorem 6.1]{MR1133743}. This gives two generalized inverses acting on edge Sobolev spaces. And by commuting any b-vector fields with the operator $dQ^{\lambda}$ and use the relation that $[\cV_{b}, \cV_{e}]\subset \cV_{b}$, we can see that the b-regularity is preserved by the generalized inverses. Hence
$$
H_{L/R}: x^{\delta}H_{e,b}^{0,k}(M;\pi_{\lambda}K)\rightarrow x^{\delta}H_{e,b}^{2,k}(M;\pi_{\lambda}K),
$$
which shows that $dQ^{\lambda}$ is Fredholm.

For a general kernel element of $dQ^{\lambda}$, we decompose it using the $SO(7,1)$ action, so it lies in the null space of the normal operator $N(dQ^\lambda)$ which is trivial by Lemma~\ref{l:inj}. Therefore the kernel is also trivial for the operator $dQ^\lambda$. The operator is therefore injective on the smaller space, and by duality surjective on the larger space.
\end{proof}

\subsection{Individual eigenvalues with $\lambda\in \Lambda$}
For those eigenvalues corresponding to indicial roots with real part equal to 3, we consider each subspace $\pi_\lambda x^{-\delta}H^{2,k}_{e,b}(M;K)$ separately. Restricted to these subspaces, the linearized operator is a 0-operator on hyperbolic space, of which the main part is the hyperbolic Laplacian $\Lap_\HH$. From Mazzeo--Melrose~\cite{MR916753} and Guillarmou~\cite{MR2153454}, the resolvent of $\Lap_\bbH-\lambda(6-\lambda)$, denoted as $R(\lambda)$, extends to a meromorphic family with finite order poles. Similarly, we want to show that $dQ$ has two generalized inverses $R_\pm$, which are the limits of the resolvent $(dQ \pm i\epsilon)^{-1}$ when $\epsilon\downarrow 0$ that extends to the spectrum.

First we show that the indicial roots become separated from the $L^2$ line by adding an imaginary perturbation.
\begin{lemma}\label{l:defo}
For $\lambda \in \Lambda$ and any $\epsilon>0$, the two indicial roots of the operator $dQ^{\lambda} \pm i\epsilon$ lie off the $\Re(s)=3$ line.
\end{lemma}
\begin{proof}
Suppose $s \in \CC$ is an indicial root for an operator $P$ on a point $p$ at the boundary, then we have $P(x^s)=O(x^{s+1})$ by definition. For $\epsilon\neq 0$, the following computation shows that $s$ is no longer an indicial root: $(P+i\epsilon)(x^s)=i\epsilon x^s + O(x^{s+1}) \neq O(x^{s+1})$. Instead, take the harmonic 4-form part which has indicial roots $3\pm 6i$ which in the indicial root computation is $dQ^{\lambda}(x^{3+\theta})=(\theta^2+36)x^{3+\theta} + O(x^{4})$, after perturbation it becomes 
$$
(dQ^{\lambda}+i\epsilon)(x^{3+\theta})=(\theta^2+i\epsilon+36)x^{3+\theta} + O(x^{4}) 
$$
so the indicial equation becomes $\theta^2=-i\epsilon-36$ which moves the two roots $3\pm \theta_{1}$ off the line of $Re(s)=3$. A similar argument  applies to other two pairs of roots $3 \pm \theta_{i}, \ i= 2, 3$.
\end{proof}

\begin{lemma}
For $\lambda\in \Lambda$ and any sufficiently small $\epsilon> 0$, the inverse $(dQ^\lambda \pm i\epsilon)^{-1}: x^\delta H_{e,b}^{2,k}(M;$ $\pi_\lambda K) \rightarrow  D_{k}(\lambda, \delta) \subset x^{-\delta}H_{e,b}^{0,k}(M;\pi_\lambda K) $ exists as a bounded operator. 
\end{lemma}
\begin{proof}
Using the indicial roots separation and same argument as in Proposition~\ref{small}, the operator $dQ^\lambda \pm i\epsilon$ is Fredholm on $x^\delta H_{e,b}^{2,k}(M;\pi_\lambda K)$, injective on the smaller space and surjective on the larger space. 
\end{proof}
Moreover with the limiting absorption principle below we show that the bound is uniform with respect to $\epsilon \downarrow 0$. Consider the reduced normal operator of $dQ^{\lambda}+i\epsilon$, which is a differential operator (parametrized by $y$ and $\epsilon$), is injective from $x^{\delta}H^2(\bbR^+) \rightarrow x^{\delta}L^2(\bbR^+)$ for any fixed $\delta>0$. This ODE operator may be extended holomorphically as $\epsilon$ approaches zero from above, 
and the solution of the ODE extends holomorphically as well. However after extending $\epsilon$ past zero, the smaller indicial root becomes the larger one, which is excluded from the solution. This is reflected in the resolvent $R_+^{\lambda}:=\lim_{\epsilon\downarrow 0} (dQ_\lambda+i\epsilon)^{-1}$ as the expansion of $R_{+}^{\lambda}u$ for $u\in x^{\delta}H_{e,b}^{0,k}(M;K)$ has only half of the indicial roots. And similarly for the other direction, the expansion of $R_{-}^{\lambda}u$ only has the other half of the indicial roots. 

To prove this limiting absorption principle, we use the specific structure of $dQ^{\lambda}$ as in~\eqref{e:318}-\eqref{e:323}, which is separated into two cases: $\lambda=0$ and $\lambda=16,40$. In the first case, the operator is exactly the Hodge Laplacian on $\HH^{7}$ acting on a 4-form, so we use the result from~\cite{Vasy} and the explicit construction of resolvents in~\cite{kantor2009eleven}. For the other case, the operator is a matrix system of which the diagonal is given by the Laplacian on $\HH^{7}$ acting on functions, and the rest is given by a constant matrix. Here we follow a similar strategy as in~\cite[Theorem 3.1.4]{DZ} and use spectral measures and contour deformation to show the bound. 

\begin{proposition}[Limiting absorption principle]\label{LAP}
For $0<\delta<\bar \delta$, $\lambda \in \Lambda$, and $\epsilon>0$ the operators $(dQ^\lambda \pm i\epsilon )^{-1}$ converges uniformly to bounded operators on weighted Sobolev spaces,
$$
\lim_{\epsilon \downarrow 0}\|(dQ^\lambda\pm i\epsilon)^{-1}-R^{\lambda}_\pm\|_{x^{\delta}H_{e,b}^{0,k}(M; \pi_\lambda K) \rightarrow  x^{-\delta} H^{2,k}_{e,b}(M; \pi_\lambda K)} =0. 
$$
\end{proposition}

\begin{proof}
We give a proof using different strategies for the following two cases: (1) $\lambda=0$; (2) $\lambda=16$ or 40.

\textbf{The $\lambda=0$ case} 

The nontrivial indicial root pair in this case comes from the system on $H_{(4,0)}$, and the equation is given by~\eqref{e:s1}.  Hence the operator $dQ^{\lambda}$ is given by $\Delta_{\HH^{7}}+36$ where the Hodge Laplacian acts on 4-forms on $\HH^{7}$. Then using the result from~\cite[Theorem 1.1]{Vasy} and the estimates (1.1) there, we have that the resolvents $(dQ^{\lambda}\pm i\epsilon)$ extend meromorphically to the real line. Moreover, from the explicit construction of $(\Delta_{\HH^{7}}+\mu)^{-1}$ in~\cite[Theorem 9.5]{kantor2009eleven} on the double space and the mapping property of the full 0-calculus from the decay rate on the boundary on $(\HH^{7})_{0}^{2}$ (using a similar proof as in the case of scalar Laplacian below), we can see that the resolvent is uniformly bounded.

\textbf{The $\lambda=16$ or $40$ case}

We discuss the $\lambda=16$ case, and the other $\lambda=40$ case is proved using the exactly same method. In this case, the system corresponds to $dQ^{\lambda}$ acting on $\Tr_{\SS^{4}}(k), \ \Tr_{\HH^{7}}(k),\ k_{(1,1)}, $ $H_{(0,4)}, \ H_{(1,3)}$. And after rearrangement it can be written as a square matrix acting on four functions on $\HH^{7}$ as in~\eqref{e:s2}. In this case, 
$$dQ^{\lambda}=(\Delta_{\HH^{7}}-9)\operatorname{Id} + A$$ where $A$ is a 4-by-4 constant matrix. From the indicial roots computation, we know that the matrix $A$ has one negative eigenvalue $-\mu_{0}^{2}=(\sigma_{2}^{\pm}-3)^{2}$, and other eigenvalues are all positive. In particular $A$ does not have 0 as an eigenvalue (otherwise there will be a pair of double indicial roots equal to 3). 

For the function Laplacian, we have the resolvent 
$$
R_{0}^{+}(s)=\left(\Delta_{\HH^{7}} -s(6-s) \right)^{-1}=\left(\Delta_{\HH^{7}} -9-\mu^{2} \right)^{-1}
,  \ i\mu=s-3,
$$
which maps from $L^{2}_{0}(\HH^{7})$ to $H^{2}_{0}(\HH^{7})$ when $\Re s>>1$ (or $\Im \mu>>1$) and is extended meromorphically to the line $\Re s=3$ (or the real line for $\mu$) to be a bounded operator from $x^{\delta}L^{2}_{0}(\HH^{7})$ to  $x^{-\delta}H^{2}_{0}(\HH^{7})$ for $\delta>0$. Similar we have the other resolvent $R_{0}^{-}$ extended from the other side. And we can prove the following uniformity result for $R_{0}^{\pm}(s)$:
\begin{equation}\label{e:R0}
\lim_{\epsilon\downarrow 0}\|(\Lap_{\bbH^{7}}-s(6-s)\pm i\epsilon)^{-1}-R_{0}^{\pm}(s)\|_{x^{\delta}L^{2}_{0}(\BB^{7})\rightarrow x^{-\delta}H^{2}_{0}(\BB^{7})}=0.
\end{equation}
The proof of this will be included in the end. Suppose we know this estimate, then by commuting with b-vector field we can get the estimate for $R_{0}^{\pm}(s)$ on the hybrid Sobolev space. 

Now we show that for the operator $(\Delta_{\HH^{7}}-9)\operatorname{Id} + A$ a similar estimate holds. Denote
$$
P=\Delta_{\HH^{7}}-9.
$$
Using spectral measure, we have the following expression for any smooth function $f$:
$$
f(P)  =\frac{1}{2\pi i}\int_{\RR} f(\mu^{2}) (P-\mu^{2})^{-1}d\mu.
$$
And we would like to apply it to two matrix-valued functions $f^{\pm}(P)=(P\otimes Id + A \pm i\epsilon \otimes Id)^{-1}=(dQ^{\lambda} \pm i\epsilon)^{-1}$. Take $f^{+}$ for example. When $\epsilon=0$, the integral above has two singularities on $\mathbb{R}$ at exactly $\pm \mu_{0}=-i(\theta_{2}^{\pm}-3)$ since this is where $f(\mu)=(\mu^{2}+A)^{-1}$ has a pole by the discussion of eigenvalues of $A$ above. And by Lemma~\ref{l:defo}, when $\epsilon>0$, these two singular points move away from $\pm\mu_{0}$, in particular, the one corresponding to $\mu_{0}$ moves into the upper half plane and the other one $-\mu_{0}$ moves into the lower half plane. So for a fixed small range of $\epsilon$, there is a curve $\gamma^{+}$ which is the real axis except near $\pm \mu_{0}$, and it goes above $\mu_{0}$ and below $-\mu_{0}$ so that it avoids all the singularities of $f^{+}$ for this fixed range of $\epsilon$ (see Figure 3.1 in~\cite{DZ} for illustration). Similarly we have an opposite contour $\gamma^{-}$ for $f^{-}$. 
Now we have the identity
$$
(dQ^{\lambda}\pm i\epsilon)^{-1}=\int_{\gamma^{\pm}} (\mu^{2}Id+A\pm i\epsilon Id)^{-1} R_{0}^{\pm}(\mu) d\mu.
$$
There is convergence at $\mu\rightarrow \pm \infty$ since away from the two singularities the integral can be completely deformed into $\Re\mu>>1$ or $-\Re \mu>>1$ where it has good invertibility. Now using the uniform estimate~\eqref{e:R0} of $R_{0}^{\pm}(\mu)$ as $\mu$ approaches the real line, we get the uniform estimate for $(dQ^{\lambda}\pm i\epsilon)^{-1}$ as $\epsilon\rightarrow 0$.
 
We now finish by proving~\eqref{e:R0}, which is essentially contained in~\cite{MR916753} with the description of the kernel of the resolvent $R(s)=(\Lap_{\bbH^{7}}-s(6-s))^{-1}$ for $s$ approaching $\Re s=3$ which corresponds to the continuous spectrum. Here we use proposition 6.2 in~\cite{MR916753} that the kernel of $R(s)$, which is a function of the hyperbolic distance $|(x,y)-(x',y')|$,
decays as $x^{\Re s}(x')^{\Re s}q(x,y,x',y')$ on the boundary of the double space and belongs to $\rho_{\ff}^{-2s}C^{\infty}((\HH^{7})_{0}^{2})$ where $\rho_{\ff}$ is a boundary defining function of the front face of $(\HH^{7})_{0}^{2}$. To show the bounds, we use the fact that $x^{3}L^{2}(\BB^{7};x^{-1}dxdy)=L^{2}(\BB^{7};x^{-7}dxdy)=L^{2}_{0}(\HH^{7})$. Therefore if we consider the kernel $\tilde R(s)$ on $L^{2}(\BB^{7}; x^{-1}dxdy)$, then $\tilde R(s,x,x')=x^{-3}R(s,x,x')(x')^{-3}$. For fixed $\delta>0$ and any $|\Re s-3|<<1$,  
$$
\begin{array}{l}
\sup_{x'\in \BB^{7}}\int_{x,y}x^{\delta}\tilde R(s,x,x',y,y')(x')^{\delta}x^{-1}dxdy<C, \\
\sup_{x\in \BB^{7}}\int_{x',y'}x^{\delta}\tilde R(s,x,x',y,y')(x')^{\delta}(x')^{-1}dx'dy'<C,
\end{array}
$$
where $C$ does not depend on $\epsilon$.
Then by Schur's lemma $\tilde R(s): x^{\delta}L^{2}(\BB^{7};x^{-1}dxdy) \rightarrow x^{-\delta}L^{2}(\BB^{7};x^{-1}dxdy)$ is bounded in the operator norm for $s$ approaching $s^{\pm}$. Transform back to the hyperbolic space, $R(s): x^{\delta}L^{2}_{0}(\HH^{7}) \rightarrow x^{-\delta}L^{2}_{0}(\HH^{7})$ is uniformly bounded. The bound on $x^{\delta}L^{2}_{0}(\HH^{7})\rightarrow x^{-\delta}H_{0}^{2}(\HH^{7})$ follows from ellipticity by commuting with the 0-elliptic operator $\Lap_{\HH^{7}}$.

\end{proof}

\subsection{Boundary data parametrization}
Combining the analysis for $\lambda$ off the $L^2$ line and on the $L^2$ line, we conclude:
\begin{proposition}
For $\delta \in (0,\bar \delta)$ and a product metric $h\times \frac{1}{4}g_{\bbS^{4}}$ on $M$ with $h$ being a Poincare--Einstein metric sufficiently close to hyperbolic metric, there are two generalized inverses $R_{\pm}: x^{\delta}H_{e,b}^{0,k}(M;K)$ $\rightarrow x^{-\delta} H_{e,b}^{2,k}(M;K)$ for operator $dQ$, such that
$$
dQ\circ R_+=Id, dQ\circ R_-=Id: x^{\delta}H_{e,b}^{0,k}(M;K) \rightarrow x^{\delta} H^{0,k}_{e,b}(M;K).
$$
\end{proposition}
\begin{proof}
When $h$ is the hyperbolic metric, we just need to combine result from Proposition~\ref{largeInd}, ~\ref{small}, and~\ref{LAP}. We will show that when the boundary conformal class of $h$, denoted as $\hat h$, is close to the standard spherical metric $\|\hat h-g_{\SS^{6}}\|_{H^{k}(\SS^{6})}<<1$ (therefore $h$ is a Poincare--Einstein metric close to $g_{\HH^{7}}$), the same result applies. 

For $dQ^{\lambda>\lambda_{0}}$, Proposition~\ref{largeInd} holds since $dQ$ is still an edge operator. And one can choose a boundary defining function on $\BB^{7}$ such that the normal operator $N(dQ^{\lambda})$ is the same as before, therefore the same analysis on $N(dQ^{\lambda})$ applies for both $\lambda$ off  and on the $L^{2}$ line. The injectivity of $dQ^{\lambda}$ for any individual $\lambda \notin\Lambda$ acting on $x^{\delta}H_{e,b}^{2, k}(M;\pi_\lambda K)$ is obtained by the same argument as in the case when $h$ is the exact hyperbolic metric (see Lemma~\ref{l:inj} and Proposition~\ref{small}), where we consider the decomposition of the kernel of $dQ^{\lambda}$ into eigenforms and eigenfunctions of $\Delta_{h}$ as the finite dimensional $L^{2}$ subspaces, which by perturbation from $\bbH^{7}$ there is not any such finite dimensional subspaces. Therefore the same analysis in $\{dQ^{\lambda}\}_{\lambda<\lambda_{0}}^{\lambda\notin\Lambda}$ applies. 

As for $\lambda\in \Lambda$, we will need to that the limit absorption principle as in Proposition~\ref{LAP} still applies for a Poincare--Einstein metric $h$ sufficiently close to the hyperbolic metric. Using the same strategy there, we separate into two cases. For the first case $\lambda=0$, it is the exact Hodge Laplacian $\Delta_{h}+36$ acting on 4-forms on an asymptotically hyperbolic manifold. Since both the result in~\cite{Vasy} and the resolvent construction~\cite[Theorem 9.5]{kantor2009eleven} cover such metrics, the same estimate holds.

For the other two cases $\lambda=16, 40$, we can use the same coutour deformation to see that $(dQ^{\lambda}\pm i\epsilon)^{-1}$ uniformly approaches $R^{\pm}$ if the estimate~\eqref{e:R0} holds. That is, we only need to see that 
$$
\lim_{\epsilon\downarrow 0}\|(\Lap_{h}-s(6-s)\pm i\epsilon)^{-1}-R_{h}(s^{\pm})\|_{x^{\delta}L^{2}_{0}(\BB^{7})\rightarrow x^{-\delta}H^{2}_{0}(\BB^{7})}=0.
$$
This can be seen by comparing the resolvent of a Poincare--Einstein metric $h$ sufficiently close to the hyperbolic metric. In particular, following the parametrix construction in~\cite{MR916753}, the resolvent is polyhomogeneous if the Poincare--Einstein metric is polyhomogeneous, and the same decay of the kernel of the resolvent on the double space is obtained. So the estimate of the resolvent kernel applies to asymptotic hyperbolic metric $h$, therefore Proposition~\ref{LAP} still holds. 

Combining the argument above, we obtain $R_{\pm}$.
\end{proof}

As a consequence, we define the following right inverse
\begin{equation}
(dQ)^{-1}:=\frac{1}{2}(R_++R_-),
\end{equation}
with the property that $dQ\circ (dQ)^{-1}=Id_{x^{\delta}H_{0,b}^{s,k}(M;K)}$. Here $R_{+}$ and $R_{-}$ are both complex-valued operators, however $(dQ)^{-1}$ being a real-valued operator which means that it has a real-valued integral kernel.

To get the main theorem, we will parametrize the domain by the boundary data, which amounts to show that there is a Poisson operator that maps boundary data into the kernel space. From the analysis in previous sections, we know that the only nontrivial kernel comes from the three pairs of special indicial roots, and therefore it is a geometric scattering problem. For the hyperbolic space, the scattering operator for Laplacian operators on functions and forms have been studied in various settings~\cite{MR916753, MR961517, mazzeo1990, melrose1995geometric,  MR1965361, MR2153454,MR2276069, lee2006fredholm}. In the context of $dQ^{\lambda}$, as we see from the computations in appendix B the three special cases all turn into problems of $(\Lap_{\bbH}-s_{\lambda})$ on functions and forms within the continuous spectrum, therefore we will use scattering operators to parametrize the kernels.

%Our idea is to use scattering matrix computation to show that, in the base case with hyperbolic metric,  the real solution space is parametrized by a real 3-form and two real functions on the 6-sphere, and therefore transversal to the range of $(R_++R_-)$. Then by a perturbation argument, we also define a domain for nearby Poincar\'{e}-Einstein metric by specifying an affine space parametrized by a real boundary data, and a translated operator $Q_v$. Then we would show that $Q_v\circ (dQ_{0,h})^{-1}$ is a local isomorphism on the range.

We start with the base case with the metric $g_{\bbH^{7}}\times \frac{1}{4}g_{\bbS^{4}}$.
\begin{lemma}
The real null space $\operatorname{Null}(dQ)\subset D_{k}(\delta)$ is parametrized by $\oplus V_{i}$ defined in Definition~\ref{bundle}. There is a Poisson operator $P: H^{k}(\SS^{6};V)\rightarrow x^{-\delta}H^{2,k}_{e,b}(M;K)$ such that $dQ\circ P=0$ and for any $v\in V$, $P(v)$ has the expansion as in~\eqref{expansion}.
\end{lemma}
\begin{proof}
For the construction of the Poisson operator, we follow Graham--Zworski~\cite{MR1965361} and Guillarmou~\cite{MR2153454} and first construct a formal solution operator $\tilde P$ by the standard asymptotic method. Given the boundary terms consisting of $v=(v_{1}^{+},v_{2}^{+}, v_{3}^{+})$ in Definition~\ref{bundle}, let $(k_{0},H_{0})$ be the leading order part of~\eqref{expansion}, which near the boundary is given by the scattering data and then extended from the boundary by a cut-off function into the interior, then we have
$$
dQ(k_{0},H_{0})=I(dQ)(k_{0},H_{0})+O(x^{3+\delta})\in x^{\delta}H_{e,b}^{0,k}(M;K)
$$
then one can solve the subsequent terms by iteratively constructing the expansion. By Borel's lemma, we arrive at a formal solution $\tilde P(v)=(\tilde k, \tilde H)$ with $dQ(\tilde k, \tilde H)=O(x^{\infty})$ where $(\tilde k, \tilde H)$ has the same leading order expansion. Then by the following proposition 3.4 in~\cite{MR1965361} there is a unique Poisson operator defined as
\begin{equation}\label{Poisson}
P=(I-R_{+}\circ dQ)\circ \tilde P
\end{equation}
with the correct mapping property. 

To show that the resulting kernel $P(v)$ is real, we use the description of scattering matrix in hyperbolic space in
Guillarmou--Naud~\cite{guillarmou2006wave}:
$$
S(s)=2^{n-2s}\frac{\Gamma(\frac{n}{2}-s)}{\Gamma(s-\frac{n}{2})}\frac{\Gamma\left(\sqrt{\Delta_{\bbS^n}+(\frac{n-1}{2})^2}+\frac{1-n}{2}+s\right)}
{\Gamma\left(\sqrt{\Delta_{\bbS^n}+(\frac{n-1}{2})^2}+\frac{n+1}{2}-s\right)},
$$
where if we put in $s=\theta_2^+$ denoted as $3+i\alpha$, we get the scattering operator
$$
S(3+i\alpha)=2^{-2\alpha i} \frac{\Gamma(-i\alpha)}{\Gamma(i\alpha)}\frac{\Gamma(\sqrt{\Lap_{\bbS^6}+\frac{25}{4}}+\frac{1}{2}+i\alpha)}{\Gamma(\sqrt{\Lap_{\bbS^6}+\frac{25}{4}}+\frac{1}{2}-i\alpha)}.
$$
 Since the scattering matrix is a function of the Laplacian on the boundary $\SS^6$, we can take the eigenvalue expansion on 6-sphere with real eigenfunction $f_\lambda$, and then consider the following expression, which is real and forms the leading order of the actual solution:
\begin{equation}\label{e:sca}
u_{\lambda}=x^{3+i\alpha}f_\lambda+x^{3-i\alpha}S(3+i\alpha)f_\lambda=x^{3+i\alpha}f_\lambda+x^{3-i\alpha}(2^{-2\alpha i}e^{2i\theta})f_\lambda.
\end{equation}
%If $f_\lambda$ is a solution, then $\bar f_\lambda$ is also a solution. However, $\bar f$ must be a complex multiple of $f$, as otherwise, there would be a two dimensional complex space of solutions instead of one. Therefore
%$$
%\bar f=cf, c\in \mathbb{C},
%$$  
%which shows
%$$
%f=e^{i\theta}|f|,
%$$
%i.e. $f$ is determined by a real function on the sphere. 
Here $\theta$ is a real number determined by
\begin{equation}\label{Gamma}
e^{2i\theta(\lambda)}=\frac{\Gamma(-i\alpha)}{\Gamma(i\alpha)}\frac{\Gamma\Big(\sqrt{\lambda+\frac{25}{4}}+\frac{1}{2}+i\alpha \Big)}{\Gamma\Big(\sqrt{\lambda+\frac{25}{4}}+\frac{1}{2}-i\alpha\Big)},
\end{equation}
by using the relation of 
$$
\Gamma(\bar z)=\overline{\Gamma(z)}
$$
so that the right hand side of~\eqref{Gamma} is a complex number with norm 1 and $\theta$ is a real number determined by $\lambda$.

Rearranging the expression, the solution in the eigenvalue $\lambda$ component in~\eqref{e:sca} is
\begin{equation}
\begin{aligned}
u_{\lambda} &=x^{3+i\alpha}f_\lambda+x^{3-i\alpha}2^{-2\alpha i}e^{2i\theta }f_\lambda\\
&=x^3 2^{-i\alpha}e^{i\theta} \left( (2x)^{i\alpha}e^{-i\theta} + (2x)^{-i\alpha} e^{i\theta} \right) f_\lambda\\
&=x^3 2^{1-i\alpha}e^{i\theta} \Re\left((2x)^{i\alpha}e^{-i \theta(\lambda)} \right) f_\lambda
\end{aligned}
\end{equation}
which is a product of a real function on $\SS^{6}$ with complex constant $2^{1-i\alpha}e^{i\theta}$. Therefore in this case, 
$$
u_{\lambda}=U_{\lambda}e^{it}, t=t(\lambda)\in \mathbb{R},
$$
which shows that $u_{\lambda}$ is given by a real function $U_{\lambda}$ for each $\lambda$. Hence the part of leading order expansion in~\eqref{expansion} given by $v_{2}^{+}x^{\theta_{2}^{+}}+S_{2}(v_{2}^{+})x^{\theta_{2}^{-}}$ can be written as $u=\sum_{\lambda}u_{\lambda}=\sum_{\lambda}U_{\lambda}e^{it(\lambda)}$ where each $U_{\lambda}$ is a real function. The same argument applies to the part given by $S_{3}$. 

For the part corresponding to $S_{1}$ which is the scattering matrix for forms on hyperbolic space, we use the explicit construction in proof of~\cite[Theorem B]{kantor2009eleven}. In particular, the construction using spherical harmonics in Chapter 7 gives the scattering matrix construction for the boundary data $v_{1}^{\pm}$, which is exactly what we need.
\end{proof}
	
And in the exact case with hyperbolic space, we can characterize the range of the Poisson operator using the two resolvents.
\begin{lemma}
The range of the Poisson operator $P$ acting on $H^{k}(\SS^{6};V)$ is the same as the range of $i(R_{+}-R_{-})$ acting on $x^{\delta}H_{e,b}^{0,k}(M;K)$.
\end{lemma}
\begin{proof}
By Stone's theorem, see for example (4.4) in~\cite{MR2276069}, the difference of $R_{+}$ and $R_{-}$ is given by $\frac{1}{2\pi i}dE$ where $dE$ is the spectral measure of $dQ$. To relate to the Poisson operator $P$, we consider its adjoint $P^{*}$ which is a map from $x^{\delta}H_{e,b}^{0,k}(M;W)$ to $H^{k}(\SS^{6}; V)$.  Since $P: H^{k}(\SS^{6}; V) \rightarrow x^{-\delta}H_{e,b}^{0,k}(M;W)$ is injective, by duality the above $P^{*}$ is surjective. We then use 
 use the formula
\begin{equation}
PP^{*}=\frac{1}{2\pi}dE=-i(R_{+}-R_{-})
\end{equation}
where the first equality can be proven using the same method from~\cite[Proposition 4.3]{Borthwick} by generalization from hyperbolic surfaces to the hyperbolic space $\HH^{7}$. Then by the surjectvity of $P^{*}$, we get that the range of $i(R_{+}-R_{-})$ acting on $x^{\delta}H_{e,b}^{0,k}(M;K)$ is the same as the range of $P$ acting on $H^{k}(\SS^{6};V)$.
%For any element $u\in x^{\delta}H_{e,b}^{0,k}(M;K)$, $dQ\circ (R_+-R_-)u=0$. 
%Since $R_+-R_-$ is imaginary, the range of $i(R_{+}-R_{-})$ is the real null space in $x^{-\delta}H_{e,b}^{2,k}(M;K)$, therefore is the same as the range of $P$. 
\end{proof}

Now we consider the perturbation from the base hyperbolic metric to a nearby Poincar\'{e}-Einstein metric $h$. As discussed before, we still have the two generalized inverses. 
\begin{lemma}\label{transv}
For a Poincar\'{e}--Einstein metric $h$ that is close to the hyperbolic metric $g_{\HH^{7}}$ with conformal infinity $\|\hat h-g_{\SS^{6}}\|_{H^{k}(\SS^{6})}<<1$, the range space of the sum of two generalized inverses $R_\pm$ is transversal to the range of their difference: $\operatorname{Range} (R_+ + R_-)$ is transversal to $\operatorname{Range} (R_+-R_-)$. 
\end{lemma}
\begin{proof}
We first show that the intersection of  $\operatorname{Range} (R_+ + R_-)$ and $\operatorname{Range} (R_+-R_-)$ in $D_{k}(\delta)$  contains only 0. If $f=(R_{+}-R_{-})u\neq 0$ is in the intersection, then $dQ( f)=0$. On the other hand, $f=(R_{+}+R_{-})v$ for some $v\neq 0$, hence $dQ( f)=dQ (R_{+}+R_{-})v =2v$ which combined with the previous line shows $v=0$ hence $f=0$.  Therefore the only element in the intersection of the two ranges is 0. 

We next show that any element $f\in D_{k}(\delta)$ can be written as the sum of two element such that 
\begin{equation}\label{e:trans}
f=u+v,  \ u\in \operatorname{Range} (R_+ + R_-), \ v\in \operatorname{Range} (R_+ - R_-).
\end{equation}
Since $f$ is an element in $D_{k}(\delta)$, it can be written as a $f=f_{0}+f'$ where $f'= O(x^{3+\delta})$ and $f_{0}$ are the leading asymptotics given by the special indicial roots $\theta_{i}^{\pm}, i=1,2,3,$ such that $f_{0}$ is formally annihilated by $dQ$. More specifically, $f_{0}$ is similar to the leading order part in~\eqref{expansion} where the scattering matrix part $S_{i}(v_{i}^{+})$ are replaced by arbitrary $v_{i}^{-} \in V_{i}$. Now we compute $R_{+}dQ(f)$ and $R_{-}dQ(f)$, and will show that $f-R_{\pm}dQ(f)$ is in the range of the Poisson operator $P$ and has incoming data $\ell_{\pm}=\{v_{i}^{\pm}\}_{i=1}^{3}$. First of all, $dQ(f-R_{+}dQ(f))=dQ(f)-dQ(f)=0$ so $u=f-R_{+}dQ(f)$ is in the kernel of $dQ$. And we can read the leading boundary asymptotic of $u$ which is given in the form of~\eqref{expansion} so that $\ell_{+}=(v_{1}^{+},v_{2}^{+},v_{3}^{+})$ are the incoming data. Now we need to show that $u=P\ell_{+}$. Since $u$ and $P\ell_{+}$ have the same incoming data and both solve $dQ(u)=dQ(P\ell_{+})=0$, we use the following pairing formula on $\ M_{\epsilon}=\HH^{7}_{\epsilon}\times \SS^{4}$ where $\HH^{7}_{\epsilon}$ is the manifold given by $\{x>\epsilon\}$ for some $\epsilon>0$:
\begin{equation}
\int_{M_{\epsilon}}u\overline{ dQ(P\ell_{+})}-\overline{dQ(u)} P\ell_{+}=\sum_{i=1}^{3}\frac{1}{2i \Im{\theta_{i}} }\int_{\partial M_{\epsilon}} v_{i}^{+}\overline {\tilde v_{i}^{-}}- v_{i}^{+}\overline {\tilde v_{i}^{-}}
\end{equation}
where $\{v_{i}^{\pm}\}$ are the incoming and outgoing data for $u$ and $\{\tilde v_{i}^{\pm}\}$ are the the incoming and outgoing data for $P\ell_{+}$. The pairing formula is derived from the usual Green's formula noting that $dQ$ is self-adjoint. And since $v_{i}^{+}=\tilde v_{i}^{+}$, by taking $\epsilon\rightarrow 0$, we get $v_{i}^{-}=\tilde v_{i}^{-}$. Therefore $u$ and $P\ell_{+}$ have the same incoming and outgoing data, and by unique continuation, we have $f-R_{+}dQ(f)=u=P\ell_{+}$. The same argument applies to $f-R_{-}dQ(f)$.

 Therefore
we write 
$$
f=R_{+}\circ dQ(f) + P \ell_{+} = R_{-}\circ dQ(f) + P \ell_{-},
$$
where $\ell_{+},\ell_{-}\in H^{k}(\SS^{6};V)$.  
From above we get 
$$
f=\frac{1}{2}(R_{+}+R_{-})\circ dQ(f)  + \frac{1}{2}P(\ell_{+}+\ell_{-}).
$$
Then we let $u=\frac{1}{2}(R_{+}+R_{-})\circ dQ(f)$ and $v=\frac{1}{2}P(\ell_{+}+\ell_{-})$. Since $u \in \operatorname{Range} (R_+ + R_-)$ by definition, and $v\in \operatorname{Range}(P)=\operatorname{Range}(R_+ + R_-)$ by the previous lemma, we get~\eqref{e:trans}.
%For the hyperbolic case, the range of $(R_{+}-R_{-})$ is the kernel of $dQ$. However, the range of $R_++R_-$ doesn't contain any element of the kernel. Otherwise And for any $f \in D_{k}(\delta)$, consider $\tilde f$ the projection of $f$ off the kernel f of $dQ$, then $\frac{1}{2}(R_{+}+R_{-})dQ(\tilde f)=\tilde f$. Therefore the two range spaces are transversal in $D_{k}(\delta)$.

Since transversality is stable under small perturbations, the result follows for nearby Poincar\'e-Einstein metrics.
\end{proof}

The Poisson operator defined in~\eqref{Poisson} (with respect to the hyperbolic metric) exists for nearby Poincar\'e-Einstein metric as well. %However, $P$ maps to a real element in the domain which is not necessarily an element in the null space, but a perturbation of the kernel. 
The range of $P$ is still the same as the range of $i(R_{+}-R_{-})$. With this we conclude:

\begin{proposition}\label{PtR}
The range of the Poisson operator $P$ acting on $H^{k}(\SS^{6};V)$ is transversal to the range of $(dQ)^{-1}=(R_++R_-)$ acting on $x^{\delta}H_{e,b}^{2,k}(M;K)$.
\end{proposition}
\begin{proof}
Since the range of $i(R_{+}-R_{-})$ is transversal to the range of $(R_{+}+R_{-})$ by lemma~\ref{transv}, the transversality in the statement follows from that the range of $P$ is the same as the range of $i(R_{+}-R_{-})$.
\end{proof}

\section{Solvability of the nonlinear operator}\label{fit}

From the discussion of the linear operator $dQ$ above, we now can apply the implicit function theorem to get results for the nonlinear operator. To do this we first need to show that the nonlinear terms are controlled. Then we will use a perturbation argument to show that for each solution with Poincar\'{e}--Einstein metric close to hyperbolic metric, the nearby solutions are parametrized by the three parameters on $\bbS^6$ as in the linear case.

To deal with the fact that the domain changes with the base metric and the boundary parameters, we will use an implicit function theorem for a map from range space to itself, and show this map is a perturbation of identity, therefore an isomorphism.

First of all we define the domain that depends on the choice of the base Poincar\'{e}--Einstein metric $h$ and the boundary parameter $v=(v_1,v_2,v_3) \in V$. From proposition~\ref{PtR}, we know that the image of $(dQ)^{-1}=\frac{1}{2}(R_++R_-)$ is transversal to the image of the Poisson operator $P$, which for a nearby Poincar\'{e}--Einstein metric is still the kernel of the linearized operator. For each fixed parameter $v$, we define the domain as an affine section of $(dQ)^{-1}(x^{\delta}H_{e,b}^{0,k}(M;K))$ translated by $Pv$.

\begin{definition} (Domain of nonlinear operator)\label{nldomain}
For a Poincar\'{e}--Einstein metric $h$ with $\|\hat h-g_{\SS^{6}}\|_{H^{k}(\SS^{6})}<<1$ and a set of parameters $v=(v_1,v_2,v_3) \in V$, the domain $D_{h,v}$ for the nonlinear operator $Q$ is defined as 
$$
D_{h,v} :=\Big\{ (dQ)^{-1}f+Pv\,\Big|\, f \in x^{\delta}H_{e,b}^{0,k}(M;K)\Big\}
$$
where $dQ=dQ_{t,W}$ is the linearization at $t=h\times \frac{1}{4}g_{\SS^{4}}$, and $P$ is the Poisson operator for $t$.
\end{definition}
Note that the domain depends on the choice of $h$ and $v$, where the dependence of $h$ comes from the construction $(dQ)^{-1}=\frac{1}{2}(R_++R_-)$.  Because of the transversality from~\ref{PtR}, $D_{h,v}$ can be viewed as a slice in $D_{k}(\delta)$:
$$
\cup_{v\in V}D_{h,v} = D_{k}(\delta).
$$
The domain has the property that, if $h=g_{\bbH^{7}}$ is the hyperbolic metric, then each slice $D_{h,v}$ is mapped by $dQ$ isomorphically back to the range space $x^\delta H_{e,b}^{0,k}(M;K)$ where the kernel in each slice is exactly $Pv$.

For nearby metric $h$, one important property of this domain is that $D_{h,v}$ is mapped surjectively to the range space $x^{\delta}H_{e,b}^{0,k}(M;K)$ by the linearized operator $dQ=dQ_{t,W}$. 

\begin{lemma}
Acting on the domain defined in~\ref{nldomain}, the linearized operator
$$dQ: D_{h,v} \rightarrow x^{\delta}H_{e,b}^{0,k}(M;K)$$
is a surjective map.
\end{lemma}
\begin{proof}
By direct computation, for any $f\in x^{\delta}H_{e,b}^{0,k}(M;K)$ and $v\in V$,
$$
dQ\left(\frac{1}{2}(R_++R_-)f + Pv\right)=dQ(dQ)^{-1}f+ dQ(Pv)=f+ dQ(Pv)=f.
$$
Here we used the fact that $R_+$ and $R_-$ are both right generalized inverses for $dQ$. And by definition of the Poisson operator $dQ(Pv) =0$ for any $v$. Since $f$ can be any element in the vector space, it follows the range of $dQ$ acting on $D_{h,v}$ is the whole space.
\end{proof}

Next we show that the nonlinear terms are well controlled, the nonlinear operator $Q$ maps $D_{h,v}$ to $x^{\delta}H^{2,k}_{e,b}(M;K)$. This is proved by showing that the difference of $Q$ and $dQ$ is small. And we only need to consider the action on the $x^{\delta}H_{e,b}^{2,k}(M;K)$ since the only part that has a worse decay is eliminated by $dQ$.
\begin{lemma}\label{l:nl1}
For sufficiently large $k$, the product type nonlinear terms: $F\circ F - d(F\circ F)$, and $F \bigwedge F - d(F\wedge F)$ are both contained in $x^{\delta} H_{e,b}^{2,k}(M;K)$.
\end{lemma}
\begin{proof}
The nonlinear parts are $F\wedge F$ and $F\circ F$ which are products of two elements in the range space $x^{\delta} H^{2,k}_{e,b}(M;K)$. Take a basis of the edge bundles, these can be considered locally as functions in $x^{\delta}H_{e,b}^{2,k}(M)$. Using proposition~\ref{product}, we know that for $r>-3$, and sufficiently large $k$, and any $f,g \in x^{r}H_{e,b}^{s,k}(M)$, the product $fg$ is also in $x^{r}H_{e,b}^{s,k}(M)$. Since in our case $\delta>0$, the result follows.
\end{proof}

The other nonlinear term is the remainder from the linearization of $\operatorname{Ric}$ operator, for which we show below that it is also contained in the range space. 
\begin{lemma}\label{l:nl2}
The nonlinear remainder of $\operatorname{Ric}$, $\operatorname{Ric}-d(\operatorname{Ric})$ acting on~$k\in x^{\delta}H_{e,b}^{0,k}(M;K)$ is contained in $x^{\delta}H_{e,b}^{0,k}(M;K)$.
\end{lemma}
\begin{proof}
We compute the linearization $d(\operatorname{Ric})$, which acting on a 2-tensor $h$ can be written as
$$
d(\Ric)[h]=\frac{-1}{2}g^{ml}(\nabla_m \nabla_l h_{jk}-\nabla_m\nabla_k h_{jl}-\nabla_l\nabla_j h_{mk}-\nabla_j\nabla_k h_{ml}).
$$
Comparing $\Ric$ and $d(\Ric)$, the difference is a 3rd order polynomial of $g,g^{-1}$ and first order derivatives of these with smooth coefficients. Since the metric component $g$ and $g^{-1}$ are smooth, hence in $x^{\delta}H_{e,b}^{s,k}(M;K)$, it follows again by the algebra property that their product is contained in $x^{\delta}H_{e,b}^{s,k}(M;K)$.
\end{proof}

We next define a nonlinear operator $Q^{h,v}$ that is defined as a translation of the original operator $Q=Q^{h,0}$ defined in~\eqref{gaugedequ}.
\begin{definition}
For a Poincare--Einstein metric $h$ such that $\|\hat h -g_{\SS^{6}}\|_{H^{k}(\SS^{6})}<<1$,  we define the parametrized nonlinear operator $Q^{h,v}$ as:
$$
Q^{h,v}: D_{k}(\delta) \rightarrow x^{\delta}H_{e,b}^{0,k}(M;K), \ 
u\mapsto Q^{h,0}(u+Pv),
$$
where $Q^{h,0}$ is the nonlinear operator defined in~\eqref{gaugedequ} and $P$ is the Poisson operator with respect to the Poincare--Einstein metric $h$.
\end{definition}

As a translation of the original operator, the linearization of $Q^{h,v}$  is closely related to the original linearized operator $dQ$:
\begin{lemma}\label{linear}
 The linearization of $Q^{h,v}$ at $(t,W)-Pv=(h\times \frac{1}{4}g_{\SS^{4}}, 6 \Vol_{\SS^{4}})-Pv$, denoted as $dQ^{h,v}_{(t,W)-Pv}$, is the same as $dQ^{h,0}_{t,W}$:
$$
dQ^{h,v}_{(t,W)-Pv}(u)=dQ^{h,0}_{(t,W)}(u), \ \forall u=(k,H)\in \eS(M) \oplus \eW{4}.
$$
\end{lemma}
\begin{proof}
Since $Q^{h,v}$ is defined as a translation of $Q^{h,0}$ by $Pv$:
$$
Q^{h,v}(\cdot)=Q^{h,0}(\cdot + Pv),
$$
and the nonlinear terms in $Q$ are all quadratic, therefore by definition of the linearization we have
$$
\begin{aligned}
dQ^{h,v}_{(t,W)-Pv}(u)=\lim_{s\rightarrow 0} \frac{Q^{h,v}((t,W)-Pv+su)}{s}=\lim_{s\rightarrow 0} \frac{Q^{h,0}((t,W)+su)}{s}
\\=dQ^{h,0}_{t,W}(u).
\end{aligned}
$$
\end{proof}

For simplicity, $dQ$ (and therefore $(dQ)^{-1}$) will be the abbreviation for $dQ_{t,W}^{h,0}$, and the linearization of $dQ^{h,v}$ will be noted separately.
The composed operator $Q^{h,v}\circ (dQ)^{-1}$ is this well-defined operator as a map on the following space:
$$
Q^{h,v}\circ (dQ)^{-1}: x^{\delta}H_{e,b}^{0,k}(M;K)\rightarrow  x^{\delta}H_{e,b}^{0,k}(M;K).
$$
$$
f \mapsto Q^{h,0}\Big(\frac{1}{2}(R_++R_-)f+Pv\Big).
$$
We now discuss the properties of this operator using the following implicit function theorem, which can be found for example in~\cite[Theorem 5.9]{MR1666820}.

\begin{lemma}[Implicit function theorem] For two Banach spaces $V$ and M, if $f$ is a smooth map $f: V\times M \rightarrow M$ near a point $(v_0,m_0) \in V\times M$ with $f(v_0,m_0)=c$, and the linearization of the map with respect to the second variable $df_2(v_0,m_0): M\rightarrow M$ is an isomorphism, then there is neighborhood $v_0\in U \subset V$ and a smooth map $g:V\rightarrow M$, such that $f(v,g(v))=c, \forall v\in U$.
\end{lemma}

\begin{theorem}\label{main}
For  $k\gg 0$, $\delta<\in (0,\bar \delta)$, there exists $\rho>0$ and $\epsilon>0$, such that, for a Poincar\'{e}--Einstein metric $h$ that is sufficiently close to the base metric $g_{\bbH^{7}}$ with $\|\hat h-g_{\SS^{6}}\|_{H^{k}(\SS^{6})}<\epsilon$, and any boundary value perturbation $v\in V$ with $\|v\|_{H^k(\bbS^{6};V)}<\rho$, there is a unique solution $u=(g, H) \in D_{v,h} \subset x^{-\delta}H_{e,b}^{s,k}(M;K)$ satisfying the gauged supergravity equations $Q(u)=0$ with the leading expansion of $(g-h\times \frac{1}{4}g_{\SS^{4}},H-6 \Vol_{\bbS^{4}})$ given by~\eqref{expansion}.
\end{theorem}

To prove the theorem, we will apply the implicit function theorem to the following operator:

$$
Q^{h,\cdot}\circ (dQ)^{-1}: H^{k}(\bbS^{6};V) \times  x^{\delta}H_{e,b}^{0,k}(M;K) \rightarrow x^{\delta}H_{e,b}^{0,k}(M;K)
$$
$$
(v, f)\mapsto Q^{h,v} \circ (dQ)^{-1}(f)
$$
From the previous discussion, this map is well defined. The following is a consequence of Lemma~\ref{linear}.

\begin{lemma}
The linearization of $Q^{h,v} \circ (dQ)^{-1}$ at point $(v,f)=(0,0) \in V \times x^{\delta}H_{e,b}^{0,k}(M;K)$ is an isomorphism.
\end{lemma}
\begin{proof}
From Lemma~\ref{linear} we know that at the point $(v,f)=(0,0) \in V\times x^{\delta}H_{e,b}^{0,k}(M;K)$ the linearization, which is the composition of linearized operators, is
$$
d(Q^{h,v} \circ (dQ)^{-1})_{(0,0)}=dQ^{h,0}_{t,W}\circ (dQ)^{-1}=\operatorname{Id}: x^{\delta}H_{e,b}^{0,k}(M;K) \rightarrow x^{\delta}H_{e,b}^{0,k}(M;K).
$$
\end{proof}

\begin{lemma}
For a given metric $h$, the map $Q^{h,v} \circ (dQ)^{-1}$ as an edge operator varies smoothly with the parameter $v\in V$.
\end{lemma}
\begin{proof}
From the construction of $(dQ)^{-1}$ we know it is an edge operator. And from the discussion for $Q^{h,0}$, this nonlinear operator is also edge. Now we we only need to show that when the nonlinear operator $Q$ applies to elements of type $f+Pv$, it varies smoothly with the parameter $v$. This follows from the algebra property and the fact that a second order elliptic edge operator maps from $H_e^s(M)$ to $H_e^{s-2}(M)$ smoothly as shown in proposition~\ref{A2}. 
\end{proof}

Now as a direct result of the implicit function theorem, we now prove the main theorem: 

\begin{proof}[Proof of Theorem~\ref{main} ]
Using the implicit function theorem, we can find neighborhoods of $0\in V$ and $f \in x^{\delta}H_{e,b}^{0,k}(M;K)$, in this case, $U_1=\{v\in V:\|v\|_{H^k(\bbS^{6};V)}<\rho\}$ and $U_2=\{f \in x^{\delta}H_{e,b}^{0,k}(M;K):\|f\|_{x^{\delta}H_{e,b}^{0,k}(M;K)}<\rho_2\}$, such that the nonlinear map $Q^{h,v}\circ (dQ)^{-1}$ is a  bijective smooth map on $U_2$ for any $v \in U_1$ . And this gives us the parametrized map $g$ from $U_1$ to $U_2$ such that
$$
Q^{h,v} (dQ)^{-1}(g(v))=0.
$$
And we can rewrite it as
$$
Q^{h,0}\left((dQ)^{-1}(g(v))+Pv\right)=0.
$$
That is, for each parameter set $v$, $u=(dQ)^{-1}(g(v))+Pv$ is the unique solution in the space $D_{v,h} \subset x^{-\delta}H_{e,b}^{s,k}(M;K)$. When the base metric is $g_{\HH^{7}}$, since the nonlinearity is quadratic, the leading order expansion of the solution is given by $Pv$. By continuity, $(dQ)^{-1}(g(v))$ is $O(x^{3+{\delta}})$ when $\|g-h\|$ and $\|v\|$ is sufficiently small. Therefore the leading order behavior is again given by $Pv$ which is~\eqref{expansion}.
\end{proof}

%main theorem
Next we show that the solution obtained above is smooth if the boundary data is smooth.
\begin{proposition}
If the boundary data $v \in C^{\infty}(\bbS^6;V)$, then the solution $u$ is in $C^{\infty}(M;K)$.
\end{proposition}
\begin{proof}
This is done by elliptic regularity. For any $k$,
$$
\|u-(h\times \frac{1}{4}g_{\SS^{4}},6 \Vol_{\SS^{4}})\|_{x^{-\delta}H^{2,k}_{e,b}(M;K)} \leq C(\|v\|_{H^k(S^6;V)}+\|Q(u)\|_{x^{\delta}H_{e,b}^{0,k}(M;K)}).
$$
For the linearized operator we use the elliptic estimate for edge operators, and 
the difference with nonlinear operator is lower order therefore can be controlled by $\|u\|_{x^{\delta}H_{e,b}^{2,k}(M;K)}$ so it is absorbed to the left hand side. Since the estimate holds for every $k$ we get the smoothness of $u$.
\end{proof}

We can also obtain the polyhomogeneous expansion of the solution. \begin{proposition}
When the boundary data $v \in C^{\infty}(\bbS^6;V)$, the solution $u$ has a classical polyhomogeneous expansion in the sense of~\cite{melrose1993atiyah}, with leading terms given by~\eqref{expansion} and the exponent of the logarithmic terms grows at most linearly with the order.
\end{proposition}
\begin{proof}
We solve the problem iteratively to obtain a formal expansion. Denote the base data $u_{0}=(h\times g_{\frac{1}{4}\SS^{4}}, 6\Vol_{\SS^{4}})$. For the first order problem, from the linearization and its inverse construction, we have $u_1=(g_{1},F_{1})$ as the leading terms of~\eqref{expansion} such that
$$
Q\left(u_{0}+u_1\right)=x^{3+\delta}e_1, \quad e_1\in \CI(M;K),
$$
which holds for any $\delta\in (0,1)$. 
The next step is to solve the following equation
$$
dQ(u_{2})=-x^{3+\delta}e_1
$$
which gives $u_{2}\in \mathcal{O}(x^{3+\delta})$. This way we solve away the $x^{3+\delta}e_1$ term and one possible $\log(x)$ term appears because of the possible appearance of the next indicial roots. By looking at the error term from Lemma~\ref{l:nl1} and~\ref{l:nl2}, we have that the error term 
$
(Q-dQ)(u_{0}+u_{1}+u_{2})
$
consists of a 3rd order polynomial in the metric and its inverse and is quadratic in the 4-form part. Therefore by using the estimate~\eqref{e:product}, we have that 
$$(Q-dQ)(u_{0}+u_{1}+u_{2})=Q(u_{0}+u_{1}+u_{2})-Q(u_{0}+u_{1})=u_{2}P_{2}(u_{1},u_{2})=x^{4+\delta}e_{2}.$$ 
Then iteratively we solve the equation
$$
dQ(u_{k})=Q(\sum_{i=1}^{k-1} u_{i})-Q(\sum_{i=1}^{k-2}u_{i})=u_{k-1}P_{k-1}(u_{1}, \dots, u_{k-1}),
$$
where $P_{k-1}$ is a polynomial. Iteratively we can show that the right hand side has a decay rate of $\mathcal{O}(x^{k+1+\delta})$ because the lowest order term $u_{k-1}u_{1}$ in the polynomial is of order at least $\mathcal{O}(x^{k+1+\delta})$. And therefore apply the inverse of $dQ$ we have $u_{k}=\mathcal{O}(x^{k+1+\delta})$. Moreover we have an expansion for $u_{k}$ such that
$$
u_k=\sum_{j}x^{s_{j,k}}(\sum_{i=0}^{\infty}x^{i}u_{i,j,k}), \ u_{i,j,k} \in (\log x)^{i}\CI(M;K),
$$
with $s_{j,k}$ being the indicial roots bigger than $k+1+\delta$ and each time the power of $\log$ increases by at most one with the power of $x$. Combining all the terms
$$
u\sim \sum_{i=0}^{\infty} u_{i}
$$
we get the polyhomogeneous expansion.
\end{proof}

Finally we prove the main theorem for the original supergravity operator. 

\begin{proof}[Proof of Theorem]
From proposition~\ref{gaugeeq}, there is a diffeomorphism $g\rightarrow \tilde g$ such that we obtain the solution to the original supergravity equations $S(\tilde g,H)=0$. The parametrization of solutions $(g,H)$ to the gauged equation is given in Theorem~\ref{main}. The regularity of $(\tilde g, H)$ is the same as $(g,H)$ because of the diffeomorphism. Since $\tilde g$ and $g$ differ by a lower order term $O(x^{3+\delta})$, they have the same leading order expansion as~\eqref{expansion}.
\end{proof}

\appendix
\section{Edge operators}\label{edge}
\begin{proposition}\label{A1}
$H_{e,b}^{s,k}(M)$ is defined independent of the order of applying edge- and b-vectors.
\end{proposition}
\begin{proof}
We prove it by induction. Take $s=k=1$, using the commutator relation $[\mathcal{V}_{e}, \mathcal{V}_{b}] \subset \mathcal{V}_{b}$, we have 
$$
V_{e}V_{b}u=V_{b}V_{e}u+V'_{b}u, \ V_{e} \in \mathcal{V}_{e}, \ V_{b}, V_{b}' \in \mathcal{V}_{b}.
$$
Therefore $\mathcal{V}_{e}\mathcal{V}_{b}u\in L^{2}(M)$ if and only if $\mathcal{V}_{b}\mathcal{V}_{e}u\in L^{2}(M)$ (since $u\in H^{1}_{b}(M)$ is implied by both sides.) For $s,k>1$, for an arbitrary order of vector fields applied to $u$, we use the commutator to reduce to the sum $\sum_{i=0}^{s}\mathcal{V}_{b}^{k}\mathcal{V}_{e}^{i}u$ and use the induction that $H_{e,b}^{s,k}\subset H_{e,b}^{i,k}$ for any $i<s$.
\end{proof}

\begin{proposition}\label{A2}
Any $m$-th order edge operator P maps $H_{e,b}^{s,k}(M) $ to $H_{e,b}^{s-m,k}(M)$, for $m\leq s$.
\end{proposition}
\begin{proof}
Locally, any m-th order edge operator P can be written as
$$
P=\sum_{j+|\alpha|+|\beta|\leq m} a_{j,\alpha.\beta}(x,y,z)(x\partial_x)^j(x\partial_y)^{\alpha}\partial_z^{\beta}
$$

If we can prove for $m=1$, P maps $H_{e,b}^{s,k}(M) $ to $H_{e,b}^{s-1,k}(M)$, then by induction we can prove for any $m$. Therefore we restrict to the case $m=1$.
 
We just need to check that, for a function $u \in H_{e,b}^{s,k}(M)$, $Pu$ satisfies 
$$V_e^{i}Pu \in H_b^k(M), 0\leq i\leq s-1.$$

The we prove the proposition by induction on k. For k=1 case, since a boundary vector field $V\in \cV_b(M)$ satisfies the commutator relation $VP=PV+[V,P]$ where the Lie bracket $[V,P] \in \cV_b$, then 
$$
VP(u)=PV(u)+V_b(u)
$$ 
by definition of $u\in H_{e,b}^{s,k}$, both $V(u)$ and $V_b(u)$ are in $H_e^s(M)$, therefore $PV(u) \in H_e^{s-1}(M)$.

If it holds for $k-1$, then by the relation
$$
V_b^{k}P(u)=V_b^{k-1}P V_b(u)+V_b^{k}(u),
$$
since $V_b(u) \in H^{s,k-1}_{e,b}$ and from induction assumption $P V_b(u) \in H_{e,b}^{s-1, k-1}$, therefore the first term $V_b^{k-1}PV_b(u) \in H_e^{s-1}(M)$, and the second term is in $H_e^s$ by definition. Therefore $Pu \in H_{e,b}^{s-1, k-1}$, which completes the induction.
\end{proof}

\begin{proposition}\label{product}
For sufficiently large $k$, $r\geq -3$, and any $s\in \bbR$, $x^rH_{e,b}^{s,k}(M)$ is an algebra. 
\end{proposition}
\begin{proof}
We first prove that, for the case $r=-3$, the b-Sobolev space $x^{-3}H_b^k(M)$ is an algebra for sufficiently large k. Working in the upper half plane model with coordinates $(x,y_1,\dots y_n,z_{1}, \dots z_{n'})$. For any element $f \in x^{-3}H_b^k(M)$, by definition, its Sobolev norm is
$$
\int_{\bbR^{+} \times \bbR^{n+n'}} |V_b^k (x^3 f)|^2x^{-7}dxdydz 
$$
Since the commutator relation satisfies $[V_b, x^3]f = x^{3}V_{b}f+3x^{3}f$, the definition of the Sobolev norm above is the same as
$$
\int |x^3 (V_b^k f)|^2x^{-7}dxdydz
$$
We do a coordinate transformation to change the problem back to $\mathbb{R}^{m}$ with $m=1+n+n'$: let $\rho=ln(x)$, then $x\partial_x=\partial_\rho$.  Therefore under the new coordinates, the b-vector fields are spanned by $(\partial_\rho, \partial_y, \partial_z)$. Let $F$ be the function after coordinate transformation
$$
F(\rho, y,z)=f(e^{\rho},y,z)
$$
then from the discussion above we can see the norm for $x^{-3}H_{b}^{k}(M)$ is characterized by
$$
\|f\|^2_{x^{-3}H_b^k}=\int_{\bbR^{+} \times \bbR^{n+n'}} |x^3 (V_b^k f)|^2x^{-7}dxdydz=\int_{\bbR^{m}} |V_b^k F|^2 d\rho dydz<\infty
$$
which means $F \in H^k(\mathbb{R}^{m})$. From~\cite{MR2744150}, the usual Sobolev space $H^{k}(\mathbb{R}^m)$ is closed under multiplication if and only if $k>\frac{m}{2}$. Therefore, take two elements $f,g \in x^{-3}H_b^k(M)$, then the corresponding functions in $\mathbb{R}^m$ satisfy $FG \in H^k(\mathbb{R}^m)$. It follows that $fg \in x^{-3}H_b^k(M)$ by taking the inverse coordinate transformation.

Then it is easy to see that $x^rH_b^k(M)$ is an algebra for $r>-3$ from the result above:
\begin{equation}\label{e:product}
\begin{gathered}
(x^{r}H_{b}^k)\cdot (x^{r}H_{b}^k)=x^{3+r}(x^{-3}H_b^k) \cdot x^{3+r}(x^{-3}H_b^k) \\
\subset\ x^{6+2r}(x^{-3}H_b^k)\ \subset x^{r}H_b^k(M).
\end{gathered}
\end{equation}

%Take the transform $\rho=ln(x)$, then $x\partial_x=\partial_\rho$.
Now that we proved $x^{r}H_b^k(M)$ is closed under multiplication, then we want to prove $x^{r}H_{e,b}^{s,k}(M)$ is also an algebra for any $s$. For any functions $f,g \in x^{r}H_{e,b}^{s,k}(M)$, by Leibniz rule,
$$
V_e^j(fg)=\sum_{i=0}^j V_e^i(f) V_e^{j-i}(g)
$$  
where by assumption, both $V_e^i(f)$ and $V_e^{j-i}(g)$ are in $x^{r}H_b^k(M)$, therefore their product is also in $x^{r}H_b^k(M)$ from the above result. Hence we proved $V_e^j(fg) \in x^{r}H_b^k(M)$ for $0\leq j \leq s$, which shows $fg \in x^{r}H_{e,b}^{s,k}(M)$.
\end{proof}

\section{Computation of the indicial roots }\label{indicial}
\subsection{Hodge decomposition}
The system contains the following equations, where the $(i,j)$ notations mean the splitting of  forms with respect to the product structure of $\bbB^7 \times \bbS^4$, i.e. $H_{(i,j)}$ has the form $\sum_{k}f_{k}\alpha_{k} \wedge \beta_{k}, \ \alpha_{k} \in {}^{e}\wedge^{i}\bbH^{7}, \ \beta_{k} \in \wedge^{j}\bbS^{4}$. Here $\Delta_{\SS^{4}}, \Delta_{\HH^{7}}$ denote the Hodge Laplacians on $\SS^{4}$ and $\HH^{7}$, and $\Delta^{rough}$ denotes the rough Laplacians. The derivation of the first two groups of equations below are straightforward from the splitting of form degrees. And for the third group we refer to the first group of equations in~\cite[Theorem 5.1]{kantor2009eleven}.
\begin{itemize}
\item From the first order equation
\begin{align}
(7,1):&\quad
\begin{gathered}
6d_{\HH^7}*_7k_{(1,1)}+3d_{\SS^4}(Tr_{H^7}(k)-Tr_{S^4}(k))\bigwedge
{}^7V\\+d_{\SS^4}*H_{(0,4)}+d_{\HH^7}*H_{(1,3)} =0
\end{gathered}
\\
(6,2):&\quad d_{\SS^4}*H_{(1,3)}+d_{\HH^7}*H_{(2,2)}+6d_{\SS^4}*_7k_{(1,1)} =0
\\
(5,3):&\quad
d_{\SS^4}*H_{(2,2)}+d_{\HH^7}*H_{(3,1)}=0
\\
(4,4):&\quad
d_{\SS^4}*H_{(3,1)}+d_{\HH^7}*H_{(4,0)}+W\wedge H_{(4,0)}=0
\end{align}

\item From $dH=0$
\begin{align}
d_{H}H_{(0,4)}+d_{S}H_{(1,3)}=0\\
d_{H}H_{(1,3)}+d_{S}H_{(2,2)}=0\\
d_{H}H_{(2,2)}+d_{S}H_{(3,1)}=0\\
d_{H}H_{(3,1)}+d_{S}H_{(4,0)}=0\\
d_{H}H_{(4,0)}=0
\end{align}

\item From the laplacian:
\begin{align}
\frac{1}{2}\Delta_{\SS^4}^{rough} k_{Ij}+\frac{1}{2}\Delta_{\HH^7}^{rough}
k_{Ij}+6k_{Ij}-3*_{\SS^4}H_{(1,3)}=0\\
\frac{1}{2}(\Delta_{\SS^4}+\Delta_{\HH^7}^{rough})k_{IJ}-k_{IJ}-6Tr_{\SS^4}(k)t_{IJ}+Tr_{\HH^7}(k)
t_{IJ}+2H_{(0,4)}t_{IJ}=0\\
\frac{1}{2}(\Delta_{\SS^4}^{rough}+\Delta_{\HH^7})k_{ij}+4k_{ij}+8Tr_{\SS^4}(k)t_{ij}-H_{(0,4)}t_{ij}=0
\end{align}
\end{itemize}

\subsection{Indicial roots}
Then we decompose further with respect to Hodge theory on the sphere, and compute the indicial roots for each part.

\begin{enumerate}
\item \textbf{$\hat k_{IJ}$: trace-free 2-tensor on $\bbH^7$}
\\
The equation is
$$
(\Lap_{\SS^4}+\Lap_{\HH^7}^{rough}-2)\hat k_{IJ}=0,
$$
and the indicial equation is
$$
(\lambda-s^2+6s)\hat k_{IJ}=0.
$$
we have indicial roots
$$
s=3\pm \sqrt{9+\lambda}
$$
The first pair of indicial roots, when $\lambda=0$, correspond to the perturbation of hyperbolic metric to Poincar\'{e}--Einstein metric.

\item \textbf{$\hat k_{ij}$: trace-free 2-tensor on $\bbS^4$}
\\
The equation is
$$
\Lap_{\SS^4}^{rough} \hat k_{ij}+\Lap_{\HH^7} \hat k_{ij}+8\hat k_{ij}=0
$$
where indicial equation is
$$
(\lambda-s^2+6s+8)\hat k_{ij}=0,
$$
indicial roots
$$
s=3\pm \sqrt{17+\lambda}.
$$

\item \textbf{$H_{(4,0)}$ with harmonic functions}
\\
We have \begin{eqnarray}
d_{\HH^7}*H_{(4,0)}+W\wedge H_{(4,0)}=0\\
d_{\HH^7} H_{(4,0)}=0
\end{eqnarray}

The second equation can be deduced from the first one. And by applying $d_{\HH^{7}}*_{\SS^{4}}$ to the first equation, we get
\begin{equation}\label{e:s1}
\Delta_{\HH^{7}}H_{(4,0)}+36 H_{(4,0)}=0.
\end{equation}

 Since the indicial operator for $d_{\HH^7}$ on a $k$-form is
$$
I[d](s)w=(-1)^k(s-k)w\wedge dx/x
$$
Let
$$
H_{(4,0)}=T+dx/x\wedge N
$$
be the decomposition with respect to tangential and normal decomposition, then the indicial equations are
$$
-(s-3)(*_{\SS^6} N)\wedge dx/x-6 dx/x\wedge N=0
$$
$$
(s-4)T\wedge dx/x=0
$$
where the first equation gives
$$
(s-3)*_{\SS^6}N-6N=0
$$
i.e. N is an eigenform of $*_{\SS^6}$ and the corresponding indicial roots are
$$
s_{3}^{-}=3-6i,  N\in \bigwedge^3(\bbS^6); *_{\SS^6}N=iN; 
$$
$$
s_{3}^{+}=3+6i: N \in \bigwedge^3(\bbS^6); *_{\SS^6} N=-iN.
$$
And plugging into the second equation, we have the vanishing of tangential form 
$$
T=0.
$$

Therefore the kernel in this case is
$$
H_{(4,0)}=dx/x \wedge N, N \in \{\bigwedge^3(\bbS^6), *_{\SS^6} N=\pm iN\}.
$$

\item \textbf{$\tau=\frac{1}{4}Tr_{\SS^4}(k), \sigma = \frac{1}{7}Tr_{\HH^7}(k), k_{(1,1)}, H_{(0,4)}, H_{(1,3)}$ on eigenfunctions / exact 1-form / coexact 3-form / exact 4-form}
\\
We have the following equations:
\begin{equation}\label{Exact1}
\begin{gathered}
6d_{\HH^7} *_{\HH^7} k_{(1,1)}^{cl}+d_{\SS^4} (3Tr_{H}(k)-3Tr_{\SS^4}(k))\bigwedge {}^7V
+d_{\SS^4}*H_{(0,4)}^{cl}+d_{\HH^7}*H_{(1,3)}^{cc}=0
\end{gathered}
\end{equation}
\begin{equation}\label{Exact2}
d_{\HH^7} H_{(0,4)}^{cl}+d_{\SS^4}H_{(1,3)}^{cc}=0
\end{equation}
\begin{equation}\label{Exact3}
d_{\HH^7} H_{(1,3)}^{cc}=0
\end{equation}
\begin{equation}\label{Exact4}
\Delta_{\SS^4} k_{(1,1)}^{cl}+\Delta_{\HH^7}
k_{(1,1)}^{cl}+12k_{(1,1)}^{cl}-6*_{\SS^4}H_{(1,3)}^{cc}=0
\end{equation}
\begin{equation}\label{Exact5}
\Lap_{\SS^4} \tau+\Lap_{\HH^7} \tau+72\tau-8*_{\SS^4}H_{0,4}^{cl} =0
\end{equation}
\begin{equation}\label{Exact6}
\Lap_{\SS^4} \sigma +\Lap_{\HH^7} \sigma +12\sigma +4*_{\SS^4}H_{0,4}^{cl}-48 \tau=0
\end{equation}

First note that \ref{Exact3} can be derived from \ref{Exact2}. Let $H_{(0,4)}^{cl}=d_{\SS^4} \eta$, here $\eta$ is a (0,3)-form. Then $H_{(1,3)}^{cc}=-d_{\HH^7} \eta$ by \ref{Exact3}. Let $f=*_{\SS^4} d_{\SS^4} \eta$. Let $k_{(1,1)}^{cl}=d_{\SS^4} w$, $w$ is (1,0)-form. Put it back to \ref{Exact1}  we get 
\begin{equation}
6 d_{\HH^7} *_{\HH^7} d_{\SS^4} w+d_{\SS^4} *_{\HH^7} (21\sigma-12\tau)+*_{\HH^7} d_{\SS^4} *_{\SS^4} d_{\SS^4} \eta-*_{\SS^4} d_{\HH^7} *_{\HH^7} d_{\HH^7} \eta=0\\
\end{equation}

Apply $*_{\HH^7}$ ($*_{\HH^7}^2=1$), we get 
\[
6*_{\HH^7} d_{\HH^7} *_{\HH^7} d_{\SS^4} w+d_{\SS^4} (21\sigma-12\tau)+d_{\SS^4} *_{\SS^4} d_{\SS^4} \eta-*_{\SS^4} *_{\HH^7} d_{\HH^7} *_{\HH^7} d_{\HH^7} \eta=0
\]
Then let $\eta=*_{\SS^4} d_{\SS^4} \xi$, $\xi$ be a function, and pull out $d_{\SS^4}$
\begin{equation}
-6\delta_{\HH^7} w+(21\sigma-12\tau)-\Lap_{\SS^4}  \xi-\Lap_{\HH^7} \xi=0
\end{equation}
and put the expression to \ref{Exact4},
\[
\Lap_{\SS^4} d_{\SS^4} w+\Lap_{\HH^7} d_{\SS^4} w+12 d_{\SS^4} w+6*_{\SS^4} d_{\HH^7} *_{\SS^4}d_{\SS^4}\xi=0 
\]
Apply $\delta_{\HH^7}$ and pull out $d_{\SS^4}$
\begin{equation}
\Lap_{\SS^4} \delta_{\HH^7} w+\Lap_{\HH^7} \delta_{\HH^7} w +12 \delta_{\HH^7} w +6\Lap_{\HH^7}\xi=0
\end{equation}
Now \ref{Exact5} becomes 
\begin{equation}
\Lap_{\SS^4}\tau+\Lap_{\HH^7} \tau+72\tau+8\Lap_{\SS^4} \xi=0
\end{equation}
And \ref{Exact6} is
\begin{equation}
\Lap_{\SS^4} \sigma+\Lap_{\HH^7} \sigma +12\sigma -4\Lap_{\SS^4} \xi -48\tau=0
\end{equation}

Putting the above four equations together, and suppose the eigenvalue of $\Lap_{\SS^4}$ is $\lambda$, we get
\[
\left( 
\begin{array}{cccc}
12+\lambda+\Lap_{\HH^7} & -48 & -4\lambda& 0\\
0 & 72+\lambda+\Lap_{\HH^7} & 8\lambda & 0\\
21 & -12 & -\lambda-\Lap_{\HH^7} & -6\\
0& 0& 6\Lap_{\HH^7} & 12+\lambda +\Lap_{\HH^7}
\end{array}
\right)
\left(
\begin{array}{c}
\sigma\\
\tau\\
\xi\\
\delta_{\HH^7} w
\end{array}
\right)=0
\]
Removing the off-diagonal $\Delta_{\HH^{7}}$ term in the last line by using the third line, we get
\begin{equation}\label{e:s2}
\left( 
\begin{array}{cccc}
12+\lambda+\Lap_{\HH^7} & -48 & -4\lambda& 0\\
0 & 72+\lambda+\Lap_{\HH^7} & 8\lambda & 0\\
-21 & 12 & \lambda+\Lap_{\HH^7} & 6\\
126& -72& -6\lambda & -24+\lambda +\Lap_{\HH^7}
\end{array}
\right)
\left(
\begin{array}{c}
\sigma\\
\tau\\
\xi\\
\delta_{\HH^7} w
\end{array}
\right)=0
\end{equation}

The determinant, after putting in the indicial operator of $\Lap_{\HH^7}$,  is
\begin{equation}
\begin{aligned}
&\lambda^4 - 4 S^2 \lambda^3 + 24 S*\lambda^3 - 90 \lambda^3 + 6 S^4 \lambda^2 - 72 S^3 \lambda^2 \\&
+  342 S^2 \lambda^2 - 756 S*\lambda^2 + 1152 \lambda^2 - 4 S^6 \lambda + 72 S^5 \lambda - 
  414 S^4 \lambda \\&
  + 648 S^3 \lambda + 1152 S^2 \lambda - 3024 S*\lambda + 10368 \lambda \\&
  + S^8 - 
  24 S^7 + 162 S^6 + 108 S^5 - 6192 S^4 \\& + 31536 S^3 - 33696 S^2 - 
  155520 S = 0
\end{aligned}
\end{equation}
Putting the lowest two eigenvalues for closed 1-form, we get the following two pairs of roots: for 
$\lambda=16
$ the indicial roots are
$\theta_2=3 \pm i \sqrt{21116145}/1655$.
with kernel
$$
\xi_{16}\in  \bigwedge^{cl}_{\lambda=16}(S)
$$
which is the closed 1-form on 4-sphere with eigenvalue 16.
and the other pair is for $\lambda=40$ then
$$\theta_3=3 \pm i3 \sqrt{582842}/20098,$$
with kernel
$$
\xi_{40} \in  \bigwedge^{cl}_{\lambda=40}(S).
$$

\item \textbf{$H_{(3,1)}, H_{(4,0)}$ with closed 1-form / eigenfunctions}
\\
We have
\begin{align}
d_{\SS^4} * H_{(3,1)}^{cl}+d_{\HH^7}*H_{(4,0)}^{cc}+6{}^4V \wedge H_{(4,0)}^{cc}&=0,\\
d_{\HH^7}H_{(3,1)}^{cl}+d_{\SS^4}H_{(4,0)}^{cc}&=0.
\end{align}

%Let $H_{(3,1)}^{cl}=d_{\SS^4}\eta$ where $\eta$ is (3,0), put into second equation to get $H_{(4,0)}^{cc}=-d_{\HH^7}\eta$. Put everything back to first equation, we get
%$d_{\SS^4} *d_{\SS^4} \eta-d_{\HH^7} *d_{\HH^7} \eta-6{}^4V\wedge d_{\HH^7}\eta =0$. Apply $*_{\HH^7} *_{\SS^4}$ get $\Lap_{\SS^4}\eta-\Lap_{\HH^7}\eta+6*_{\HH^7}d_{\HH^7}\eta=0$. Let $\Lap_{\SS^4} \eta=\lambda \eta$ and suppose the half laplacian $*_{\HH^7}d_{\HH^7} \eta=t\eta$, $\Lap_{\HH^7} \eta=t^2 \eta$, then $\lambda-t^2+6t=0$.  The roots of t is $3 \pm \sqrt{9+\lambda}$. When $\lambda$ large, t is either negative or large positive.%

Let 
\[H_{(3,1)}^{cl}=d_{\SS^4}\eta\] 
where $\eta$ is (3,0), put into second equation to get 
\[H_{(4,0)}^{cc}=-d_{\HH^7}\eta\]
Put everything back to first equation, we get
\[d_{\SS^4} *d_{\SS^4} \eta-d_{\HH^7} *d_{\HH^7} \eta-6{}^4V\wedge d_{\HH^7}\eta =0.\] 
Apply $*_{\SS^4}$, and note $*_{\SS^4}^2=(-1)^{k(4-k)}=1$, $\delta_{\SS^4}=(-1)^{4(k+1)+1}*_{\SS^4}d_{\SS^4}*_{\SS^4}=-*_{\SS^4}d_{\SS^4}*_{\SS^4}$, $\Lap_{\SS^4}=d\delta+\delta d$,
\[
*_{\HH^7} (-\delta_{\SS^4})d_{\SS^4}\eta-d_{\HH^7}*_{\HH^7}d_{\HH^7}\eta +d_{\HH^7}\eta=0
\] 
Then apply $*_{\HH^7}$, note $(*_{\HH^7})^2=1$, get
$$-\Lap_{\SS^4}\eta-*_{\HH^7}d_{\HH^7}*_{\HH^7}d_{\HH^7}\eta+6*_{\HH^7}d_{\HH^7}\eta=0.$$
Let $\Lap_{\SS^4} \eta=\lambda \eta$ 
\[
-\lambda \eta-\Lap_{\HH^7}\eta+6*_{\HH^7}d_{\HH^7}\eta=0
\]
The indicial equation: using $I[d](s)w=(-1)^k(s-k)w\wedge \frac{dx}{x}$,
\[
-\lambda \eta+(s-3)^2\eta+6(s-3)*_{\SS^6}\eta=0
\]
that is 
\[
(s-3)^2\pm 6i(s-3)-16=0
\]
with roots 
$$s=3\pm \sqrt{7}\pm 3i.$$

\item \textbf{$k_{(1,1)}, H_{(1,3)}, H_{(2,2)}$ with coclosed 1-form / closed 3-form / coclosed 2-form}

We have
\begin{equation}\label{coexact1}
6d_{\HH^7} *_{\HH^7}k_{(1,1)}^{cc}+d_{\HH^7}*H_{(1,3)}^{cl}=0
\end{equation}
\begin{equation}\label{coexact2}
d_{\SS^4} *H_{(1,3)}^{cl}+d_{\HH^7}*H_{(2,2)}^{cc}+6d_{\SS^4} *_{\HH^7} k_{(1,1)}^{cc}=0
\end{equation}
\begin{equation}\label{coexact3}
d_{\HH^7} H_{(1,3)}^{cl}+d_{\SS^4}H_{(2,2)}^{cc}=0
\end{equation}
\begin{equation}\label{coexact4}
\frac{1}{2}\Delta_{\SS^4} k_{(1,1)}^{cc}+\frac{1}{2}\Delta_{\HH^7}
k_{(1,1)}^{cc}+6k_{(1,1)}^{cc}-\frac{1}{2}*_{\SS^4}H_{(1,3)}^{cl}=0
\end{equation}
First note that (\ref{coexact1}) can be derived from (\ref{coexact2})
Let $H_{(1,3)}^{cl}=d_{\SS^4}\eta$, where $\eta$ is (1,2)-form. 
Then $H_{2,2}^{cc}=-d_{\HH^7}\eta$ from (\ref{coexact3}). Put it to  (\ref{coexact2}),  $d_{\SS^4} * d_{\SS^4} \eta-d_{\HH^7} *d_{\HH^7} \eta+6d_{\SS^4} *_{\HH^7} k_{1,1}^{cc}=0$. Apply $*_{\SS^4}, *_{\HH^7}$, get $-\Lap_{\SS^4} \eta-\Lap_{\HH^7} \eta +6 *_{\SS^4} d_{\SS^4} k_{(1,1)}^{cc}=0$. Apply $*_{\SS^4} d_{\SS^4}$ again, get $-\Lap_{\SS^4} (*_{\SS^4}d_{\SS^4} \eta)-\Lap_{\HH^7} (*_{\SS^4}d_{\SS^4} \eta)-6\Lap_{\SS^4} k_{(1,1)}^{cc}=0$. Combining with (\ref{coexact4}), and let $\lambda$ be the eigenvalue for $\Lap_{\SS^4}$ on coclosed 1-form,  we get 
\[
\left(\begin{array}{cc}
-\lambda -\Lap_{\HH^7} & -6\lambda\\
-1 & \lambda+\Lap_{\HH^7} +12
\end{array} \right) 
\left(
\begin{array}{c}
*_{\SS^4} d_{\SS^4} \eta\\
k_{(1,1)}^{cc}
\end{array}
\right)=0
\]
The indicial equation is
$$
\lambda^2-(36+(s-1)(s-5)+s^2-6s-1)\lambda -(s-1)(s-5)(-s^2+6s+1)=0.
$$
With smallest eigenvalue for coclosed 1-form to be $\lambda=24$, indicial roots are
$$
s=3\pm \sqrt{\pm 3\sqrt{97}+31}
$$

\item \textbf{$H_{(2,2)}, H_{(3,1)}$ with closed 2-form / coclosed 1-form}
\begin{eqnarray}
d_{\SS^4}*H_{(2,2)}^{cl}+d_{\HH^7}*H_{(3,1)}^{cc}=0\\
d_{\HH^7}H_{(2,2)}^{cl}+d_{\SS^4} H_{(3,1)}^{cc}=0
\end{eqnarray}

Apply $d_{\HH^7}$ and $d_{\SS^4}$ to the equations, we have
\begin{equation}
d_{\HH^7} d_{\SS^4}*H_{(2,2)}^{cl}=0, d_{\SS^4} d_{\HH^7} H_{(2,2)}^{cl}=0
\end{equation}
let $H_{(2,2)}^{cl}=d_{\SS^4}\eta$ where $\eta$ is a coclosed (2,1)-form, Putting it back, and using $d_{\SS^4}$ is an isomorphism, $d_{\HH^7} \eta=-H_{(3,1)}^{cc}$. Then from first equation, $d_{\SS^4}*d_{\SS^4} \eta-d_{\HH^7}*d_{\HH^7}\eta=0$, which is  $-*_{\HH^7}*_{\SS^4}\Lap_{\SS^4}\eta-*_{\SS^4} *_{\HH^7}\Lap_{\HH^7}\eta=0$ then it requires $\Lap_{\HH^7} \eta=-\lambda \eta$. Putting $\lambda=4(k+2)(k+3)$, the result is
$$
s=3\pm \sqrt{17}.
$$
\end{enumerate}

\newcommand{\etalchar}[1]{$^{#1}$}


\begin{thebibliography}{CDF{\etalchar{+}}84}

\bibitem[BFOP02]{blau2002penrose}
Matthias Blau, Jose Figueroa-O'Farrill, and George Papadopoulos.
\newblock Penrose limits, supergravity and brane dynamics.
\newblock {\em Classical and Quantum Gravity}, 19(18):4753, 2002.

\bibitem[BST87]{bergshoeff1987supermembranes}
E~Bergshoeff, Ergin Sezgin, and Paul~K Townsend.
\newblock Supermembranes and eleven-dimensional supergravity.
\newblock {\em Physics Letters B}, 189(1):75--78, 1987.

\bibitem[Bor07]{Borthwick}
David Borthwick. 
\newblock Spectral theory of infinite-area hyperbolic surfaces. 
\newblock{\em Basel: Birkhäuser,} 2007.


\bibitem[CDF{\etalchar{+}}84]{castellani1984bosonic}
L~Castellani, Riccardo D'Auria, P~Fre, K~Pilch, and P~Van~Nieuwenhuizen.
\newblock The bosonic mass formula for freund-rubin solutions of d= 11
  supergravity on general coset manifolds.
\newblock {\em Classical and Quantum Gravity}, 1(4):339, 1984.

\bibitem[CS77]{MR0441129}
E.~Cremmer and J.~Scherk.
\newblock Spontaneous compactification of extra space dimensions.
\newblock {\em Nuclear Phys. B}, 118(1--2):61--75, 1977.

\bibitem[Del02]{delay2002essential}
Erwann Delay.
\newblock Essential spectrum of the lichnerowicz laplacian on two tensors on
  asymptotically hyperbolic manifolds.
\newblock {\em Journal of Geometry and Physics}, 43(1):33--44, 2002.

\bibitem[DZ18]{DZ}
Semyon Dyatlov and Maciej Zworski. 
\newblock Mathematical theory of scattering resonances. 
\newblock book in progress.

\bibitem[FR80]{freund1980dynamics}
Peter~GO Freund and Mark~A Rubin.
\newblock Dynamics of dimensional reduction.
\newblock {\em Physics Letters B}, 97(2):233--235, 1980.

\bibitem[GL91]{MR1112625}
C.~Robin Graham and John~M. Lee.
\newblock Einstein metrics with prescribed conformal infinity on the ball.
\newblock {\em Adv. Math.}, 87(2):186--225, 1991.


\bibitem[GZ03]{MR1965361}
C.~Robin Graham and Maciej Zworski.
\newblock Scattering matrix in conformal geometry.
\newblock {\em Invent. Math.}, 152(1):89--118, 2003.



\bibitem[GN06]{guillarmou2006wave}
Colin Guillarmou and Fr{\'e}deric Naud.
\newblock Wave 0-trace and length spectrum on convex co-compact hyperbolic
  manifolds.
\newblock {\em Communications in analysis and geometry}, 14(5):945--967, 2006.

\bibitem[Gui05]{MR2153454}
Colin Guillarmou.
\newblock Meromorphic properties of the resolvent on asymptotically hyperbolic
  manifolds.
\newblock {\em Duke Math. J.}, 129(1):1--37, 2005.

\bibitem[Gui06]{MR2276069}
Colin Guillarmou.
\newblock Scattering and resolvent on geometrically finite hyperbolic manifolds
  with rational cusps.
\newblock In {\em Seminaire: {E}quations aux {D}\'eriv\'ees {P}artielles.
  2005--2006}, S\'emin. \'Equ. D\'eriv. Partielles, pages Exp. No. III, 17.
  \'Ecole Polytech., Palaiseau, 2006.

\bibitem[Kan09]{kantor2009eleven}
Joshua~M Kantor.
\newblock {\em Eleven Dimensional Supergravity on Edge Manifolds}.
\newblock PhD thesis, University of Washington, 2009.

\bibitem[Lan85]{MR772023}
Serge Lang.
\newblock {\em Differential manifolds}.
\newblock Springer-Verlag, New York, second edition, 1985.

\bibitem[Lan99]{MR1666820}
Serge Lang.
\newblock {\em Fundamentals of differential geometry}, volume 191 of {\em
  Graduate Texts in Mathematics}.
\newblock Springer-Verlag, New York, 1999.

\bibitem[Lau03]{Lauter}
Robert Lauter. 
\newblock{\em Pseudodifferential analysis on conformally compact spaces.} No. 777.
\newblock American Mathematical Soc., 2003.

\bibitem[Lee06]{lee2006fredholm}
John~M Lee.
\newblock {\em Fredholm operators and Einstein metrics on conformally compact
  manifolds}.
\newblock American Mathematical Soc., 2006.

\bibitem[Maz88]{MR961517}
Rafe~R. Mazzeo.
\newblock The {H}odge cohomology of a conformally compact metric.
\newblock {\em J. Differential Geom.}, 28(2):309--339, 1988.

\bibitem[Maz91]{MR1133743}
Rafe~R. Mazzeo.
\newblock Elliptic theory of differential edge operators. {I}.
\newblock {\em Comm. Partial Differential Equations}, 16(10):1615--1664, 1991.

\bibitem[Maz91II]{MazUni}
Rafe~R. Mazzeo.
\newblock Unique continuation at infinity and embedded eigenvalues for asymptotically hyperbolic manifolds.
\newblock {\em American Journal of Mathematics} 113, no. 1: 25-45, 1991.

\bibitem[Mel93]{melrose1993atiyah}
Richard~B. Melrose.
\newblock {\em The {A}tiyah-{P}atodi-{S}inger index theorem}.
\newblock Wellesley: AK Peters, 1993.

\bibitem[Mel95]{melrose1995geometric}
Richard~B. Melrose.
\newblock {\em Geometric scattering theory}.
\newblock Cambridge University Press, 1995.

\bibitem[MM87]{MR916753}
Rafe~R. Mazzeo and Richard~B. Melrose.
\newblock Meromorphic extension of the resolvent on complete spaces with
  asymptotically constant negative curvature.
\newblock {\em J. Funct. Anal.}, 75(2):260--310, 1987.

\bibitem[MP90]{mazzeo1990}
Rafe~R. Mazzeo and Ralph~S. Phillips.
\newblock Hodge theory on hyperbolic manifolds.
\newblock {\em Duke Math. J.}, 60(2):509--559, 04 1990.

\bibitem[Nah78]{nahm1978supersymmetries}
Werner Nahm.
\newblock Supersymmetries and their representations.
\newblock {\em Nuclear Physics B}, 135(1):149--166, 1978.

\bibitem[Nas11]{Nastase:2011aa}
Horatiu Nastase.
\newblock Introduction to supergravity.
\newblock {\em arXiv preprint arXiv:1112.3502}, 2011.

\bibitem[Tay11]{MR2744150}
Michael~E. Taylor.
\newblock {\em Partial differential equations {I}. {B}asic theory}, volume 115
  of {\em Applied Mathematical Sciences}.
\newblock Springer, New York, second edition, 2011.

\bibitem[VN81]{van1981supergravity}
Peter Van~Nieuwenhuizen.
\newblock Supergravity.
\newblock {\em Physics Reports}, 68(4):189--398, 1981.

\bibitem[VN85]{van1985complete}
Peter Van~Nieuwenhuizen.
\newblock The complete mass spectrum of d= 11 supergravity compactified on $\SS^{4}$
  and a general mass formula for arbitrary cosets $M^{4}$.
\newblock {\em Classical and Quantum Gravity}, 2(1):1, 1985.

\bibitem[Vas12]{Vasy}
Andr\'as Vasy. 
\newblock Analytic continuation and high energy estimates for the resolvent of the Laplacian on forms on asymptotically hyperbolic spaces. \newblock arXiv preprint arXiv:1206.5454 (2012).


\bibitem[Wit97]{Witten:1996md}
Edward Witten.
\newblock {On flux quantization in M theory and the effective action}.
\newblock {\em J.Geom.Phys.}, 22:1--13, 1997.

\end{thebibliography}
\end{document}